\documentclass[twoside,11pt]{article}

\usepackage{amsmath}
\usepackage{algorithm} 
\usepackage{algpseudocode} 
\usepackage{booktabs}
\usepackage{adjustbox}
\usepackage[colorlinks=false]{hyperref}
\usepackage{tabularx,array,makecell,threeparttable,ragged2e}
\usepackage{caption}
\usepackage[nohyperref, preprint]{jmlr2e}

\newcommand{\RateCell}[1]{\makecell[r]{$\widetilde{\mathcal O}(#1)$}}
\newcommand{\Dash}{\makecell[c]{--}}
\newtheorem{assumption}[theorem]{Assumption}

\usepackage{lastpage}
\jmlrheading{27}{2026}{1-\pageref{LastPage}}{5/25; Revised 4/26}{4/26}{25-1016}{Axel Friedrich Wolter and Tobias Sutter}

\ShortHeadings{Two-Timescale Primal-Dual RL}{Wolter and Sutter}
\firstpageno{1}

\begin{document}

\title{A Two-Timescale Primal-Dual Framework for Reinforcement Learning via Online Dual Variable Guidance}

\author{\name Axel Friedrich Wolter \email axel-friedrich.wolter@uni-konstanz.de\\
        \addr Department of Computer Science\\
              University of Konstanz, Germany\\
              Universitätsstrasse 10, 78464 Konstanz\vspace{2mm} \\ 
\name Tobias Sutter \email tobias.sutter@unisg.ch\\
        \addr Department of Economics\\
              University of St.Gallen, Switzerland\\  
              Rosenbergstrasse 20-22, 9000 St. Gallen
              }
           
\editor{Nan Jiang}
\maketitle

\begin{abstract}%
We study reinforcement learning by combining recent advances in regularized linear programming formulations with the classical theory of stochastic approximation. 
Motivated by the challenge of designing algorithms that leverage off-policy data while maintaining on-policy exploration, we propose PGDA-RL, a novel primal-dual projected gradient descent-ascent algorithm for solving regularized Markov decision processes (MDPs). PGDA-RL integrates experience replay-based gradient estimation with a two-timescale decomposition of the underlying nested optimization problem. 
The algorithm operates asynchronously, interacts with the environment through a single trajectory of correlated data, and updates its policy online in response to the dual variable associated with the occupancy measure of the underlying MDP. We prove that PGDA-RL converges almost surely to the optimal value function and policy of the regularized MDP. Our convergence analysis relies on tools from stochastic approximation theory and holds under weaker assumptions than those required by existing primal-dual RL approaches, notably removing the need for a simulator or a fixed behavioral policy. 
Under a strengthened ergodicity assumption on the underlying Markov chain, we establish a last-iterate finite-time guarantee with $\widetilde{\mathcal O}(k^{-2/3})$ mean-square convergence, aligning with the best-known rates for two-timescale stochastic approximation methods under Markovian sampling and biased gradient estimates.

\end{abstract}


\section{Introduction}
Recent years have seen many new promising developments in the theory and application of reinforcement learning (RL) as a framework for sequential decision making with notable successes of regularization-based algorithms, for example, the Soft Actor-Critic \citep{ref:haarnoja2018:soft} or Trust Region Policy Optimization \citep{ref:schulman2015:trust}. The regularization typically takes the form of a convex penalty term in the policy, which is a setting that has been subsequently formalized by \cite{geist2019theory,ref:neu-17}. The adjustment of the objective aims at inducing desirable exploratory behavior of the learned policies and stability in their training \citep{ref:schulman2017:equivalence}, and allows the application of dynamical systems analysis that takes advantage of the adjusted problem's strongly convex structure \citep{li2024accelerating}. 

As another development, the linear programming (LP) reformulation of Markov decision processes (MDPs), originally dating back to \cite{ref:Manne-60,ref:Hernandez-96,ref:Borkar-2002}, has recently received renewed interest. It has led to LP-based algorithms shown to be on par with broadly applied benchmarks, for example, the Dual Actor-Critic \citep{ref:Dai-18:boosting}. 
The LP approach to RL allows the reformulation of the Bellman optimality equation, which characterizes an optimal policy, as a min-max saddle-point problem that can be addressed with gradient methods. The analysis of the various algorithms proposed in this line of research often relies on simplifying assumptions such as access to a sampling generator \citep{chen2016}, a known model \citep{li2024accelerating}, or the uniqueness of the optimal policy \citep{ref:lee2019:stochastic}. A further challenge in solving the min-max RL problem formulation lies in the nested loop structure that arises when the inner optimization is first solved approximately during each outer iteration, which complicates the convergence analysis \citep{ref:gabbianelli2024:offline,ref:Dai-18:boosting}.

The asymptotic convergence under general model-access assumptions is well understood for (approximate) dynamic programming algorithms; (see \cite{ref:tsitsiklis1994:asynchronous, ref:jaakkola1993:convergence}). However, finite-time guarantees for the \emph{last iterate} remain scarce for LP-based primal-dual schemes driven by Markovian single-trajectory data and iterate-dependent (on-policy) exploration, particularly when gradients are formed from replay buffers and are therefore biased at finite times. This perspective is complementary to regret-minimization in online RL: rather than optimizing regret or minimax sample complexity via exploration bonuses, we study the stability and last-iterate behavior of a primalreplay-buffer biasdual stochastic approximation method for an LP saddle-point formulation in a realistically correlated-data regime.

In this work, we combine recent advances on regularized LP formulations of RL \citep{li2024accelerating} with the classical theory of stochastic approximation \citep{borkar2023stochastic} and modern finite-time analyses for two-timescale stochastic approximation \citep{ref:zeng2024:two}. Our goal is to design a provably convergent algorithm that maintains on-policy exploration while performing off-policy updates using previously observed transitions stored in an experience replay buffer, and to quantify its last-iterate finite-time rate. The main contributions of this paper are summarized as follows:

\begin{enumerate}
    \item We introduce a novel primal-dual Projected Gradient Descent-Ascent (PGDA-RL) algorithm for regularized MDPs. The algorithm combines experience replay-based gradient estimation with a two-timescale decomposition of the nested optimization problem and leverages on-policy exploration. PGDA-RL operates asynchronously and interacts with the environment through a single trajectory of data, generated under a policy that evolves in response to the dual variable associated with the occupancy measure.
    \item We establish that PGDA-RL converges almost surely to the optimal value function and policy of the regularized MDP. Our convergence analysis is based on tools from the stochastic approximation literature and requires significantly weaker assumptions than existing primal-dual RL methods. Notably, we do not assume access to a simulator or a fixed behavior policy.
    \item Under a stronger ergodicity assumption on the underlying Markov chain, we establish a \emph{last-iterate} finite-time guarantee for the asynchronous PGDA-RL scheme despite replay-buffer biased gradient estimates. Concretely, on a high-probability visitation event $\mathcal G_\delta$ ensuring linear growth of the least-visited state-action count after a burn-in time, the dual iterate satisfies $\mathbb{E}[\|\rho_k-\rho^\star\|_2^2\mid\mathcal{G}_\delta]=\widetilde{\mathcal{O}}(k^{-2/3})$. This yields corresponding $\widetilde{\mathcal{O}}(k^{-2/3})$ mean-square convergence guarantees for the regularized optimal value and policy. 
\end{enumerate}

\subsection{Related Work}
Recent research has explored several directions in reinforcement learning that build upon primal-dual methods and the regularized Markov decision processes framework.
Table~\ref{table:related:work} provides a concise overview of the most closely related LP-based primal-dual RL methods and highlights the key differences compared to our approach regarding underlying model access, function approximation setting, and convergence analysis.

\cite{ref:Dai-18:boosting} introduced a dual actor-critic algorithm rooted in the LP formulation that incorporates multi-step Bellman equations, path regularization, and a nested stochastic dual-ascent update. While their approach shares similarities with ours in starting from an LP/Lagrangian viewpoint and employing convexity-inducing regularizers, we adopt a two-timescale stochastic approximation framework that yields an \emph{incremental} update rule and enables an almost sure convergence analysis under Markovian sampling. 
Similarly, works by \cite{chen2016,ref:wang2017:primal,ref:chen2018:scalable} study LP formulations for discounted, average reward, and finite-horizon MDPs through projected stochastic primal-dual methods (including incremental coordinate updates). These works provide finite-sample or asymptotic convergence results under the crucial assumption of access to a sampling generator. In contrast, we focus on the \emph{last iterate} and establish \emph{almost sure convergence} assuming only access to a single trajectory of state-action pairs (Markovian model access).
\cite{ref:lee2019:stochastic} extends \cite{chen2016} by incorporating a Q-function estimation step. Their analysis allows time-varying behavioral policies under conditions requiring sublinear convergence of the induced state-action distribution, which can be difficult to verify a priori and makes fully on-policy analysis nontrivial. Moreover, their algorithm relies on samples drawn from the stationary state-action distribution of the current policy.
In~\cite{ref:bas-serrano2021:logistic} and the subsequent thesis~\cite{ref:BasSerrano2022:LagragianDuality}, the authors introduce the logistic Bellman error derived from an entropy-regularized LP formulation and establish finite-sample rates in the tabular setting under stationary-distribution generator access.
In comparison, our stochastic approximation analysis further allows us to avoid stationary distribution sampling and to analyze iterate-dependent (on-policy) exploration under Markovian single-trajectory access. Under a strengthened ergodicity condition on the Markov chain induced by the iterate-dependent behavioral policy, we additionally establish a \emph{last-iterate} finite-time mean-square convergence rate.
\begin{table}[t]
\centering
\refstepcounter{table}
\label{table:related:work}
\scriptsize
\setlength{\tabcolsep}{3.2pt}
\renewcommand{\arraystretch}{1.10}
\begin{threeparttable}
\begin{tabularx}{\linewidth}{@{}
>{\RaggedRight\arraybackslash}p{0.22\linewidth}
>{\RaggedRight\arraybackslash}p{0.16\linewidth}
>{\RaggedRight\arraybackslash}p{0.12\linewidth}
>{\RaggedRight\arraybackslash}p{0.16\linewidth}
>{\RaggedRight\arraybackslash}X
>{\RaggedLeft\arraybackslash}p{0.10\linewidth}
@{}}
\toprule
\textbf{Paper} & \textbf{Model Access} & \textbf{Func.\ Appr.} & \textbf{Analysis} & \textbf{Guarantee}& \textbf{Rate} \\
\midrule
\cite{chen2016} &
Generator &
Tabular &
PAC &
Duality gap &
\RateCell{\epsilon^{-2}} \\

\cite{ref:Dai-18:boosting} &
Markovian &
Non-linear &
Sketch &
\Dash &
\Dash \\

\cite{ref:lee2019:stochastic} &
\makecell[l]{Stationary\\Distr.\ Generator}&
Tabular &
PAC &
Duality gap &
\RateCell{\epsilon^{-2}} \\

\cite{ref:BasSerrano2022:LagragianDuality} &
\makecell[l]{Stationary\\Distr.\ Generator} &
Tabular &
PAC &
\makecell[l]{$\epsilon$-optimal\\policy} &
\RateCell{\epsilon^{-9}} \\

\cite{ref:gabbianelli2024:offline} &
Offline &
Linear &
PAC &
\makecell[l]{$\epsilon$-optimal\\policy} &
\RateCell{\epsilon^{-2}} \\

\cite{li2024accelerating} &
Full model &
Tabular &
\makecell[l]{Deterministic\\asymptotic} &
\Dash &
\Dash \\

\textbf{Our paper} &
Markovian &
Tabular &
\makecell[l]{Almost sure,\\Convergence rate} &
\makecell[l]{Dual iterate MSE \\ after burn-in} &
\RateCell{\epsilon^{-3/2}} \\
\bottomrule
\end{tabularx}
\caption{Overview of LP-based RL methods and their guarantees.
We summarize related approaches based on linear-programming or Lagrangian formulations of MDPs, organized by the model access assumption, the function-approximation regime, the type of analysis, and the corresponding guarantee and reported rate.
The ``Rate'' column is stated with respect to the corresponding guarantee, and is therefore not directly comparable across rows when the guarantees or access models differ.
}
\end{threeparttable}
\end{table}

Two-timescale stochastic approximation has been applied to dynamic programming-based actor-critic algorithms \citep{ref:konda1999:actor} and to constrained MDPs \citep{ref:borkar2005:actor}, where safety constraints are added to the reward maximization objective. The Lagrangian relaxation of the constrained optimization problem motivates the application of various primal-dual methods \citep{ref:chen2024:primal-dual,ref:hong2024:primal,ref:li2023:double}.

The LP approach has also been applied to Offline RL problems, where near-optimal policies are learned from static data sets of transitions \citep{ref:zhan2022:offline,ref:ozdaglar-23a,ref:gabbianelli2024:offline}. Even though the analysis objective differs, since we consider the online setting, their saddle-point problem resulting from the Lagrangian dual formulation of the LP is closely related to our optimization target. 

Parallel investigations into regularized MDPs have also been influential. \cite{li2024accelerating} explore a primal-dual formulation of entropy-regularized MDPs with additional convexity in the primal variable, providing asymptotic convergence analysis using natural gradient methods, albeit under the assumption of full model knowledge. \cite{ref:neu-17} and \cite{geist2019theory} formalize the entropy-regularized MDP framework, outlining the conditional entropy regularizer and the contraction properties of the regularized Bellman evaluation operator. \cite{ying2020note} offers a comprehensive overview of various LP formulations in both regularized and unregularized contexts, emphasizing their equivalence to Bellman equations and policy gradient methods. Furthermore, convex Q-learning \citep{meyn2022control,ref:Prashant-20,ref:lu2023convex} is motivated by the idea of deriving a convex relaxation of the squared Bellman error based on the LP approach. The method involves solving a constrained convex optimization problem over the primal variable, the Q-function. It is primarily designed for analyzing the non-tabular setting, where the value function is learned via function approximation, which distinguishes this approach from our own.

Finally, our finite-time analysis builds on recent last-iterate convergence-rate results for two-timescale stochastic approximation \cite{ref:zeng2024:two}. We face additional challenges in our setting, namely asynchronous updates and replay-buffer gradient estimates that introduce bias. Despite this, on the high-probability good visitation event $\mathcal G_\delta$ (after a burn-in period), we prove a last-iterate mean-square rate $\widetilde{\mathcal O}(k^{-2/3})$ for the dual iterates, implying corresponding mean-square rates for the induced regularized value function and policy. Equivalently, achieving $\epsilon$ accuracy in mean-square error requires $k(\epsilon)=\widetilde{\mathcal O}(\epsilon^{-3/2})$ iterations to ensure $\mathbb E[\|\rho_k-\rho^\star\|_2^2\mid\mathcal G_\delta]\le \epsilon$ (up to logarithmic factors).

\textit{Structure.} 
The structure of the paper is as follows. Section~\ref{sec:problem:formulation} formally introduces regularized MDPs and the associated saddle-point formulation that underlies our proposed algorithms. Section~\ref{sec:synchronous:setting} presents a first algorithm in a synchronous setting, assuming access to a generative model. In Section~\ref{sec:asynchronous:setting}, we extend this setting to the asynchronous case, where learning occurs from a single trajectory of data generated under a potentially changing policy. The section further contains the almost sure convergence result { and our convergence rate analysis.}  Section~\ref{sec:numerical:results} demonstrates the effectiveness of the proposed asynchronous algorithm on a standard reinforcement learning benchmark.  
Finally, Section~\ref{sec:conclusion} concludes the paper. The proofs details and technical results are deferred to the appendix to improve the readability of the main text.

\textit{Notation.} For any state $s\in\mathcal{S}$, we denote by $\delta_s$ the Dirac measure in $s$.
For a variable $v(s,a)$ depending on state-action pairs $(s,a)\in\mathcal{S}\times\mathcal{A}$, we denote the state marginal by $\tilde v(s):=\sum_{a\in\mathcal A}v(s,a)$.
For a finite set $\mathcal{S}$ with cardinality $n$, the probability simplex over $\mathcal{S}$ is denoted by $\Delta_\mathcal{S}=\{x \in \mathbb{R}_{+}^n: \sum_{(s,a)\in\mathcal S\times\mathcal A} x_i=1\}$. 
We use the symbol $e_i\in\mathbb{R}^{|\mathcal{S}|}$ for the $i$-th standard basis vector on $\mathbb{R}^{|\mathcal{S}|}$.
For any logical expression $\mathcal{E}$, the indicator function $\mathbf{1}\{\mathcal{E}\}$ evaluates to $1$ if $\mathcal{E}$ is true and to 0 otherwise. 
To denote sequences, we use $\{\beta_k\}$ and $\{\beta(k)\}$ interchangeably for $k\in\mathbb N$. 
For $a,b\in \mathbb{R}$ we denote with $a\lor b$ the larger of $a$ and $b$ and with $\lfloor a\rfloor$ the integer part of $a$. 
For a compact set $H\subset\mathbb{R}^d$, $d\in\mathbb N$, we denote the Euclidean projection operator as $\Pi_H$, where $\Pi_H:\mathbb{R}^d\rightarrow H$ is defined by $\Pi_H(x)=\arg\min_{y\in H}\|x-y\|_2$. 
We denote the total variation norm as $\|\cdot\|_{\mathrm{TV}}$. 
For kernels $\mathcal P,\,\mathcal Q$ on $\mathcal X$, we define \(\|\mathcal P-\mathcal Q\|_{\infty}\;:=\;\sup_{x\in\mathcal X}\,\|\mathcal P(\cdot|x)-\mathcal Q(\cdot|x)\|_{\mathrm{TV}}.\)

\section{Problem Formulation} \label{sec:problem:formulation}
We consider a Markov decision process $\mathcal{M} = (\mathcal{S}, \mathcal{A}, \mathcal{P}, r, \gamma)$, where $\mathcal{S}$ and $\mathcal{A}$ denote finite state and action spaces, respectively. The transition kernel $\mathcal{P}(\cdot | s, a) \in \Delta_{\mathcal{S}}$ specifies the distribution over next states given the current state $s \in \mathcal{S}$ and action $a \in \mathcal{A}$. Each state-action pair $(s, a)$ is associated with a reward $r(s, a) \in \mathbb{R}$, and $\gamma \in (0, 1)$ denotes the discount factor.
Our objective is to characterize and subsequently compute optimal randomized stationary Markovian policies $\pi : \mathcal{S} \rightarrow \Delta_{\mathcal{A}}$ that maximize the expected total discounted reward. We denote by $\Pi = (\Delta_{\mathcal{A}})^{|\mathcal{S}|}$ the set of all stationary Markovian policies.
We assume a positive and bounded deterministic reward, that is, $r(s,a)\in[0,C_r]$ for all $(s,a)\in\mathcal{S}\times\mathcal{A}$ and some positive constant $C_r<\infty$. 
The unregularized value function of a stationary Markovian policy $\pi\in\Pi$ is defined as
$$V^\pi_{ur}(s)=\mathbb{E}^\pi_{s}\left[\sum_{k=0}^\infty\gamma^k r(s_k,a_k)\right],$$ where the notation $\mathbb{E}^\pi_{\bar s}$ refers to the expectation under $s_0=s, a_k\sim\pi(\cdot|s_k),$ $s_{k+1}\sim \mathcal{P}(\cdot| s_k,a_k)$. 
The optimal value of the unregularized MDP is $V^\star_{ur}:=\max_{\pi\in\Pi} V^\pi_{ur}$, the maximum value of all stationary Markovian policies. 
The Bellman optimality equation states that the optimal value function satisfies $$V^\star_{ur}(s)=\max_{\pi\in\Pi}\mathbb{E}^\pi_{s} \left[r( s,a)+\gamma\sum_{s'\in\mathcal{S}} V^\pi_{ur}\left(s^{\prime}\right) \mathcal{P}(s'|s,a)\right], \qquad \forall s\in\mathcal{S},$$ 
which admits a unique solution. 
We approach the unregularized MDP by its formulation as a linear program, \citet[Section 6.9]{puterman1994markov}. 
Reformulating the Bellman optimality equations as inequalities yields that a bounded function $V:\mathcal{S}\rightarrow\mathbb{R}$ that satisfies 
\begin{align*}
V(s)\geq r(s,a)+\gamma \sum_{s'\in\mathcal{S}} V\left(s^{\prime}\right) \mathcal{P}(s'|s,a) \quad \forall (s,a)\in\mathcal{S}\times\mathcal{A},
\end{align*}
is an upper bound to the value $V^\star_{ur}$ of the MDP. 

Building on this finding, the primal formulation of the Bellman LP is then given by
\begin{align}
    \label{LP-p}
    P: \begin{cases}\min _{V \in \mathbb{R}^{|\mathcal{S}|}} & \mu^\top V \\ \text {s.t.} & 0 \geq -V(s)+r(s,a)+\gamma \sum_{s'\in\mathcal{S}} V\left(s^{\prime}\right) \mathcal{P}(s'|s,a)\quad \forall (s,a)\in\mathcal{S}\times\mathcal{A},\end{cases}
\end{align}
where $\mu\in\mathbb{R}_+^{|\mathcal{S}|}$ is any strictly positive vector such that $\sum_{s\in\mathcal{S}}\mu(s)=1$.
It is well-known \cite[Section~6.9]{puterman1994markov} that the optimizer to the primal LP~\eqref{LP-p} is exactly the optimal value function $V^\star_{ur}$.
The corresponding dual linear problem is defined as
\begin{align}
    \label{LP-d}
    D: \begin{cases}    \max _{\rho\geq0} &\sum_{(s,a)\in\mathcal{S}\times\mathcal{A}} \rho(s, a) r(s, a) \\ \text {s.t.} & \sum_{a\in\mathcal{A}}\rho(s, a)-\gamma \sum_{s^{\prime}\in\mathcal{S}, a^{\prime}\in\mathcal{A}}\rho\left(s^{\prime}, a^{\prime}\right)\mathcal{P}\left(s | s^{\prime}, a^{\prime}\right)=\mu(s) \quad \forall s\in\mathcal{S},\end{cases}
\end{align}
where $\mu$ is the vector specified in the primal Bellman LP~\eqref{LP-p}.
A dual feasible variable $\rho$ corresponds to the \textit{discounted occupancy measure} of a stationary randomized policy $\pi_\rho$ defined for $s\in\mathcal{S}$ as 
\begin{align}
\label{eq:dual:policy}
\pi_\rho(a | s)=\frac{\rho(s, a)}{\tilde\rho(s)}, \qquad \forall a\in\mathcal{A},
\end{align} 
where $\tilde\rho(s):=\sum_{a'\in\mathcal{A}} \rho(s, a')$ denotes the state marginal of the dual variable.
Let $\rho^\star_{ur}$ denote a solution to the dual problem~\eqref{LP-d}, then $\pi_{\rho_{ur}^\star}$ is an optimal policy, that is, $V^{\pi_{\rho^\star_{ur}}}_{ur}=V^\star_{ur}$, see \cite[Theorem 6.9.4]{puterman1994markov}.
The Lagrangian formulation of the primal LP~\eqref{LP-p} is
\begin{align}
    \label{eq: clLP}
    \max_{\rho\geq0} \min_{V\in \mathbb{R}^{|\mathcal{S}|}} L_{ur}(V, \rho):= \mu^\top V+\sum_{(s,a)\in\mathcal{S}\times\mathcal{A}} \rho(s, a) \Delta[V](s, a),
\end{align}
where $$\Delta[V](s, a)=-V(s)+r(s, a)+\gamma \sum_{s'\in\mathcal{S}} V\left(s^{\prime}\right) \mathcal{P}(s'|s,a)$$ denotes the Bellman error. It can be directly seen that the primal variable $V$ can be assumed without loss of generality to belong to the compact set $\mathcal{V} := \{v\in\mathbb{R}_+^{|\mathcal{S}|} : v_i \leq \frac{C_r}{1-\gamma} \ \forall i=1,\dots |\mathcal{S}|\}$.

\subsection{Regularized MDP}
Our approach is based on the regularized MDP formulation as formally described in \citet{geist2019theory}. Let $G:\mathbb{R}^{|\mathcal{A}|}_+\rightarrow\mathbb{R}$ be a strongly concave function and $\eta_\rho>0$ a regularization parameter. For a stationary policy $\pi\in\Pi$, the regularized Bellman operator $T_{\pi,G}:\mathbb{R}^{|\mathcal{S}|}\rightarrow\mathbb{R}^{|\mathcal{S}|}, V\mapsto T_{\pi,G}V$ is defined state-wise as 
\begin{align}
\label{eq:regBellman}
    [T_{\pi,G}V](s):=\mathbb{E}^\pi _{\bar s}&\left[r(s,a)+\gamma\sum_{s'\in\mathcal{S}} V\left(s^{\prime}\right) \mathcal{P}(s'|s,a)\right]+\eta_\rho G(\pi(\cdot|s)).
\end{align}
The regularized Bellman optimality operator is defined accordingly as
\begin{align}
    \label{eq:regBellman:opt}
    T_{\star,G}:\mathbb{R}^{|\mathcal{S}|}\rightarrow\mathbb{R}^{|\mathcal{S}|}&, V\mapsto T_{\star,G}V:=\max_{\pi\in\Pi}T_{\pi,G}V,
\end{align}
where the maximum is state-wise. Both operators are $\gamma$-contractions in the supremum norm, \cite[Prop. 2]{geist2019theory}. The regularized value function of a policy $\pi$, denoted $V_{r}^\pi$, is defined as the fixed point of the operator $T_{\pi,G}$. The regularized value function of the regularized MDP, denoted $V^\star_{r}$, is defined as the fixed point of the regularized optimality operator $T_{\star,G}$. Due to the strong concavity of $G$, the regularized Bellman optimality operator has a unique maximizing argument.
\begin{proposition}\citep[Theorem 1]{geist2019theory}
    \label{prop: optdual-reg}
    The policy $\pi^\star_{r}:=\arg\max_{\pi\in\Pi} T_{\pi,G}V^\star_r$ is the unique optimal regularized policy in the sense that, state-wise, $V^{\pi^\star_r}_r=V^\star_r\geq V^\pi_r$ for all policies $\pi\in\Pi$.
\end{proposition}

The optimal regularized policy $\pi^\star_{r}$ is a strictly positive and exploratory softmax policy for entropy-type $G$. The Lagrangian formulation of the regularized MDP \citep{ying2020note} includes the strongly concave function of the dual variable $g(\rho):=\sum_{s\in\mathcal{S}}G(\rho_s)$ with $\rho_s=\rho(s,\cdot)\in\mathbb{R}^{|\mathcal{A}|}_+$ and is given as
\begin{align}
    \label{eq: dual-reg LP}
    \max_{\rho\geq0} \min_{V\in \mathbb{R}^{|\mathcal{S}|}} L_{r}(V, \rho):= \mu^\top V+\sum_{(s,a)\in\mathcal{S}\times\mathcal{A}} \rho(s, a) \Delta[V](s, a)+\eta_\rho g(\rho).
\end{align}
The following corollary is a consequence of Proposition~\ref{prop: optdual-reg} and characterizes the solution to the regularized Lagrangian. 
\begin{corollary}
    \label{cor: dual-reg}
    Problem~\eqref{eq: dual-reg LP} has the unique saddle-point $(V^\star_{r},\rho^\star_{r})$. By the uniqueness of the optimal regularized policy, we have $\pi^\star_r = \pi_{\rho^\star_r}$.
\end{corollary}
The proof is provided in Appendix~\ref{app:additional:proofs}.
The regularization affects the optimal value, that is, the regularized value function $V^\star_r$ differs from the unregularized value function $V^\star_{ur}$. However, the regularized policy $\pi^\star_r$ approximates the unregularized problem well. More precisely, the suboptimality of $\pi_r^\star$ for the unregularized MDP can be quantified and controlled.

\begin{proposition}\citep[Theorem 2]{geist2019theory}
    Assume that constants $L_G$ and $U_G$ exist such that for all $\pi\in\Pi$ we have $L_G\leq \max_{s\in\mathcal{S}} G(\pi(\cdot|s))\leq U_G$. Then, $V^\star_{ur}(s)-\eta_\rho\frac{U_G-L_G}{1-\gamma}\leq V_{ur}^{\pi^\star_r}(s)\leq V^\star_{ur}(s)$ for all $s\in\mathcal{S}$.
\end{proposition}

For the Shannon entropy regularizer $G(\pi(\cdot|s)):=\sum_{a\in\mathcal A}\pi(a|s)\log(\pi(a|s)),$ we have $0\le G(\pi(\cdot|s))\le\log|\mathcal A|,$ for all $\pi\in\Pi$ and we can take $L_G=0,\ U_G=\log|\mathcal A|$ in the above stated bounds. 
To ensure better computational properties, we seek to make the regularized problem~\eqref{eq: dual-reg LP}, defined by the regularized Lagrangian \( L_r \), strongly convex in \( V \). Since \( L_r \) is originally linear in \( V \), we propose an additional modification to problem~\eqref{eq: dual-reg LP}. Specifically, we replace the linear term \( \mu^\top  V \) with a quadratic term and introduce a specific dual regularization function \( g \) of the form  
\begin{equation}
\label{eq:Entropy:regularizer}
    g(\rho)=-\sum_{(s,a)\in\mathcal{S}\times\mathcal{A}}\rho(s,a)\log\left(\frac{\rho(s,a)}{\tilde{\rho}(s)}\right).
\end{equation}  
With these modifications, we arrive at the main formulation of our study, a further regularized problem that takes the form
\begin{align}
    \label{eq: regLP}
    \max_{\rho\geq0} \min_{V\in \mathbb{R}^{|\mathcal{S}|}} L(V, \rho):= \frac{\eta_V}{2}\|V\|_2^2+\sum_{(s,a)\in\mathcal{S}\times\mathcal{A}} \rho(s, a) \left(\Delta[V](s, a)-\eta_\rho\log\left(\frac{\rho(s,a)}{\tilde{\rho}(s)}\right)\right).
\end{align}
The closely related recent work \citet{li2024accelerating} also studies the convexified regularized MDP~\eqref{eq: regLP}. They prove that the convex modification does not affect the optimal regularized value and only changes the optimal dual variable by a scaling factor that does not influence the optimal regularized policy.
\begin{proposition}\citep[Theorem 2.1]{li2024accelerating} \label{thm:Li-et-al}
    The unique solution $(V^\star,\rho^\star)$ to the double-regularized problem~\eqref{eq: regLP} is such that $V^\star=V^\star_{r}$ and the optimal policy induced by $\rho^\star$ coincides with $\pi^\star_{r}$, that is, $\pi_{\rho^\star}=\pi_{\rho_{r}^\star}=\pi^\star_r$.
\end{proposition}
The optimal dual variable $\rho^\star$ of the double-regularized problem~\eqref{eq: regLP} depends on the convexity parameter $\eta_V$ and does not form a state-action occupancy measure. The following technical lemma provides bounds for the optimal dual variable. Its proof is deferred to Appendix~\ref{app:additional:proofs}.
\begin{lemma}[Optimal dual variable]
    \label{lemma: dualBound}
    The optimal dual variable of problem~\eqref{eq: regLP} with \\bounded reward is strictly positive and bounded from above, that is, there exist $C^L,C^U>0$ such that $C^L<\rho^\star(s,a)< C^U$ for all $(s,a)\in\mathcal{S}\times\mathcal{A}$.
\end{lemma}
Let $C^L$ and $C^U$ be the explicit bounds derived in the constructive proof of Lemma~\ref{lemma: dualBound}. Define $H$ as
\begin{align}
        \label{eq: def.H}
        H:=\{\rho\in\mathbb{R}^{|\mathcal{S}|\cdot|\mathcal{A}|} \ | \ C^L\leq \rho(s,a)\leq C^U, \ \forall (s,a)\in\mathcal{S}\times\mathcal{A}\}.
\end{align}
It holds that $\rho^\star\in \operatorname{int} H$. This compact set $H$ will be used in the algorithms to specify the projection region of the dual variable iterates without impacting the optimal solution.

\section{Synchronous Algorithm}\label{sec:synchronous:setting}

We now propose our solution to the saddle-point problem~\eqref{eq: regLP} that is based on a projected stochastic gradient descent-ascent method summarized in Algorithm~\ref{alg1}. The algorithm does not assume to know the transition probabilities of the MDP but has access to a generative model \citep{ref:kakade2003:sample,chen2016}, also known as a sampling oracle, which provides sampling access to any state in the environment. Given a state-action pair $(s,a)$ as input, the generative model provides the reward $r(s,a)$ and a follow-up state $s'$ with $s'\sim\mathcal{P}(\cdot| s,a)$. The almost sure convergence analysis is based on the two-timescale stochastic approximation framework developed in \citet{borkar1997stochastic}, \citet{borkar2023stochastic}. The dynamics of the projected recursion are based on results of \citet{dupuis1987large} and \citet{dupuis1993dynamical}.

Since the true gradients of the Lagrangian $L$ depend on the unknown transition probabilities, we use stochastic gradient approximations, denoted by $\widehat{\nabla}_VL$ and $\widehat{\nabla}_\rho L$, to perform the updates. To ensure that the iterates remain bounded, we project the dual iterates onto the feasible set $H$ using the bounds on the optimal dual variable established in Lemma~\ref{lemma: dualBound}. In the synchronous setting, we fully update the primal and dual variables in each iteration. We construct the stochastic gradients from generative model transition samples obtained for each state-action pair. 

Specifically, let \(\widetilde\Xi := \{\tilde s(s,a)\}_{(s,a)\in\mathcal S\times\mathcal A}\) be a collection of random variables such that
\[
\tilde s(s,a)\sim \mathcal P(\cdot| s,a),
\quad\text{and }\{\tilde s(s,a)\}_{(s,a)\in\mathcal S\times\mathcal A}\text{ are independent.}
\]
Recall the definition of the state-marginal \(\tilde\rho(s):=\sum_{a\in\mathcal A}\rho(s,a)\).
Given \(\widetilde\Xi\), the stochastic gradients used in Algorithm~\ref{alg1} are
\begin{subequations}\label{eq:sync:stochgrad}
\begin{align}
\label{eq:VStochgrad}
\widehat{\nabla}_{V}L(V,\rho)(s')
&= \eta_V V(s')-\tilde\rho(s')
+ \gamma\sum_{(s,a)\in\mathcal S\times\mathcal A}\rho(s,a)\,\mathbf 1\{\tilde s(s,a)=s'\},
\qquad s'\in\mathcal S,
\\
\label{eq:rstochgrad}
\widehat{\nabla}_{\rho}L(V,\rho)(s,a)
& = -V(s)+r(s,a)+\gamma V\bigl(\tilde s(s,a)\bigr)
-\eta_\rho\log\Bigl(\frac{\rho(s,a)}{\tilde\rho(s)}\Bigr),
\ (s,a)\in\mathcal S\times\mathcal A.
\end{align}
\end{subequations}

\begin{algorithm}[t]
    \caption{Synchronous Projected Gradient Descent-Ascent RL} 
    \label{alg1}
    \begin{algorithmic}[1]
            \Require $H$ as in~\eqref{eq: def.H}, $V_0\in \mathbb{R}^{|\mathcal{S}|},\rho_0\in H,\eta_\rho,\eta_V>0,\{\alpha_i\}_{i=1}^K,\{\beta_i\}_{i=1}^K\subset\mathbb{R}_{>0}$
        \For {$k=1,2,\ldots, K$}
            \State Sample state-transitions $s'\sim\mathcal{P}(\cdot| s,a)$ for all $(s,a)\in\mathcal{S}\times\mathcal{A} $
            \State Update $V$ with an SGD step $V_k=V_{k-1}-\alpha_k\widehat{\nabla}_{V}L(V_{k-1},\rho_{k-1})$
            \State Update $\rho$ with a projected SGA step $\rho_k=\Pi_H\left[\rho_{k-1}+\beta_k\widehat{\nabla}_{\rho}L(V_{k-1},\rho_{k-1})\right]$
        \EndFor
    \end{algorithmic} 
        \Return $V_K, \rho_K$
\end{algorithm}

Algorithm~\ref{alg1} is a stochastic approximation algorithm based on the primal-dual gradient method for saddle-point estimation \citep{ref:arrow1958:studies}. \cite{chen2016}, \cite{ref:lee2019:stochastic}, among others, study versions of a related approach for unregularized LP formulations, where the optimal policy's uniqueness must be assumed. 
They provide sample complexity bounds by applying a concentration bound to the average iterate. In contrast, we take the dynamical systems viewpoint \citep{borkar2023stochastic} on the stochastic approximation algorithm and analyze the asymptotic almost sure convergence. This change in perspective is reflected in the choice of stepsizes, where our approaches feature diminishing stepsizes operating on two timescales, such that the last iterates asymptotically approach the limiting ODEs. In contrast, they rely on a fixed stepsize scheme, which suffices for the average iterate analysis.

We first present this synchronous algorithm to introduce the convergence analysis and generalize the setting in Section~\ref{sec:asynchronous:setting} to allow asynchronous updates based on observations generated by a single trajectory of state-action pairs under a fixed or varying behavioral policy. 

\subsection{Almost Sure Convergence Analysis of Algorithm~\ref{alg1}}\label{sec:convergence:algo1}

\begin{theorem}[Convergence of Algorithm~\ref{alg1}]\label{thm:Conv:Alg1}
Consider the stochastic gradients as defined in~\eqref{eq:sync:stochgrad} and the projection region $H$ defined in~\eqref{eq: def.H}. Let the stepsize sequences $\{\alpha_k\},\{\beta_k\}$ satisfy $$\sum_{k=1}^\infty \alpha_k=\sum_{k=1}^\infty\beta_k=\infty, \quad \sum_{k=1}^\infty\left(\alpha_k^2+\beta_k^2\right)<\infty, \quad \lim_{k\rightarrow\infty}\frac{\beta_k}{\alpha_k}=0.$$ Let $\{(V_k,\rho_k)\}$ be a sequence generated by Algorithm~\ref{alg1} with an arbitrary initialization $V_0\in \mathbb R^{|\mathcal S|},\ \rho_0\in H$. Then the last iterates converge $\mathbb{P}$-a.s. to the regularized saddle-point of~\eqref{eq: regLP}, that is, $$\lim_{K\rightarrow\infty}(V_K,\rho_K)=(V^\star,\rho^\star)\quad \mathbb{P}\text{-a.s.}$$
\end{theorem}

Examples of suitable stepsize sequences are $\alpha_k=\frac{1}{k^q}, q\in(\frac{1}{2},1)$ and $\beta_k=\frac{1}{k}$ or $\alpha_k=\frac{1}{k}$, $\beta_k=\frac{1}{1+k\log k}$.
The proof of Theorem~\ref{thm:Conv:Alg1} builds on the two-timescale stochastic approximation framework of~\cite{borkar1997stochastic}. This framework takes the dynamical systems viewpoint on the stochastic iteration and applies to the regularized MDP problem addressed with Algorithm~\ref{alg1} as follows.

The recursion in Algorithm~\ref{alg1} corresponds to the following coupled iterations
\begin{subequations}\label{eq:sync:recursion}
\begin{align}
    \label{eq:sync:V:recursion}
    V_k &= V_{k-1}-\alpha_k(\nabla_VL(V_{k-1},\rho_{k-1})+M^{(1)}_k),\\
    \label{eq:sync:rho:recursion}
    \rho_k &= \Pi_H\left[\rho_{k-1}+\beta_k(\nabla_\rho L(V_{k-1},\rho_{k-1})+M^{(2)}_k)\right],
\end{align}
\end{subequations}
where $\nabla_VL:\mathbb{R}^{|\mathcal{S}|}\times H\rightarrow\mathbb{R}^{|\mathcal{S}|}$, $\nabla_\rho L:\mathbb{R}^{|\mathcal{S}|}\times H\rightarrow\mathbb{R}^{|\mathcal{S}|\cdot|\mathcal{A}|}$ are the gradients of the regularized Lagrangian and $\{M^{(1)}_k\}$, $\{M^{(2)}_k\}$ are the sequences of differences between the true gradients and their noisy estimates given in~\eqref{eq:sync:stochgrad}, that is,
\begin{align*}
    M^{(1)}_k = \widehat{\nabla}_{V}L(V_{k-1},\rho_{k-1}) - \nabla_{V}L(V_{k-1},\rho_{k-1}),\, 
     M^{(2)}_k = \widehat{\nabla}_{\rho}L(V_{k-1},\rho_{k-1}) - \nabla_{\rho}L(V_{k-1},\rho_{k-1}). 
\end{align*}
The iterate sequences of Algorithm~\ref{alg1} fulfill the following properties regarding the gradients and noise.

\begin{proposition}[Gradient and noise properties of Algorithm~\ref{alg1}] \label{prop: sync:gradient:noise}
The following\\ conditions are satisfied.
    \begin{enumerate}
        \item[\textbf{(A1)}]\label{prop-A1} The regularized Lagrangian, defined in~\eqref{eq: regLP}, has Lipschitz continuous gradients.
        \item[\textbf{(A2)}]\label{prop-A2} The noise terms $\{M_k^{(i)}\}$ for $i=1,2$ arising from the iterates of Algorithm~\ref{alg1} are martingale difference sequences with respect to the increasing $\sigma$-fields
    $$\mathcal{F}_k:=\sigma(V_\ell,\rho_\ell,M^{(1)}_\ell,M^{(2)}_\ell,\ell\leq k),\quad k\geq0,$$
    and satisfy $$\mathbb{E}\left[\|M^{(i)}_k\|^2|\mathcal{F}_{k-1}\right]\leq C(1+\|V_{k-1}\|^2+\|\rho_{k-1}\|^2),\quad i=1,2,$$ for a constant $C>0$ for $k\geq 1$.
    \end{enumerate}
\end{proposition}
The proof is provided in Appendix~\ref{sec: synchronous-alg-conv-proof}.

Together with the two-timescale stepsize sequences assumed in Theorem~\ref{thm:Conv:Alg1}, these gradient and noise properties ensure that the coupled system~\eqref{eq:sync:recursion} asymptotically approximates a dynamical system with a fast transient $V(t)$ and a slowly evolving dual variable $\rho(t)$. The separation of the timescales allows the fast transient to view the dual variable as quasi-constant, and in turn, the dual variable to assume the transient to be equilibrated. The convergence of the stochastic recursion is derived from the asymptotic properties of the limiting dynamical system. The ODE describing the dynamics of the primal variable $V$ is the gradient flow 
\begin{subequations}
\label{eq:sync:ODEs}
    \begin{align}
        \label{eq: Vode}
        \dot{V}(t)=-\nabla_VL(V(t),\rho)
    \end{align}
     for fixed $\rho\in H$. Since the mapping $V\mapsto L(V,\rho)$ is strictly convex, the set of equilibrated primal variables given a fixed $\rho\in H$ is a singleton. We denote this ``optimal response'' as $$\lambda(\rho):={\arg\min}_{V\in\mathbb{R}^{|\mathcal S|}} L(V,\rho)=\{V\in\mathbb{R}^{|\mathcal{S}|}: \nabla_V L(V,\rho) = 0\}.$$ 
     
    The dual variable's dynamic is restricted to the compact set $H$. We follow the notation of \cite[Section 4.3]{kushner2003} for the projected dynamics. For $\rho\in \partial H$, the boundary of $H$, let $C(\rho)$ be the infinite convex cone generated by the outer normals at $\rho$ of the faces on which $\rho$ lies. For $\rho$ in the interior of $H$, let $C(\rho)=\{0\}$. The projected ODE describing the flow of the dual variable is then given by 
    \begin{align}
        \label{eq: rhoode}
        \dot{\rho}(t)=\nabla_\rho L(\lambda(\rho(t)),\rho(t))+\zeta(t), \quad \zeta(t)\in-C(\rho(t)),
    \end{align}
    with the \textit{projection term} $\zeta$, the minimum force needed to keep $\rho$ in $H$. To obtain almost sure convergence of the stochastic recursion, the limiting dynamical systems need to meet asymptotic stability criteria that are fulfilled by the recursion of Algorithm~\ref{alg1}.
\end{subequations}

\begin{proposition}[Dynamical systems properties of Algorithm~\ref{alg1}]\label{prop: ODEs}
The following\\ conditions are satisfied.
    \begin{enumerate}
        \item[\textbf{(A3)}]\label{prop-A3} For each $\rho\in H$, the primal ODE (\ref{eq: Vode}) has a globally asymptotically stable equilibrium $\lambda(\rho)$, with $\lambda:H\rightarrow\mathbb{R}^{|\mathcal{S}|}$ a Lipschitz map.
        \item[\textbf{(A4)}]\label{prop-A4} The projected dual ODE (\ref{eq: rhoode}) has a globally asymptotically stable equilibrium $\rho^\star$.
        \item[\textbf{(A5)}]\label{prop-A5} Let $\{(V_k,\rho_k)\}$ be a sequence of iterates generated by~\eqref{eq:sync:recursion} with a stepsize sequence fulfilling the assumption of Theorem~\ref{thm:Conv:Alg1}. Then $\sup_k(\|V_k\|+\|\rho_k\|)<\infty,$ a.s.
    \end{enumerate}
\end{proposition}

\paragraph{Proof Sketch} 
The full proof is stated in Appendix~\ref{sec: synchronous-alg-conv-proof}. The main steps are as follows.
To prove the convergence of the fast transient, $V(t)\rightarrow \lambda(\rho)$ for a fixed dual variable $\rho$ stated in~\hyperref[prop-A3]{(A3)}, we provide a strict Lyapunov function argument.
We show the convergence of the slow component given the equilibrated fast transient,~\hyperref[prop-A4]{(A4)}, by stating the well-posed projected dynamic system using results found in \citet{dupuis1993dynamical} and again constructing a strict Lyapunov function to establish global asymptotic stability with \citet[Theorem 4.6]{dupuis1987large}. 
Lastly, to prove the iterate stability property,~\hyperref[prop-A5]{(A5)}, we apply the scaling limit or ODE@$\infty$ method, established in \citet{borkar2000ode}. The dual estimators are bounded by projection to the compact hyperrectangle $H$ with $\rho^\star\in H$.
\hfill $\blacksquare$

\vspace{11pt}
Given the above stability properties, existing stochastic approximation theory by \cite{borkar2023stochastic} ensures the asymptotic convergence as follows.

\paragraph{Proof of Theorem~\ref{thm:Conv:Alg1}}
Given Propositions~\ref{prop: sync:gradient:noise}, and~\ref{prop: ODEs}, and the two-timescale stepsize sequences assumed in Theorem~\ref{thm:Conv:Alg1}, \citep[Theorem 8.1]{borkar2023stochastic} applies and the coupled iterates $(V_k,\rho_k)$ generated by the recursions~\eqref{eq:sync:recursion}, and hence by Algorithm~\ref{alg1}, converge $\mathbb{P}\text{-a.s.}$ to $(\lambda(\rho^\star),\rho^\star)$.
\hfill $\blacksquare$

\section{Asynchronous Algorithm}\label{sec:asynchronous:setting}
Assuming access to a generative model, as in Algorithm~\ref{alg1}, can be restrictive in applications where samples can only be obtained by interacting with the environment. We therefore extend the approach to the standard online setting in which data are generated along a trajectory induced by a (fixed or time-varying) behavioral policy. This covers classical online methods such as Q-learning~\citep{ref:watkins1992:q}. The fixed behavioral policy setting is also related to offline RL where policies are learned from a fixed set of past experiences, see~\cite{ref:levine2020:offline} for an overview and~\cite{ref:gabbianelli2024:offline,ref:hong2024:primal} for LP-based formulations. 

Our algorithm retains the projected stochastic gradient descent-ascent structure, but performs \emph{asynchronous} updates: only the state-action components visited by the behavioral policy are updated at each iteration. The almost-sure convergence analysis follows the asynchronous stochastic approximation framework of~\cite{ref:konda1999:actor, ref:perkins2013:asynchronous}, while our finite-time rate bounds build on the two-timescale analysis developed in~\cite{ref:zeng2024:two}.

In addition to the related works mentioned in the discussion of Algorithm~\ref{alg1}, the asynchronous algorithm shares features with the algorithms SBEED \citep{ref:dai2018:sbeed} and Dual-AC \citep{ref:Dai-18:boosting}. They approximate different saddle-point formulations of the RL problem with behavioral policy access to the stochastic environment by sampling transitions from a growing experience replay buffer. In their convergence analysis, they assume a fixed behavioral policy. Our proposed algorithm combines the experience replay-based gradient estimation with the separation of the nested loop minimization by the two-timescale method and on-policy exploration. This combination allows us to apply the asynchronous stochastic approximation analysis in a way that has not been proposed for the regularized LP approach to reinforcement learning.

The asynchronous stochastic approximation setting features two sources of randomness: first, the index selection for the updates, which is based on the random transitions following the behavioral policy, and second, the stochastic gradient estimators that incorporate information from previous transitions. 

We propose the following procedure for the index selection at a given update step in Algorithm~\ref{alg2}. Let $\pi^b_k\in\Pi$ for $k\ge 0$ denote the behavioral policy at iteration $k$, and $\{X_k\}\subset\mathcal X$ the state-action process generated by the sequence of behavioral policies $\{\pi_k^b\}$ where $$X_k:=(s_{k-1},a_{k-1}),\ s_k\sim \mathcal P(\cdot|s_{k-1},a_{k-1}),\ a_k\sim \pi_{k-1}^b(\cdot|s_k).$$
At iteration $k$, the primal and dual variables are updated with projected stochastic gradient steps in the entries corresponding to the state $s_k$ and state-action pair $X_k$, respectively.

Without the generative model, the trajectory $\{X_k\}$ is our only interaction with the environment. Since we cannot obtain new unbiased samples of every transition in each iteration under this Markovian observational model, we propose gradient estimators based on a version of experience replay to obtain asymptotically unbiased gradient directions. The experience replay buffer stores the observed transitions, which we use to estimate the transition probabilities featured in the true gradients as follows. 

At each iteration, we add the current transition of the behavioral policy to the experience replay buffer, denoted as $\mathcal{D}_k$ in iteration $k$. To efficiently leverage the gathered experience, we opt for a structured buffer with a separate list $\mathcal{D}_k(s,a)$ for each state-action pair $(s,a)$. Furthermore, for each state $s'$, we maintain a set, denoted as $\mathcal{D}_{\mathrm{inc}}(s')$, containing the ``incoming" state-action pairs that have previously led to state $s'$.\footnote{The set of incoming state-action pairs does not contain duplicates. In contrast, the replay buffer consists of lists that allow duplicate entries.} Hence, if $X_k\notin\mathcal{D}_{\mathrm{inc}}(s_k)$, we update $\mathcal{D}_{\mathrm{inc}}(s_k):=\mathcal{D}_{\mathrm{inc}}(s_k)\cup\{X_k\}$.

We estimate the transition probability vector $\mathcal P(\cdot|s_{k-1},a_{k-1})$ in $\nabla_\rho L(V,\rho)(s_{k-1}, a_{k-1})$ with the unbiased new transitions $s_k$ as the one-hot vector $e_{s_k}$. For the primal variable gradient entry $\nabla_VL(V,\rho)(s_k)$, we use buffer-based estimates of $\mathcal{P}(s_k|s,a)$ for state-action pairs with a previously observed transition to $s_k$.
Concretely, for each $(s,a)\in\mathcal D_{\mathrm{inc}}(s_k)$ we draw one next-state sample $\xi_k(s,a)\sim \mathcal P_{\mathcal D_k}(\cdot| s,a)$ from the corresponding buffer, where $\mathcal P_{\mathcal D_k}(s'|s,a):=\nu_k(s,a)^{-1}\sum_{\bar s\in\mathcal D_k(s,a)}\mathbf{1}\{\bar s=s'\}$, and set $$\widehat{\mathcal P}_k(s_k| s,a)=\mathbf 1\{\xi_k(s,a)=s_k\}.$$
For $(s,a)\notin\mathcal D_{\mathrm{inc}}(s_k)$ no draw is made and we define $\widehat{\mathcal P}_k(s_k| s,a)=0$.
Denote the combined buffer-draw random variable as $\Xi_k$. The structured buffer draw can be interpreted as sampling from the empirical transition kernel. Under infinite visitation, this empirical kernel converges a.s. to the true kernel, hence the estimator is consistent.

\begin{algorithm}[t]
    \caption{Asynchronous PGDA-RL with Structured Experience Replay} 
    \label{alg2}
    \begin{algorithmic}[1]
            \Require $H$ as in~\eqref{eq: def.H}, $V_0\in \mathbb{R}^{|\mathcal{S}|},\rho_0\in H,\eta_\rho,\eta_V>0,\{\alpha_k\}_{k\in\mathbb N},\{\beta_k\}_{k\in\mathbb N}\subset\mathbb{R}_{>0},\pi^b_0\in\Pi,\mu\in\Delta_\mathcal{S}, \{\epsilon_k\}_{k=1}^K\subset[0,1]$
        \State Initialize visitation counters: $\tilde\nu_0=\textbf{0}\in\mathbb{R}^{|\mathcal{S}|},\nu_0=\textbf{0}\in\mathbb{R}^{|\mathcal{S}|\cdot|\mathcal{A}|}$
        \State Initialize the structured experience replay buffer: For all $(s,a)\in\mathcal{S}\times\mathcal{A}$ set $\mathcal{D}_0(s,a) = \emptyset$
        \State Initialize incoming states sets: For all $s\in\mathcal{S}$ set $\mathcal{D}_{\mathrm{inc}}(s)=\emptyset$
        \State Sample initial state-action pair: $s_0\sim\mu,a_0\sim\pi^b_0(\cdot| s_0)$
        \For {$k=1,2,\ldots, K$}
            \State Set $X_k=(s_{k-1},a_{k-1})$
            \State Sample state-transition: $s_k\sim\mathcal{P}(\cdot| X_k)$
            \State Sample action with behavioral policy: $a_{k}\sim\pi^b_{k-1}(\cdot| s_k)$
            \State Update visitation counters:
            $\tilde\nu_k=\tilde\nu_{k-1}+e_{s_{k}}, \quad \nu_k=\nu_{k-1}+e_{X_k}$
            \State Update experience replay buffers: $$\mathcal{D}_k(X_k)=\mathcal{D}_{k-1}(X_k)\cup\{s_k\} \text{ (and } \mathcal{D}_k(s,a)=\mathcal{D}_{k-1}(s,a) \text{ for } (s,a)\neq X_k)$$
            \State \textbf{if} $X_k\notin \mathcal D_{\mathrm{inc}}(s_k)$ \textbf{then} $\mathcal D_{\mathrm{inc}}(s_k)=\mathcal D_{\mathrm{inc}}(s_k)\cup\{X_k\}$
            \State For each $(s,a) \in \mathcal{D}_{\mathrm{inc}}(s_k)$ draw a transition $\xi_k(s,a)$ uniformly from $\mathcal{D}_{k}(s,a)$, set $$\widehat{\mathcal{P}}_k(s_k|s,a)=\mathbf{1}\{\xi_k(s,a)=s_k\}$$ 
            \State Update $V$ with an SGD step in the entry $s_{k}$ with $\widehat g_k$ as in~\eqref{eq: Vstochgrad_alg2}: 
            $$V_k=V_{k-1}-\alpha\big(\tilde\nu_k(s_{k})\big)\widehat g_k(\rho_{k-1},V_{k-1};X_k,\Xi_k)$$
            \State Update $\rho$ with a projected SGA step in the entry $X_k$ with $\widehat h_k$ as in~\eqref{eq: rstochgrad_alg2}:
            $$\rho_k=\Pi_H\left[\rho_{k-1}+\beta\big(\nu_k(X_k)\big)\widehat h_k(\rho_{k-1},V_{k-1};X_k)\right]$$
            \State \textit{(Optional)} Update the behavioral policy:\Comment{On-Policy Exploration}           
            $$\pi_k^b(a|s)=(1-\epsilon_k)\frac{\rho_k(s,a)}{\sum_{a\in\mathcal{A}}\rho_k(s,a)}+\epsilon_k\hspace{1pt}\mathcal{U}(\mathcal{A)},\qquad (s,a)\in\mathcal{S}\times\mathcal{A}$$
        \EndFor
    \end{algorithmic} 
        \Return $V_K, \rho_K$
\end{algorithm}
Summarizing the above, the stochastic gradients for the asynchronous algorithm at iteration $k$ are 
\begin{subequations}\label{eq:async:stochgrad}
    \begin{align}
    \label{eq: Vstochgrad_alg2}
        \widehat g_k(\rho,V;X_k,\Xi_k)&:=e_{s_k}\left(\eta_VV(s_k)-\tilde\rho(s_k)+\gamma\sum_{(s,a)\in\mathcal{D}_{\mathrm{inc}}(s_k)}\rho(s,a)\mathbf{1}\{\xi_k(s,a)=s_k\}\right)\\
    \label{eq: rstochgrad_alg2}
        \widehat h_k(\rho, V; X_k)&:=e_{X_k}\left(-V(s_{k-1})+r(X_k)+\gamma V(s_{k})-\eta_\rho\log\big(\tfrac{\rho(X_k)}{\tilde\rho(s_{k-1})}\big)\right),
    \end{align}
\end{subequations}
where $\widehat g_k$ and $\widehat h_k$ are the estimates of $\nabla_VL$ and $\nabla_\rho L$, respectively.   

To ensure that each entry receives the same sequence of stepsizes when updates occur, we adjust the effective stepsizes based on visitation counts. We denote the number of visits to $(s,a)\in\mathcal X$ until iteration $k$ as $\nu_k(s,a):=\sum_{\ell\le k}\mathbf{1}\{X_\ell=(s,a)\}$ and correspondingly the number of visits to state $s$ is the state marginal $\tilde \nu_k(s):=\sum_{a\in\mathcal A}\nu_k(s,a)$. Then, the asynchronous update at iteration $k$ applies the random scalar stepsize $\bar{\alpha}_k:=\alpha(\tilde{\nu}_k(s_k))$ and $\bar{\beta}_k:=\beta(\nu_k(X_k))$ for the primal and dual updates respectively.\footnote{Note that $\beta(k)$ and $\beta_k$ both denote the $k^{th}$ stepsize sequence element. Hence $\beta(\nu_k(s,a))$ is entry number $\nu_k(s,a)\in\mathbb N$ of the stepsize sequence $\{\beta(k)\}$.} The visitation counts act as ``local clocks''. 

Rather than performing updates along the trajectory of a fixed behavior policy, it is often desirable to explore the state-action space using the current policy iterate in an online manner. The updating of the behavioral policy in Algorithm~\ref{alg2} allows for exploiting the learned policy during training while maintaining sufficient exploration through the projection of the dual iterate. While data collection follows an on-policy scheme, the updates are still performed off-policy using past experiences stored in the replay buffer. This approach ensures more stable off-policy learning with effective use of observed data.

\subsection{Almost Sure Convergence Analysis of the Algorithm~\ref{alg2}}
\label{ssec: async-as-conv}
We make the following assumptions to ensure sufficient exploration of the state-action space and obtain convergence results.

\begin{assumption}[Ergodicity]\label{ass:ergodicity}
    For any strictly exploratory policy $\pi$ (that is, $\pi(a|s)>0$ for all $(s,a)\in\mathcal{S}\times\mathcal{A}$), the induced Markov chain over the state space is ergodic. 
\end{assumption}

This is weaker than assuming ergodicity for all stationary policies. It further ensures that Markov chains on the state-action space induced by strictly exploratory policies are also ergodic. 

\begin{assumption}[Asynchronous Step-Size Sequences]\label{ass:async:stepsize}
    The stepsize sequences $\{\alpha_k\}$,\\$\{\beta_k\}$ satisfy 
    \begin{enumerate}
            \item $\sum_{k=1}^\infty\alpha_k=\sum_{k=1}^\infty\beta_k=\infty, \quad \sum_{k=1}^\infty\left(\alpha_k^2+\beta_k^2\right)<\infty, \quad \lim_{k\to\infty} \frac{\beta_k}{\alpha_k} = 0.$
            \item For all $k$, for both $a_k=\alpha_k$ and $a_k=\beta_k$, $a_k\geq a_{k+1}$ and for all $x\in(0,1)$, there exists $A_x$ with $\sup_k\frac{a_{[xk]}}{a_k}<A_x<\infty$.
    \end{enumerate}     
\end{assumption}

With the above assumptions, we obtain the following almost sure convergence result for both versions of Algorithm~\ref{alg2}.

\begin{theorem}[Almost Sure Convergence of Algorithm~\ref{alg2}]\label{thm:Conv:Alg2}
    Consider the stochastic\\ gradients defined in~\eqref{eq:async:stochgrad}. Let Assumption~\ref{ass:ergodicity} hold and choose stepsize sequences $\{\alpha_k\}$,$\{\beta_k\}$ that satisfy Assumption~\ref{ass:async:stepsize}.  Let $\{(V_k,\rho_k)\}$ be a sequence generated by Algorithm~\ref{alg2} either with on-policy exploration or with a fixed strictly exploratory behavioral policy. Then, the last iterates converge $\mathbb{P}$-a.s. to the regularized saddle-point of~\eqref{eq: regLP}, that is,
    \begin{align*}
        \lim_{K\to\infty} (V_K,\rho_K) = (V^\star,\rho^\star) \quad \mathbb{P}\text{-a.s.}
    \end{align*}
\end{theorem}

The proof of Theorem~\ref{thm:Conv:Alg2} is based on the asynchronous stochastic approximation theory by \cite{ref:borkar1998:asynchronous}, \cite{ref:konda1999:actor}, and the generalization in \cite{ref:perkins2013:asynchronous} that clarifies the conditions on the update frequencies. 

To analyze the dynamics of the iterates generated by Algorithm~\ref{alg2}, we define the filtrations 
$$\mathcal F_k:=\sigma(V_\ell,\rho_\ell,s_\ell,a_\ell,\Xi_\ell;\ell\leq k)\quad \text{ and }\quad\mathcal F_k^-:=\sigma(\mathcal F_{k-1},s_k,a_k),\ k\ge 1.$$
The first filtration contains the complete history of the algorithm, while the latter filtration captures the information prior to sampling $\Xi_k$ from the buffer $\mathcal{D}_k$ at iteration $k$. Indeed, $\mathcal{D}_k$ is a measurable function of $\{(s_\ell,a_\ell)\}_{\ell\le k}$, hence $\mathcal{D}_k$ is $\mathcal F_k^-$-measurable. Clearly, $\mathcal F_{k-1}\subset\mathcal F_k^-\subset\mathcal F_k$.

We redefine the noise sequences $\{M_k^{(1)}\}$ and $\{M_k^{(2)}\}$ as 
\begin{subequations}\label{eq:SA:def:MG}
\begin{align}
    \label{eq:SA:def:MD1}
    M^{(1)}_k&:=\gamma\sum_{(s,a)\in\mathcal{S}\times\mathcal{A}}\rho_{k-1}(s,a)\left(\widehat{\mathcal{P}}_{k}(\cdot| s,a)-\mathbb{E}[\widehat{\mathcal{P}}_k(\cdot|s, a)|\mathcal F_k^-]\right),\\
    \label{eq:SA:def:MD2}
    M_k^{(2)}&:=e_{X_k} \gamma \left(e_{s_k}-\mathcal{P}(\cdot| X_k)\right)^\top V_{k-1},
\end{align} 
\end{subequations}
where $\widehat {\mathcal P}_k(\cdot|s,a)\equiv0$ for pairs not drawn at time $k$. 
Note that $V_{k-1}$ and $X_k=(s_{k-1},a_{k-1})$ are $\mathcal F_{k-1}$-measurable and, since
$s_k\sim \mathcal P(\cdot| X_k)$, we have $\mathbb E[e_{s_k}| \mathcal F_{k-1}]=\mathcal P(\cdot| X_k)$.
Consequently, $\mathbb E[M_k^{(2)}| \mathcal F_{k-1}]=0$. 
We further introduce the bias sequence 
\begin{equation}\label{eq:SA:def:bias}
    \begin{aligned}
        \mathcal{E}_k:=\gamma\sum_{(s,a)\in\mathcal{S}\times\mathcal{A}}\rho_{k-1}(s,a)\left(\mathbb{E}[\widehat{\mathcal{P}}_k(\cdot|s, a)|\mathcal F_k^-]-\mathcal{P}(\cdot|s,a)\right).
    \end{aligned} 
\end{equation}
Note that the noise and bias sequences depend on the current iterates. When necessary, this dependence is made explicit.
Denote the square diagonal matrices with the update index selectors as the diagonal entries as $\Lambda^{(1)}_k:=\text{diag}(e_{s_k}), \Lambda^{(2)}_k:=\text{diag}(e_{X_k})$. 
The asynchronous stochastic gradients~\eqref{eq:async:stochgrad} satisfy
\begin{subequations}
\label{eq:async:recursion}
\begin{align}
    \label{eq:async:V:recursion}
    \widehat g_k(\rho_{k-1},V_{k-1};X_k,\Xi_k)&=\Lambda^{(1)}_k\big(\nabla_VL(V_{k-1},\rho_{k-1})+M^{(1)}_k+\mathcal{E}_k\big),\\
    \label{eq:async:rho:recursion}
    \widehat h_k(\rho_{k-1}, V_{k-1}; X_k)&=\Lambda^{(2)}_k\big(\nabla_\rho L(V_{k-1},\rho_{k-1})+M^{(2)}_k\big),\qquad \text{for all }k\in\mathbb N.
\end{align}    
\end{subequations}

Similarly to the synchronous case, the asymptotic analysis of the asynchronous stochastic iteration consists of analyzing the equilibrium properties of the dynamical systems that describe the limiting behavior of the iteration. 
In the asynchronous setting, the random selection of update components impact the limiting dynamics. 
Following the notation of~\cite{ref:perkins2013:asynchronous}, define for general $\epsilon>0,d\in\mathbb{N}$, $\Omega^\epsilon_d:=\{\text{diag}(w_1,\dots,w_d)|w_i\in[\epsilon,1],i=1,\dots,d\}$. 
For the fast iterate $V$, we are concerned with the ODE
\begin{subequations}\label{eq:async:ODE}
\begin{align}
    \label{eq:async:V:ODE}
    \dot{V}(t)=-\Lambda^{(1)}(t)\nabla_VL(V(t),\rho),
\end{align}
where $\Lambda^{(1)}:\mathbb{R}_{\geq 0}\to\Omega^\epsilon_{|\mathcal{S}|}$ for some $\epsilon>0$ is a measurable function. 
As in the synchronous case, we denote by $\lambda(\rho)$ the equilibrated value of $V$ for a given $\rho$, that is, $\lambda(\rho) = \{V\in\mathbb{R}^{|\mathcal{S}|}: \nabla_V L(V,\rho) = 0\}$ and recall that since the mapping $V\mapsto L(V,\rho)$ is strictly convex, $\lambda(\rho)$ is unique.
For the slower dual iterate, the ODE of interest is
\begin{align}
    \label{eq:async:rho:ODE}
    \dot{\rho}(t)=\Lambda^{(2)}(t)\nabla_\rho L(\lambda(\rho(t)),\rho(t))+\zeta(t), \quad \zeta(t)\in C(\rho(t)),
\end{align}
\end{subequations}
with reflection term $\zeta$ and where $\Lambda^{(2)}:\mathbb{R}_{\geq 0}\to\Omega^\epsilon_{|\mathcal{S}|\cdot|\mathcal{A}|}$ for some $\epsilon>0$ is a measurable function.

If we do not specify the diagonal matrix-valued functions $\Lambda^{(1)}(\cdot)$ and $\Lambda^{(2)}(\cdot)$ further, we can equivalently look at the properties of the class of dynamics captured by the differential inclusions 
\begin{subequations}
    \label{eq:diff:incl}
    \begin{align}
    \label{eq:async:V:diff:incl}
        \dot{V}(t)&\in-\Omega^\epsilon_{|\mathcal{S}|}\cdot\nabla_VL(V(t),\rho),\\
   \label{eq:async:rho:diff:incl}
        \dot{\rho}(t)&\in\Omega^\epsilon_{|\mathcal{S}|\cdot|\mathcal{A}|}\cdot\nabla_\rho L(\lambda(\rho(t)),\rho(t))+\zeta(t), \quad \zeta(t)\in C(\rho(t)).
    \end{align} 
\end{subequations}
Statements regarding the asymptotic properties of any solution to the differential inclusions defined in~\eqref{eq:diff:incl} also apply to the non-autonomous systems~\eqref{eq:async:ODE}. 

To relate the stochastic iterations to the dynamical systems and obtain convergence, we use the following properties of the (stochastic) gradients, the dynamical systems, and the update scheme. 

\begin{proposition}[Asynchronous SA]\label{prop:async:SA}\hspace{2pt}
Let $\{(V_k,\rho_k)\}$ be a sequence generated by Algorithm~\ref{alg2}. The following properties hold.
    \begin{enumerate}
        \item[\textbf{(B1)}]\label{prop-B1} For the noise terms in~\eqref{eq:SA:def:MG} and the bias in~\eqref{eq:SA:def:bias} it holds that
        \begin{enumerate}
            \item $\{M_k^{(1)}\}$ and $\{M_k^{(2)}\}$ are martingale-difference sequences with respect to $\{\mathcal F_k\}$, that is, $\mathbb E[M_k^{(i)}|\mathcal F_{k-1}]=0, \ i=1,2$.
            \item For a constant $C>0$ and $k\geq1$ 
            \begin{align*}
                &\mathbb{E}\left[\|M^{(1)}_k\|^2|\mathcal F_k^-\right]\leq C(1+\|V_{k-1}\|^2+\|\rho_{k-1}\|^2),\\
                &\mathbb{E}\left[\|M^{(2)}_k\|^2|\mathcal F_{k-1}\right]\leq C(1+\|V_{k-1}\|^2+\|\rho_{k-1}\|^2).
            \end{align*}
            \item The bias sequence $\{\mathcal{E}_k\}$ is a bounded random sequence satisfying $\lim_{k\rightarrow\infty}\mathcal{E}_k=0$  $\mathbb P$-$a.s.$
        \end{enumerate}
        
        \item[\textbf{(B2)}]\label{prop-B2} For each $\rho\in H$ and every $\epsilon>0$, for all $\Omega^\epsilon_{|\mathcal{S}|}$-valued measurable functions $\Lambda^{(1)}(\cdot)$ the ODE (\ref{eq:async:V:ODE}) has a unique globally asymptotically stable equilibrium $\lambda(\rho)$. The function $\lambda:H\rightarrow\mathbb{R}^{|\mathcal{S}|}$ is Lipschitz continuous. 
        
        \item[\textbf{(B3)}]\label{prop-B3} For every $\epsilon>0$, for all $\Omega^\epsilon_{|\mathcal{S}|\cdot |\mathcal{A}|}$-valued measurable functions $\Lambda^{(2)}(\cdot)$, the projected ODE (\ref{eq:async:rho:ODE}) has a unique globally asymptotically stable equilibrium $\rho^\star$.
        
        \item[\textbf{(B4)}]\label{prop-B4} $\sup_k(\|V_k\|+\|\rho_k\|)<\infty, \ \mathbb P\text{-}a.s.$ 
    
        \item[\textbf{(B5)}]\label{prop-B5} For $(s',a'),(s,a)\in\mathcal X$ and any $\rho\in H$, let $\mathcal{Q}_{(s,a),(s',a')}(\rho):=\mathcal{P}\left(s'|s,a\right)\,\pi^b_\rho(a'|s')$. Then
        \begin{enumerate}
            \item $\mathbb{P}\left((s_{k+1},a_{k+1})=(s',a')|\mathcal F_k\right)=\mathcal{Q}_{(s_k,a_k),(s',a')}(\rho_k)$;
            \item For each $\rho\in H$, the transition probabilities $\mathcal{Q}_{(s,a),(s',a')}(\rho)$ define an aperiodic, irreducible Markov chain over $\mathcal X$. 
            \item The mapping $\rho\mapsto\mathcal{Q}_{(s,a),(s',a')}(\rho)$ is Lipschitz continuous.
        \end{enumerate}
    \end{enumerate}
\end{proposition}
The proof is stated in Appendix~\ref{sec:async:PropProof}.

The properties~\hyperref[prop-B1]{(B1)}-a,b and~\hyperref[prop-B4]{(B4)} do not change compared to the respective assumptions for the synchronous scheme. Under property~\hyperref[prop-B1]{(B1)}-c, the bias term $\mathcal E_k$ vanishes and does not affect the limiting dynamics (see Remark 3 in~\cite[Section 2.2]{borkar2023stochastic}). Properties~\hyperref[prop-B2]{(B2)} and~\hyperref[prop-B4]{(B4)} are the asynchronous counterparts of~\hyperref[prop-A4]{(A4)} and~\hyperref[prop-A5]{(A5}. Property~\hyperref[prop-B5]{(B5)} ensures that the state-action index process $(s_k,a_k)$ is an ergodic Markov chain on $\mathcal X$, which implies that every state-action component of $\rho$ and every state component of $V$ is updated with positive asymptotic frequency (see also \cite[Lemma 3.4]{ref:perkins2013:asynchronous}).

\paragraph{Proof of Theorem~\ref{thm:Conv:Alg2}}
Let Assumptions~\ref{ass:ergodicity} and~\ref{ass:async:stepsize} hold, and let the stochastic gradients be as stated in~\eqref{eq:async:stochgrad}.
Given the properties listed in Proposition~\ref{prop:async:SA} and with the asynchronous step sizes defined in Assumption~\ref{ass:async:stepsize}, the fast iterates' interpolated trajectory converges to the trajectory of a solution of the differential inclusion~\eqref{eq:async:V:diff:incl}~\cite[Lemma 4.3]{ref:perkins2013:asynchronous}, and the slow iterate asymptotically tracks a solution to~\eqref{eq:async:rho:diff:incl} almost surely~\cite[Theorem 4.7]{ref:perkins2013:asynchronous}. 

Finally, the almost sure convergence of the last iterates of Algorithm~\ref{alg2} follows from the global asymptotic equilibrium property of the projected ODE~\eqref{eq:async:rho:ODE},~\hyperref[prop-B3]{(B3)}, by applying \citep[Corollary 4.8]{ref:perkins2013:asynchronous}.
\hfill $\blacksquare$

\subsection{Convergence Rate Analysis of Algorithm~\ref{alg2}}
\label{ssec: async-conv-rate}
In this section, we derive the convergence rate for the two-timescale asynchronous PGDA algorithm. 
For the derivation we impose one change compared to the stated Algorithm~\ref{alg2}, by projecting the primal variable iterates $V_k$ to the compact set $\mathcal V_r\subset \mathbb{R}_+^{|\mathcal{S}|}$, where 
\begin{equation}
\label{eq: primal-set-Vr}
    \mathcal{V}_r: = \{v\in\mathbb{R}_+^{|\mathcal{S}|} : v_i \leq \frac{C_r+\eta_\rho U_G}{1-\gamma} \ \forall i=1,\dots |\mathcal{S}|\}.
\end{equation}
As discussed in the proof of Corollary~\ref{cor: dual-reg}, this projection does not affect the optimal primal variable. 

To establish the almost sure convergence result in Theorem~\ref{thm:Conv:Alg2}, we assumed that for any strictly exploratory policy, the induced Markov chain is ergodic (Assumption~\ref{ass:ergodicity}). For the derivation of the convergence rate, we strengthen this assumption slightly as follows.

We denote Dobrushin's ergodic coefficient for a Markov kernel \(Q\) on \(\mathcal X\) as
\[
\delta(Q)
:=
\sup_{x,y\in\mathcal X}
\|\delta_x Q-\delta_y Q\|_{\mathrm{TV}}.
\]
It measures the maximal contraction of the kernel in total variation distance; see \citet{ref:dobrushin1956:central} and \citet[Definition 4.3.11]{ref:bremaud2020:markov}.

\begin{assumption}[Uniform mixing of dual-induced chains]
\label{ass:UGE:dual:policies}
Let \(\{\mathcal P_\rho\}_{\rho\in H}\) denote the \\family of state-action transition kernels
induced by the dual-induced policies \(\{\pi_\rho\}_{\rho\in H}\).
\begin{enumerate}
    \item[(i)] There exist constants \(C_{\mathcal X}>0\) and \(\varrho\in(0,1)\) such that,
    for every \(\rho\in H\), the Markov chain with kernel \(\mathcal P_\rho\) has a unique
    stationary law \(\mu_\rho\) and
    \[
    \sup_{x\in\mathcal X}
    \|\delta_x \mathcal P_\rho^k - \mu_\rho\|_{\mathrm{TV}}
    \le
    C_{\mathcal X}\varrho^k,
    \qquad k\ge 0.
    \]
    \item[(ii)] There exist constants \(m_\star\in\mathbb N\) and \(\kappa\in(0,1)\) such that, for every
    sequence \(\{\rho_t\}_{t\in\mathbb N}\subset H\),
    \[
    \sup_{k\ge 0}
    \delta\!\left(
    \prod_{t=k}^{k+m_\star-1}\mathcal P_{\rho_t}
    \right)
    \le
    \kappa.
    \]
\end{enumerate}
\end{assumption}

\begin{remark}
\label{rem:UGE:dual:policies}
    A sufficient condition for Assumption~\ref{ass:UGE:dual:policies} is the existence of a deterministic state-feedback policy \(\pi_d:\mathcal S\rightarrow\mathcal A\) whose induced state chain is uniformly geometrically ergodic, that is, it has a unique stationary law and converges to it geometrically fast in total variation, uniformly over initial states; see Lemma~\ref{lemma:UGE:dual:policies}.
\end{remark}

The convergence rate analysis presented here builds on the recent work of \cite{ref:zeng2024:two}. The analysis applies with a specific choice of the two-timescale stepsize sequences. In the following let the stepsize sequences be 
\begin{equation}\label{eq:async:stepsize:cvg_rate}
\alpha_k=\frac{\alpha_0}{(k+1)^{2/3}},\qquad \beta_k=\frac{\beta_0}{k+1},\qquad
\alpha_0,\,\beta_0>0.
\end{equation}
It can be easily verified that this choice satisfies Assumption~\ref{ass:async:stepsize}.

We extend and adapt the convergence rate analysis of \cite{ref:zeng2024:two} for our LP-based algorithm framework and asynchronous updating scheme in two important nontrivial aspects. 

First, while \cite{ref:zeng2024:two} considers only the global decaying stepsize setting, our focus is on the asynchronous regime with local stepsizes, which introduces randomness and non-monotonicity to the realized stepsize sequence. To handle the non-monotonicity, we switch the analysis to the currently least frequently visited state or state-action. This coordinate still receives updates linear in iteration counter $k$ eventually, with a high-probability bound on the time when the linear growth rate applies. The proof is based on the uniform geometric ergodicity of the dual policy induced chains and Lipschitz properties of the dual-induced Markov kernels.

Second, unlike \cite{ref:zeng2024:two}, the stochastic gradient estimator based on the replay buffer is biased in our setting. We find that the bias decays at a rate fast enough to not become dominating.

These properties of the asynchronous PGDA iterates as well as further regularity and boundedness properties form the basis of our convergence rate proof and are stated in the following proposition. 
Throughout, let $p_\star$ denote a uniform lower bound on the dual-variable induced stationary distributions over the state-action space,\footnote{An explicit expression for $p_\star$ is derived in the proof of Lemma~\ref{lemma: visitation-floor-envelopes}.}  valid for all $\rho \in H$.
Denote the minimal number of visits to any state-action pair until iteration k as $\underline{\nu}_k:=\min_{(s,a)\in\mathcal X}\nu_k(s,a)$. It holds that $\underline{\nu}_k\le \nu_k(s,a)\le k$ for all $(s,a)\in\mathcal X$. Denote the reduced Lagrangian objective as 
\begin{equation}
\label{eq: red-objective}
    f(\rho)=\min_V L(V,\rho).
\end{equation}
Recall that $\bar\alpha_k=\alpha(\tilde\nu_k(s_k))$ and $\bar\beta_k=\beta(\nu_k(X_k))$.

\begin{proposition}
\label{prop: collected-facts-A2}
There exist positive constants $\mu_{\mathrm{opt}},\,L$, $B$, $D$, and for each $\delta\in(0,1)$,
an event $\mathcal G_\delta$ with $\mathbb P(\mathcal G_\delta)\ge 1-\delta$ and an index $K(\delta)\in\mathbb N$ such that on $\mathcal G_\delta$ and for all $k\ge K(\delta)$,
\begin{itemize}
    \item[(i)] (Visitation floor) It holds that $\underline{\nu}_{k}\ge (p_\star/2)\,k$, and there exists a constant $c$ such that $c\cdot\alpha(\underline{\nu}_{k})\le\bar\alpha_k\le\alpha(\underline{\nu}_k),$ and $ c\cdot \beta(\underline{\nu}_{k})\le\bar\beta_k\le\beta(\underline{\nu}_k)$ with probability at least $1-\delta$ for $k\ge K(\delta)$.
    
    \item[(ii)] (Replay-buffer bias) There exists a constant $C_{\mathrm{buf}}(\delta)>0$, such that for all $\rho\in H$ the buffer-induced bias $\mathcal E_k(\rho)$, see~\eqref{eq:SA:def:bias}, decays as \[\|\mathbb{E}_{s\sim\tilde{\mu}_\rho}[e_s\,\mathcal E_k(\rho)(s)]\|_\infty \le C_{\mathrm{buf}}(\delta) \sqrt{\tfrac{1}{k}}\] with $C_{\mathrm{buf}}(\delta)=O(\log(1/\delta))$.
\end{itemize}
Moreover, the following regularity and boundedness properties hold.
\begin{itemize}
    \item[(iii)] (Regularity) $V\mapsto L(V,\rho)$ is $\eta_V$-strongly convex for all $\rho\in H$; 
    $\rho\mapsto\lambda(\rho)$ is Lipschitz on $H$; 
    $\rho\mapsto f(\rho)$ is $L$-smooth and $\mu_{\mathrm{opt}}$-strongly concave on $H$.
    
    \item[(iv)] (Uniform bounds) The stochastic gradients are uniformly bounded over $H\times\mathcal V_r$ for all $k\ge1$:
    $\sup_{\rho\in H, \, V\in\mathcal V_r}\|\widehat g_k(\rho,V;\cdot,\cdot)\|\le B$, $\sup_{\rho\in H, \, V\in\mathcal V_r}\|\widehat h_k(\rho,V;\cdot)\|\le  D$. 
\end{itemize}
\end{proposition}
The proof of Proposition~\ref{prop: collected-facts-A2} is provided in Section~\ref{ssec:Proof:convergence:rate}. 

Before stating the convergence rate theorem, we provide the necessary additional notation.
We focus the analysis on the least-visited state or state-action entry of the iterates. With the above stated visitation floor $\underline{\nu}_k\ge(p_\star/2)k$, we can bound the minimal visits to all state-action pairs and introduce the monotone decreasing upper envelope on the primal and dual stepsizes respectively as 
\begin{align}
\label{eq: stepsize-envelope}
    \alpha_{\mathrm{env},k}:=\alpha(\lfloor p_\star k/2\rfloor),\qquad \beta_{\mathrm{env},k}:=\beta(\lfloor p_\star k/2\rfloor).
\end{align}
Note that $\tfrac{\alpha_{\mathrm{env},k}}{\alpha_k}\le \tfrac{4}{p_\star}=:C_{\mathrm{env}}$. 
We define the mixing-time function as the worst-case mixing time across all admissible dual iterates. For $\zeta>0$
$$\tau(\zeta):=\sup_{\rho\in H}\min\{\ell\in\mathbb N:\sup_{x\in\mathcal X}\|\delta_x\mathcal P_\rho^\ell-\mu_\rho\|_\mathrm{TV}\le\zeta\}, \quad \tau_k:=\tau\big(\beta_{\mathrm{env},k}\big).$$ 
The uniform ergodicity of the Markov chains corresponding to the strictly positive dual policies now implies that there exists a constant $C_\tau>0$ such that $$\tau_k\le C_\tau \log(k+1)$$ and hence $\lim_{k\rightarrow\infty}\tau_k^2\alpha_k=0$. Consequently, there exists $\mathcal K\in\mathbb N$ such that for all $k\ge\mathcal K$ 
\begin{equation}
    \label{eq: K_ZDR}
    \tau_k^2\alpha_{\mathrm{env},k-\tau_k}\le\tfrac{\tilde \eta_V}{18\,C_1\,C_{\mathrm{env}}},\hspace{2pt} \tau_k\alpha_{\mathrm{env},k-\tau_k}\le\min\{\tfrac{1}{\sqrt{3}D},\,\tfrac{1}{\sqrt{3}LB}\}, \hspace{2pt} \alpha_k\le\tfrac{1}{\tilde\eta_V}, \hspace{2pt} \tfrac{p_\star k}{2}\ge 1,
\end{equation}
where 
\begin{equation}
    \label{eq: C1-eta_tilde}
    C_1:=2(LDB+LB+LD+D^2+\tfrac{\gamma |\mathcal X|C^U}{\beta_0}), \qquad\tilde\eta_V:=\eta_V\,p_\star|\mathcal A|.
\end{equation} 
Lastly, there exists $c_\tau\in(0,1)$, such that for all $k\ge\tau_k$, $c_\tau(k+1)\le k+1-\tau_k.$

With these definitions at hand, we can quantify the convergence rate of the asynchronous PGDA iterates as follows.

\begin{theorem}[Convergence rate of Algorithm~\ref{alg2}]\label{thm:B1}
Let Assumption~\ref{ass:UGE:dual:policies} hold.
Let $\{\rho_k\}$, $\{V_k\}$ denote the iterates of Algorithm~\ref{alg2} with the projection on the primal iterates to $\mathcal V_r$. Fix $\delta\in(0,1)$ and let $\mathcal G_\delta$, and $K(\delta)$ be as in Proposition~\ref{prop: collected-facts-A2}. 
Let the stepsize schedules be as in~\eqref{eq:async:stepsize:cvg_rate} with
\begin{equation*}
\beta_0\ge \frac{2}{p_\star\mu_{\mathrm{opt}}},\ \ \frac{\beta_0}{\alpha_0}\le \frac{\tilde\eta_V}{8p_\star\mu_{\mathrm{opt}}}.
\end{equation*}

Then, on $\mathcal G_\delta$, for all $k\ge k_0:=\max\{K(\delta),\mathcal K\}$,
\begin{align*}
    \mathbb E\!\left[f(\rho^\star)-f(\rho_k)| \mathcal G_\delta\right] &\le \left(\mathbb E\!\left[f(\rho^\star)-f(\rho_{k_0})| \mathcal G_\delta\right]+\frac{16 L^2}{p_\star^3\tilde\eta_V}\frac{\beta_{k_0}}{\alpha_{k_0}}\mathbb E\!\left[\|V_{k_0}-\lambda(\rho_{k_0})\|_2^2|\mathcal G_\delta\right]\right)\frac{k_0+1}{k+1}\\
&~~~~~~~~~~~~~~~~~~~~~~~~~~~~~~~~+ \frac{2C_2(\delta)\log^2(k+1)}{3(k+1)^{2/3}},
\end{align*}
where $L$ is a Lipschitz parameter and $C_2(\delta)$ is a constant depending only on model parameters and $\delta$ but not on $k$.
Moreover,
\[
\mathbb E\!\left[\|\rho_k-\rho^\star\|^2_2| \mathcal G_\delta\right]
\le \tfrac{2}{\mu_{\mathrm{opt}}}\,\mathbb E\!\left[f(\rho^\star)-f(\rho_k)| \mathcal G_\delta\right]=\tilde{\mathcal{O}}(k^{-2/3}).
\]

\end{theorem}

The rate of convergence of the order $\widetilde{\mathcal{O}}(k^{-2/3})$ recovers the existing rate of \cite{ref:zeng2024:two} for two-timescale stochastic gradient schemes under Markovian model access, albeit under the more challenging setting of buffer-based biased gradient estimates and accounting for asynchronous updating with visitation-adjusted stepsizes. 

As a consequence of the finite-time convergence rate of the dual iterates (Theorem~\ref{thm:B1}), we obtain the following bound on the unregularized value suboptimality of the dual-induced policies; see Appendix~\ref{ssec:Aux_results:convergence:rate} for the proof.
\begin{corollary}[Value suboptimality of the induced policy]\hspace{2pt}
\label{cor:value-suboptimality}
Under the assumptions of Theorem~\ref{thm:B1}, let \(L_\pi>0\) satisfy
\[
\max_{s\in\mathcal S}\|\pi_\rho(\cdot|s)-\pi_{\hat\rho}(\cdot|s)\|_{1}\le L_\pi\|\rho-\hat\rho\|_2,
\qquad \forall\,\rho,\hat\rho\in H,
\]
and define
\[
L_V:=\frac{C_r}{(1-\gamma)^2}L_\pi,
\qquad
B_{\mathrm{reg}}:=\frac{\eta_\rho\log|\mathcal A|}{1-\gamma}.
\]
Then, for every \(s\in\mathcal S\) and all \(k\ge k_0\),
\[
\mathbb E\!\left[
\left(
V^\star_{ur}(s)-V^{\pi_{\rho_k}}_{ur}(s)-B_{\mathrm{reg}}
\right)_+^2
\;\middle|\;
\mathcal G_\delta
\right]
\le
L_V^2\,
\mathbb E\!\left[
\|\rho_k-\rho^\star\|_2^2
\;\middle|\;
\mathcal G_\delta
\right].
\]
In particular, the excess unregularized value suboptimality beyond the deterministic
regularization bias \(B_{\mathrm{reg}}\) inherits the conditional mean-square rate
\(\widetilde{\mathcal O}(k^{-2/3})\) from Theorem~\ref{thm:B1}.
\end{corollary}

The proof of Theorem~\ref{thm:B1} builds on two contraction inequalities, one each for the primal and the dual variable, which are then connected by the two-timescale lemma~\cite[Lemma 4]{ref:zeng2024:two}. Recall that we overload the Lipschitz constant notation and denote all Lipschitz constants by the envelope $L$.

\begin{proposition}
    \label{prop: primal-tracking}
    For $\delta\in(0,1)$ let $\mathcal G_\delta$ denote the high-probability visitation event of Proposition~\ref{prop: collected-facts-A2}. Let $C_1$ be as in~\eqref{eq: C1-eta_tilde}. Then, for all $k\ge \max\{K(\delta),\mathcal K\}$, the primal variable tracking error $\mathbb E\left[\|V_k-\lambda(\rho_k)\|_2^2|\mathcal G_\delta\right]$ satisfies
    \begin{align*}
        \mathbb E\left[\|V_{k+1}-\lambda(\rho_{k+1})\|_2^2|\mathcal G_\delta\right]\;\le&\;(1-\tfrac{\tilde\eta_V}{4}\alpha_{k})\mathbb E\left[\|V_k-\lambda(\rho_k)\|_2^2|\mathcal G_\delta\right]+r_k(\delta),
    \end{align*}
    where 
    $$r_k(\delta):=\tfrac{9}{2}C_1\tau_k^2\alpha_{\mathrm{env},k-\tau_k}\alpha_{\mathrm{env},k}+\tfrac{3}{2}\tfrac{C_{\mathrm{buf}}(\delta)^2|\mathcal S|}{\tilde\eta_V}\tfrac{\alpha_{\mathrm{env},k}}{k}+\tfrac{3}{2}D^2\alpha_{\mathrm{env},k}^2+L^2B^2\beta_{\mathrm{env},k}^2+\tfrac{2L^2B^2}{\tilde\eta_V}\tfrac{\beta_{\mathrm{env},k}^2}{\alpha_k}.$$ 
    \end{proposition}

\begin{proposition}
    \label{prop: dual-variable-error}
     For $\delta\in(0,1)$ let $\mathcal G_\delta$ denote the high-probability visitation event of Proposition~\ref{prop: collected-facts-A2}. Then for all $k\ge\max\{K(\delta),\,\mathcal K\}$, the reduced objective error $\mathbb E[f(\rho^\star)-f(\rho_{k})|\mathcal G_\delta]$ satisfies 
     \begin{align*}
        \mathbb E[f(\rho^\star)-f(\rho_{k+1})|\mathcal G_\delta]&\le(1-p_\star\mu_\mathrm{opt}\,\beta_k)\mathbb E[f(\rho^\star)-f(\rho_{k})|\mathcal G_\delta]+\frac{L^2\beta_{\mathrm{env},k}}{2p_\star}\mathbb E[\|V_k-\lambda(\rho_k)\|_2^2|\mathcal G_\delta]\\
        &~~~~~~~~~~+\frac{25 L^2B^3}{2}\tau_k^2\beta_{\mathrm{env},k-\tau_k}\beta_{\mathrm{env},k}, 
    \end{align*}
    where $\mu_\mathrm{opt}>0$ is the strong concavity modulus of the reduced objective $f$.
\end{proposition}
The proofs of the two propositions are stated in Appendix~\ref{ssec: contraction-proofs}.

The proof of Theorem~\ref{thm:B1} now consists of combining the recursive contractions from Propositions~\ref{prop: primal-tracking} and \ref{prop: dual-variable-error} using the two-timescale lemma~\cite[Lemma 4, Case 1]{ref:zeng2024:two}, restated here for completeness.
\begin{lemma}[Two-timescale lemma,~\citealt{ref:zeng2024:two}]\label{lemma: two-timescaleZDR}
    Let $\left\{a_k, b_k, c_k, d_k, e_k, f_k\right\}$ be\\ non-negative sequences satisfying $\frac{a_{k+1}}{d_{k+1}} \leq \frac{a_k}{d_k}<1$ for all $k \geq 0$. Let $\left\{x_k\right\},\left\{y_k\right\}$ be two nonnegative sequences. Suppose that $x_k, y_k$ satisfy the following coupled inequalities:
        $$
        x_{k+1} \leq\left(1-a_k\right) x_k+b_k y_k+c_k, \quad y_{k+1} \leq\left(1-d_k\right) y_k+e_k x_k+f_k .
        $$
    In addition, assume that there exists a constant $A \in \mathbb{R}$ such that
        $$
        A a_k\big(1-\frac{a_k}{d_k}\big)- b_k\ge 0 \quad \text { and } \quad \frac{A e_k}{d_k} \leq \frac{1}{2} \quad \forall k \geq 0 .
        $$
    Then we have for all $0 \leq \tau \leq k$
        $$
        x_k \leq\left(x_\tau+\frac{A a_\tau}{d_\tau} y_\tau\right) \prod_{t=\tau}^{k-1}\left(1-\frac{a_t}{2}\right)+\sum_{\ell=\tau}^{k-1}\left(c_{\ell}+\frac{A a_{\ell} f_{\ell}}{d_{\ell}}\right) \prod_{t=\ell+1}^{k-1}\left(1-\frac{a_t}{2}\right) .
        $$
\end{lemma}

\paragraph{Proof of Theorem~\ref{thm:B1}}

Fix $\delta\in(0,1)$ and work on $\mathcal G_\delta$ defined in Proposition~\ref{prop: collected-facts-A2}.
Let $k_0:=\max\{K(\delta),\mathcal K\}$ be the burn-in index such that the visitation floor bounds and the technical conditions in~\eqref{eq: K_ZDR} hold for all $k\ge k_0$.
Set
\[
x_k := \mathbb E\left[f(\rho^\star)-f(\rho_k)| \mathcal G_\delta\right],
\qquad
y_k := \mathbb E\left[\|V_k-\lambda(\rho_k)\|_2^2| \mathcal G_\delta\right].
\]
The proof proceeds in four steps.

\emph{(I.)} Following the two coupled inequalities of Proposition~\ref{prop: primal-tracking} and Proposition~\ref{prop: dual-variable-error} it holds that for all $k\ge k_0$
\begin{equation}
\label{eq: x-rec-ZDR}
x_{k+1}\le (1-a_k)x_k + b_k y_k + c_k,
\end{equation}
where
\[
a_k:=p_\star\mu_{\mathrm{opt}}\beta_k,\qquad
b_k:=\frac{L^2}{2p_\star}\beta_{\mathrm{env},k},\qquad
c_k:=\frac{25L^2B^3}{2}\,\tau_k^2\,\beta_{\mathrm{env},k-\tau_k}\beta_{\mathrm{env},k},
\]
and
\begin{equation}
\label{eq: y-rec-ZDR}
y_{k+1}\le (1-d_k)y_k + f_k,
\end{equation}
where
\[
d_k:=\frac{\tilde\eta_V}{4}\alpha_k,\qquad e_k:=0,\qquad f_k:=r_k(\delta).
\]
Inserting the stepsize sequences into $r_k(\delta)$ (see Proposition~\ref{prop: primal-tracking}) and adding the constant factors, we find that, with constant 
$C_R(\delta):=\tfrac{12}{p_\star\tilde\eta_V}C_\mathrm{buf}(\delta)^2|\mathcal S|\alpha_0+\tfrac{12}{p_\star^2}\alpha_0^2(\tfrac{6C_1C_\tau^2}{c_\tau}+2D^2)+\tfrac{16}{p_\star^2}L^2B^2\beta_0^2(1+\tfrac{2}{\alpha_0\tilde\eta_V})$, for all $k\ge k_0$,
\begin{equation}
\label{eq: f-bound}
f_k \le \frac{C_R(\delta)\,\log^2(k+1)}{(k+1)^{4/3}}.
\end{equation}

\emph{(II.)} We proceed to determine a feasible constant $A\in\mathbb R$ to apply Lemma~\ref{lemma: two-timescaleZDR}. Let $k\ge k_0$. Inserting the chosen stepsize scheme~\eqref{eq:async:stepsize:cvg_rate} into the fraction $a_k/d_k$ gives
\begin{equation}
\label{eq: ratio-ad}
\frac{a_k}{d_k}
=\frac{p_\star\mu_{\mathrm{opt}}\beta_k}{(\tilde\eta_V/4)\alpha_k}
=\frac{4p_\star\mu_{\mathrm{opt}}}{\tilde\eta_V}\frac{\beta_k}{\alpha_k}=\frac{4p_\star\mu_{\mathrm{opt}}}{\tilde\eta_V}\frac{\beta_0}{\alpha_0(k+1)^{1/3}}.
\end{equation}
By the imposed conditions on the stepsize parameters, $\beta_0\ge \tfrac{2}{p_\star\mu_{\mathrm{opt}}}$ and $\alpha_0$ is set such that $\frac{\beta_0}{\alpha_0}\le \frac{\tilde\eta_V}{8p_\star\mu_{\mathrm{opt}}}$.
Since $\tfrac{\beta_k}{\alpha_k} = \tfrac{\beta_0}{\alpha_0(k+1)^{1/3}}$ is non-increasing in $k$, this implies $\tfrac{a_k}{d_k}\le\tfrac{1}{2}(k+1)^{-1/3}\le \tfrac{1}{2}$.
Moreover, by the stepsize envelope construction (see~\ref{eq: stepsize-envelope}), 
$$\beta_{\mathrm{env},k}\le\beta(\lfloor\tfrac{p_\star k}{2}\rfloor)=\tfrac{\beta_0}{\lfloor\tfrac{p_\star k}{2}\rfloor+1}\le \tfrac{4}{p_\star}\,\beta_k.$$
We set $A := \frac{4L^2}{p_\star^3\mu_{\mathrm{opt}}}.$
Then, using $\tfrac{a_k}{d_k}\le\tfrac{1}{2}$ and $b_k=\tfrac{L^2}{2p_\star}\beta_{\mathrm{env},k}\le \tfrac{2L^2}{p_\star^2}\beta_k$, we obtain 
\[
A a_k\Bigl(1-\frac{a_k}{d_k}\Bigr)-b_k
\ge \frac{A}{2}a_k - b_k
\ge \frac{2L^2}{p_\star^2\mu_{\mathrm{opt}}}\mu_{\mathrm{opt}}\beta_k - \frac{2L^2}{p_\star^2}\beta_k
=0,
\]
and $A e_k/d_k=0$ since $e_k=0$. Hence, Lemma~\ref{lemma: two-timescaleZDR} applies for the coupled system~\eqref{eq: x-rec-ZDR}--\eqref{eq: y-rec-ZDR} and gives
\begin{equation}
\label{eq: combined-x-ineq}
    x_k \leq\left(x_{k_0}+A \frac{a_{k_0}}{d_{k_0}} y_{k_0}\right) \prod_{t={k_0}}^{k-1}\left(1-\frac{a_t}{2}\right)+\sum_{\ell={k_0}}^{k-1}\left[\prod_{t=\ell+1}^{k-1}\left(1-\frac{a_t}{2}\right)\left(c_{\ell}+A \frac{a_{\ell}}{d_{\ell}} f_{\ell}\right)\right].
\end{equation}

\emph{(III.)}
To derive the convergence rate result, it remains to bound the terms of the above inequality. First, using $\tau_k\le C_\tau\log(k+1)$ and $\beta_{\mathrm{env},k-\tau_k}\le \tfrac{4}{c_\tau p_\star}\beta_k$,
\begin{equation*}
c_k \le \frac{25L^2B^3}{2c_\tau}\,(C_\tau^2\log^2(k+1))\,(\tfrac{4}{p_\star}\beta_k)^2= \frac{200L^2B^3C_\tau^2\beta_0^2}{c_\tau p_\star^2}\cdot\frac{\log^2(k+1)}{(k+1)^2}.
\end{equation*}
Second, combine~\eqref{eq: f-bound} and~\eqref{eq: ratio-ad} to get
\begin{equation*}
A\frac{a_k}{d_k}f_k
\le \frac{2L^2C_R(\delta)}{p_\star^3\mu_{\mathrm{opt}}}\cdot \frac{\log^2(k+1)}{(k+1)^{5/3}}.
\end{equation*}
Therefore, the multiplier in the second term of~\eqref{eq: combined-x-ineq} satisfies
\begin{equation}\label{eq: gk-bound}
c_k + A\frac{a_k}{d_k}f_k
\le
C_2(\delta)\frac{\log^2(k+1)}{(k+1)^{5/3}},
\qquad
C_2(\delta):=\frac{2L^2C_R(\delta)}{p_\star^3\mu_{\mathrm{opt}}}+\frac{200L^2B^3C_\tau^2\beta_0^2}{c_\tau p_\star^2},
\end{equation}
since $(k+1)^{-2}\le (k+1)^{-5/3}$ for all $k\ge 1$.

\emph{(IV.)}
Using $1-u\le e^{-u}$ and the harmonic-sum bound $\sum_{t=k_1}^{k_2}\tfrac{1}{t+1}\ge\log(\tfrac{k_2+2}{k_1+1})$, for all $k_0\le \ell<k$,
\begin{align*}
\prod_{t=\ell+1}^{k-1}\Bigl(1-\frac{a_t}{2}\Bigr)
&\le \exp\Bigl(-\frac12\sum_{t=\ell+1}^{k-1}a_t\Bigr)
= \exp\Bigl(-\frac{p_\star\mu_{\mathrm{opt}}\beta_0}{2}\sum_{t=\ell+1}^{k-1}\frac{1}{t+1}\Bigr)\\
&\le \Bigl(\frac{\ell+2}{k+1}\Bigr)^{\tfrac{p_\star\mu_{\mathrm{opt}}\beta_0}{2}}\le\frac{\ell+2}{k+1},
\end{align*}
where the last inequality is due to the condition $\beta_0\ge\tfrac{2}{p_\star\mu_{\mathrm{opt}}}$. With the same reasoning $\prod_{t={k_0}}^{k-1}\left(1-\frac{a_t}{2}\right)\le\tfrac{k_0+1}{k+1}$. 
Inserting the harmonic sum bounds and~\eqref{eq: gk-bound} in~\eqref{eq: combined-x-ineq} and using $\frac{\ell+2}{k+1}\le2\tfrac{\ell+1}{k+1}$ yields, for all $k\ge k_0$,
\[
x_k \le \left(x_{k_0}+\frac{16L^2}{p_\star^3\tilde\eta_V}\frac{\beta_{k_0}}{\alpha_{k_0}}y_{k_0}\right)\frac{k_0+1}{k+1}
+ \frac{2C_2(\delta)\log^2(k+1)}{k+1}\sum_{\ell=k_0}^{k-1}\frac{1}{(\ell+1)^{2/3}},
\]
where for the second summand, we use the monotonicity $\log^2(\ell+1)\le\log^2(k+1)$ for $\ell\le k$.
We apply an integral bound to the sum. Since $y\mapsto y^{-2/3}$ is decreasing on $(0,\infty)$, we have $(\ell+1)^{-2/3}\le \int_{\ell}^{\ell+1} y^{-2/3}dy$. Hence, $\sum_{\ell=k_0}^{k-1}(\ell+1)^{-2/3}\le \int_{k_0}^{k} y^{-2/3}dy = 3\big(k^{1/3}-k_0^{1/3}\big)\le 3(k+1)^{1/3}.$ Inserting the integral bound finally yields the claimed rate
\begin{align*}
\mathbb E\left[f(\rho^\star)-f(\rho_k)| \mathcal G_\delta\right] 
\le 
\big(
\mathbb E\left[f(\rho^\star)-f(\rho_{k_0})| \mathcal G_\delta\right]
&+\frac{16 L^2}{p_\star^3\tilde\eta_V}\frac{\beta_{k_0}}{\alpha_{k_0}}
\mathbb E\left[\|V_{k_0}-\lambda(\rho_{k_0})\|_2^2|\mathcal G_\delta\right]
\big)
\frac{k_0+1}{k+1}\\
&~~~+ \frac{6C_2(\delta)\log^2(k+1)}{(k+1)^{2/3}}.
\end{align*}
Since $\rho^\star\in \mathrm{int}H$, it holds that $\nabla f(\rho^\star)=0$, hence with the strong concavity of $f$ (see Proposition~\ref{prop: collected-facts-A2}-(iv) and the proof in Lemma~\ref{lemma: strong-concavity-f}) 
\begin{align*}
    f(\rho_k)&\le f(\rho^\star)+\langle\nabla f(\rho^\star), \rho_k-\rho^\star\rangle-\tfrac{\mu_{\mathrm{opt}}}{2}\|\rho_k-\rho^\star\|_2^2\\
    &\Rightarrow \|\rho_k-\rho^\star\|_2^2\le \tfrac{2}{\mu_{\mathrm{opt}}}\left( f(\rho^\star)-f(\rho_k)\right),
\end{align*}
\hfill $\blacksquare$


\section{Numerical Results}\label{sec:numerical:results}

We evaluate PGDA-RL on the FrozenLake-v1 environment~\citep{ref:towers2024gymnasium}, a standard reinforcement learning benchmark with discrete states and actions.\footnote{All experiments are implemented in Python, and the code for reproducing all numerical results is available from \url{https://github.com/AxelFW/tt-primal-dual-rl}.} 
The environment consists of a $4\times 4$ grid in which the agent must reach a goal state while avoiding holes. The action space contains four actions (up, down, left, right) and transitions are stochastic due to the ``slippery'' dynamics. To adapt this episodic task to our infinite-horizon setting, we treat reaching the goal or a hole as a transition back to the start state, allowing continuous learning. The reward for reaching the goal is $100$, and all other transitions yield zero reward.

We set the model parameters to $\{\gamma,\eta_\rho,\eta_V\}=\{0.9,0.1,0.1\}$. The asynchronous stepsizes are
$\alpha_k=(1+\tilde k)^{-2/3}$ and $\beta_k=(1+\tilde k)^{-1}$ with $\tilde k=9+k/100$.
We apply Algorithm~\ref{alg2} with on-policy exploration and an additional $\epsilon$-greedy strategy with linearly decaying $\epsilon$ ($[\epsilon_0,\epsilon_K]=[1,0.1]$, $K=10^5$) to mitigate infrequent updates of rarely visited state-action pairs. Replay buffers are capped at size $10^3$. In the numerical implementation, we additionally impose a lower bound of $10^{-12}$ on the entries of the dual iterate after projection to avoid numerical instabilities when projected components become extremely small.

We track value-function and policy errors throughout training. In particular, the reported relative root-MSE (rRMSE) is the relative $\ell_2$ error
$\|V-\bar V\|_2/\|\bar V\|_2$ evaluated on non-terminal states, where $\bar V$ denotes the corresponding reference value.
\begin{figure}[t]
\captionsetup{justification=justified,singlelinecheck=false}
    \centering
    \includegraphics[width=\textwidth]{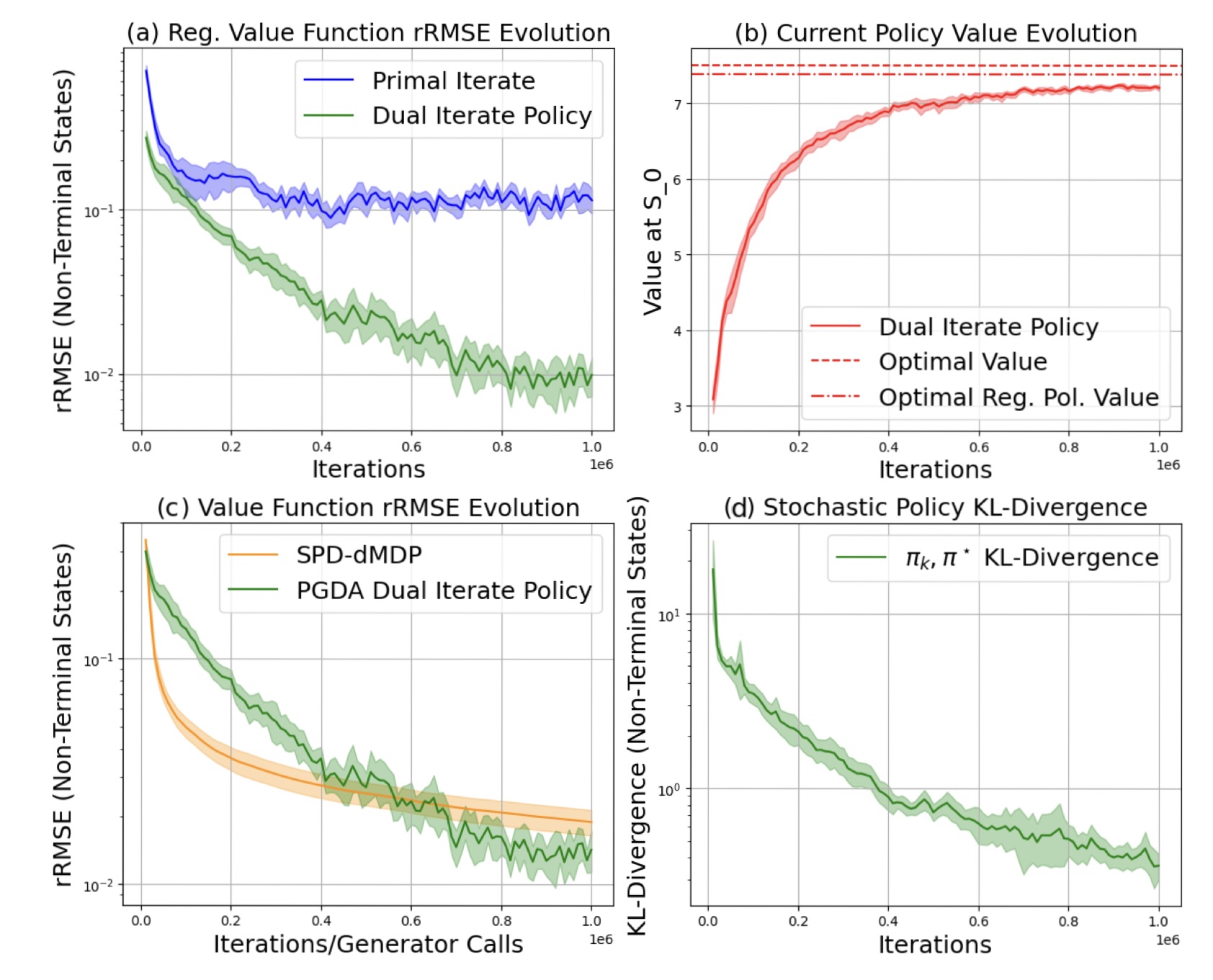}
    \caption{Value function and policy evaluation during training. The plots show the mean (solid lines) and shaded $\pm 2$ standard errors over $10$ independent training runs. 
    (a) rRMSE to the optimal regularized value. 
    (b) Initial-state value under the dual-iterate policy. 
    (c) rRMSE to the optimal (unregularized) value and comparison to the generator-based benchmark \textit{SPD-dMDP}~\citep{chen2016}. 
    (d) KL-divergence between the learned policy $\pi_{\rho_k}$ and the optimal regularized policy $\pi_r^\star$.}
    \label{fig:training_results}
\end{figure}

Figure~\ref{fig:training_results} summarizes the learning dynamics. Panel (a) shows that both the primal value iterate and the value induced by the current dual policy steadily reduce their relative value error on non-terminal states. The dual-policy value improves faster and reaches a lower error level over the shown horizon. 
This separation is consistent with the two-timescale design: while $\rho_k$ evolves slowly toward $\rho^\star$, the value update primarily tracks the moving best response $\lambda(\rho_k)$, so the policy-induced value can stabilize earlier than the primal iterate over the finite training horizon.
Panel (b) tracks the initial-state value under the dual-iterate policy: after a rapid increase during early training, the curve approaches the initial state value of the optimal regularized policy. The separation between the optimal unregularized value and the value of the optimal regularized policy shows the suboptimality induced by the entropy regularization. 
Panel (c) compares PGDA-RL against the generator-based benchmark SPD-dMDP~\citep{chen2016} that builds on the LP approach to (unregularized) MDPs. SPD-dMDP improves more quickly early on, whereas PGDA-RL closes the gap over longer training and achieves comparable final accuracy over the shown horizon. A key difference is the sampling model: SPD-dMDP draws transitions from uniformly sampled state-action pairs, while our method relies on trajectory data generated by the current behavior policy. Moreover, SPD-dMDP reports the averaged iterate as output, which yields a visibly smoother curve than the last-iterate behavior shown for PGDA-RL. 
Finally, panel (d) reports the KL divergence between the learned dual iterate policy and the optimal regularized policy summed uniformly over non-terminal states, $\mathsf{D}(\pi_r^\star \parallel \pi_{\rho_k}) = 
\sum_{s \in \mathcal{S},a\in\mathcal{A}}\pi_r^\star(a|s) \log \left( \frac{\pi_r^\star(a|s)}{\pi_{\rho_k}(a|s)} \right)$, which decreases consistently over training, indicating convergence in policy space alongside value improvement.
In summary, the simulations support our theoretical findings by demonstrating convergence of the dual-induced policy $\pi_{\rho_k}$ toward the unique optimal regularized policy $\pi_r^\star$ on all non-terminal states.
\section{Concluding Reflections}\label{sec:conclusion}
We proposed a new LP-based primal-dual RL algorithm for entropy-regularized MDPs, suing a single on-policy trajectory generated by a behavioral policy that may evolve with the dual iterate. Using the dynamical-systems approach to stochastic approximation \citep{borkar2023stochastic}, we establish almost sure last-iterate convergence of the primal and dual iterates to the optimal regularized value function and policy, respectively, under an asynchronous two-timescale scheme with diminishing stepsizes. To the best of our knowledge, such almost sure guarantees have not previously been established for LP-based methods in the entropy-regularized setting.

Under a strengthened ergodicity condition, we additionally derive a finite-time last-iterate guarantee, showing that the expected dual objective gap, and consequently the dual mean-squared error, decays at rate $\widetilde{\mathcal O}(k^{-2/3})$. This matches the best-known last-iterate behavior for two-timescale stochastic approximation under Markovian sampling~\citep{ref:zeng2024:two}, despite replay-buffer bias and asynchronous, visitation-adjusted updates along a single trajectory.

Our finite-time analysis relies on a strengthened ergodicity condition, namely the existence of a deterministic ergodic policy, and a high-probability visitation event ensuring sufficient exploration along the trajectory. A natural future direction is to weaken these requirements and to derive last-iterate rates under milder mixing conditions, potentially by combining refined concentration tools for inhomogeneous Markov chains with adaptive exploration schedules. Moreover, while our guarantees are stated in terms of the dual objective gap (and the induced mean-square error), it would be valuable to connect primalreplay-buffer biasdual convergence more directly to control-centric criteria such as (regularized) return suboptimality and constraint-violation measures, thereby enabling sharper comparisons with alternative performance guarantees in online RL.

While this work focuses on the tabular setting, extending PGDA-RL to function approximation remains a pressing challenge. One promising route is to combine our primal-dual viewpoint and stochastic approximation analysis with LP-based frameworks under linear function approximation, as studied for offline settings by \citet{ref:gabbianelli2024:offline}, and to develop corresponding finite-time guarantees under Markovian single-trajectory data. Such extensions will likely require reformulating the asynchronous update scheme and the structured experience replay mechanism to accommodate approximation error and stability issues. Finally, it is an interesting open question to integrate acceleration ideas for convexified regularized MDP formulations, which have been shown to speed up convergence under full model knowledge \citep{li2024accelerating}, into the present Markovian, replay-based setting.

\acks{The authors were supported by the German Research Foundation (DFG) under the grant ``SU 1433/2-1".}

\appendix
\section{Synchronous Almost Sure Convergence: Proof Details}
\label{sec: synchronous-alg-conv-proof}
This section contains the proofs of Propositions~\ref{prop: sync:gradient:noise} and \ref{prop: ODEs} that are at the heart of the almost sure convergence proof of Algorithm~\ref{alg1}.

\paragraph{Proof of Proposition~\ref{prop: sync:gradient:noise}}
\label{sec: GradMart}
The true gradients of the regularized Lagrangian~\eqref{eq: regLP} are given by
\begin{subequations}\label{eq: gradients}
    \begin{align}   
        \label{eq: Vgrad}
        \nabla_{V}L(V,\rho)({s'})=&\eta_VV(s')+\tilde\rho(s')+\gamma\sum_{(s,a)\in\mathcal{S}\times\mathcal{A}}\rho(s,a)\mathcal{P}(s' | s, a), \quad s'\in\mathcal{S} \\
        \label{eq: rgrad}
        \nabla_{\rho}L(V,\rho)(s,a)=&-V(s)+r(s, a)+\gamma \sum_{s'\in\mathcal{S}} V\left(s^{\prime}\right) \mathcal{P}(s'|s,a)-\eta_\rho\log\left(\frac{\rho(s,a)}{\tilde{\rho}(s)}\right), 
        \end{align}
\end{subequations}
where $(s,a)\in\mathcal{S}\times\mathcal{A}.$

The Lagrangian's gradient with respect to $V$,~\eqref{eq: Vgrad}, is linear in $V$ and therefore Lipschitz continuous. Continuously differentiable functions are Lipschitz continuous on compact sets. Hence, the projection of the dual variable to the compact hyper-rectangle $H$ ensures the Lipschitz continuity of the gradient~\eqref{eq: rgrad} and property~\hyperref[prop-A1]{(A1)}
The noise terms resulting from the stochastic gradients~\eqref{eq:sync:stochgrad} with the iterates $\rho_{k-1},\ V_{k-1}$ inserted at iteration $k$ are
\begin{subequations}\label{eq:sync:noise}
    \begin{align}
        \label{eq: noise1}
        M^{(1)}_k(s')&=\gamma\sum_{(s,a)\in\mathcal S\times\mathcal A}\rho_{k-1}(s,a)\left(\mathbf 1\{\tilde{s}(s,a)=s'\}-\mathcal{P}(s' | s,a)\right), \quad s'\in\mathcal{S}\\
        \label{eq: noise2}
        M^{(2)}_k(s,a)&=\gamma\left(V_{k-1}\left(\tilde{s}(s,a)\right)-\sum_{s'\in\mathcal{S}} V_{k-1}(s^{\prime}) \mathcal{P}(s'|s,a)\right),\ (s,a)\in\mathcal{S}\times \mathcal{A}.
        \end{align}
\end{subequations}
Recall the definition of the filtration $\mathcal F_k:=\sigma(V_\ell,\rho_\ell,M^{(1)}_\ell,M^{(2)}_\ell,\ell\le k),\ k\ge 0$.
The martingale difference properties of the noise terms with respect to the filtration $\{\mathcal F_k\}$ are now a direct consequence of the access to the generative model to obtain the state transitions $\tilde s(s,a)\sim\mathcal{P}(\cdot|s,a)$ for all state-action pairs in each update step. 

We proceed with the bounds on the second moment and begin with $M_k^{(1)}$
\begin{align*}
    \mathbb{E}\left[\|M^{(1)}_k\|_2^2| \mathcal{F}_{k-1}\right]&=\mathbb{E}\left[\sum_{s'\in\mathcal{S}}(M^{(1)}_k(s'))^2| \mathcal{F}_{k-1}\right]\\
    =\gamma^2&\sum_{s'\in\mathcal{S}}\mathbb{E}\left[\left(\sum_{(s,a)\in\mathcal S\times\mathcal A}\rho_{k-1}(s,a)\left(\mathbf 1\{s'=\tilde{s}(s,a)\}-\mathcal{P}(s'| s,a)\right)\right)^2|\mathcal{F}_{k-1}\right]\\
    \leq\gamma^2&|\mathcal{S}|\mathbb{E}\left[\left(\sum_{(s,a)\in\mathcal S\times\mathcal A}\rho_{k-1}(s,a)\right)^2|\mathcal{F}_{k-1}\right]\\
    =\gamma^2&|\mathcal{S}|\left(\sum_{(s,a)\in\mathcal S\times\mathcal A}\rho_{k-1}(s,a)\right)^2 \\
    =\gamma^2&|\mathcal{S}|\|\rho_{k-1}\|_1^2\ \leq \ \gamma^2|\mathcal{S}|^2|\mathcal{A}|\|\rho_{k-1}\|_2^2.
\end{align*}
The first inequality holds since $\mathbf 1\{s'=\tilde{s}(s,a)\}-\mathcal{P}(s'| s,a)\in[-1,1]$. 
Now consider $M_k^{(2)}$

\begin{align*}
    \mathbb{E}\left[\|M^{(2)}_k\|_2^2| \mathcal{F}_{k-1}\right]&=\mathbb{E}\left[\sum_{(s,a)\in\mathcal S\times\mathcal A}(M^{(2)}_k(s,a))^2| \mathcal{F}_{k-1}\right]\\
    &=\gamma^2\ \mathbb{E}\left[\sum_{(s,a)\in\mathcal S\times\mathcal A}\left(V_{k-1}(\tilde{s}(s,a))-\sum_{s'\in\mathcal{S}} V_{k-1}(s^{\prime}) \mathcal{P}(s'|s,a)\right)^2| \mathcal{F}_{k-1}\right]\\
    &\leq\gamma^2\sum_{(s,a)\in\mathcal S\times\mathcal A}\mathbb{E}\left[V_{k-1}(\tilde{s}(s,a))^2|\mathcal{F}_{k-1}\right]\\
    &=\gamma^2\sum_{(s,a)\in\mathcal S\times\mathcal A}\sum_{s'\in\mathcal{S}}\mathcal{P}(s'|s,a)V_{k-1}(s')^2\\
    &\leq\gamma^2|\mathcal S||\mathcal A|\|V_{k-1}\|_2^2,
\end{align*}
which completes the proof of property~\hyperref[prop-A2]{(A2)}. 
\hfill $\blacksquare$

\paragraph{Proof of Proposition~\ref{prop: ODEs}}
We start with the proof of the global asymptotic stability property~\hyperref[prop-A3]{(A3)} for the primal-dual gradient flow.

Let $\rho\in H$ be fixed. Recall that due to the strict convex-concave structure of the regularized Lagrangian, there exists a unique saddle-point $(V^\star,\rho^\star)$. With the first-order optimality condition based on~\eqref{eq: Vgrad}, the best response $\lambda(\rho)$ is
$$
\begin{aligned}
\nabla_V L(V(t), \rho)(s')=\eta_V V(t)(s^{\prime})-\tilde\rho(s')+\gamma\sum_{(s,a)\in\mathcal{S}\times\mathcal{A}} \rho(s, a)\mathcal{P}(s^{\prime} | s, a) \stackrel{!}{=} 0 \\
\Leftrightarrow \lambda(\rho)(s^{\prime})=\frac{1}{\eta_V}\left(\tilde\rho(s')- \gamma \sum_{(s,a)\in\mathcal{S}\times\mathcal{A}} \rho(s, a)\mathcal{P}(s^{\prime} | s, a)\right).
\end{aligned}
$$
Since $\tilde \rho(s')=\sum_{a\in\mathcal A}\rho(s',a)$, the mapping $\lambda$ is linear in $\rho$, hence Lipschitz continuous. 

Next, we establish the global asymptotic stability by providing a strict Lyapunov function, $\varphi_\rho:\mathbb R^{|\mathcal S|}\rightarrow\mathbb R_+$. For fixed $\rho$ set $\varphi_\rho(V(t))=\frac{1}{2}\|V(t)-\lambda(\rho)\|_2^2$. Then,
$$
\begin{aligned}
\dot{\varphi}_\rho(V(t))&=\dot{V}(t)^{\top}(V(t)-\lambda(\rho)) \\
& =\nabla_V L(V(t), \rho)^\top (\lambda(\rho)-V(t)) \\
& \leqslant L(\lambda(\rho), \rho)-L(V(t), \rho) \leq 0 \quad \text { equality only for } V(t)=\lambda(\rho),\\
\end{aligned}
$$
where the first inequality is due to the convexity and the last inequality is due to the optimality of $\lambda(\rho)$. Due to the strong convexity, the last inequality is strict for any $V(t)\neq \lambda(\rho)$.

We next address the property~\hyperref[prop-A4]{(A4)}. The projected dynamic system~\eqref{eq: rhoode} has to be well-posed to make claims about its stability. The Lipschitz continuity of the gradient $\nabla_\rho L(\lambda(\cdot),\cdot)$ follows from the linearity of the $\lambda(\rho)$ and the Lipschitz continuity of $\nabla_\rho L(V,\cdot)$ on the compact set $H$. To clarify the properties of the projection term in (\ref{eq: rhoode}), we present the derivation of the projected dynamic following \citet{dupuis1987large}. Define the projection of vector $v\in\mathbb{R}^{|\mathcal{S}|\cdot|\mathcal{A}|}$ at $\rho\in H$ as 
\begin{align*}
    \pi(\rho,v):=\lim_{\delta\searrow0}\frac{\Pi_H[\rho+\delta v]-\rho}{\delta}.
\end{align*}
The projected dynamics of $\rho$~\eqref{eq: rhoode} are then given by 
\begin{align}
    \label{eq:rho:ode:pi}
    \dot{\rho}=\pi(\rho,v), \quad \text{with}\quad v=\nabla_\rho L(\lambda(\rho),\rho).
\end{align} 
Let $N(\rho)$ be the set-valued function that gives the outward normals at $\rho\in\partial H$, that is, $N(\rho)=\{\zeta \in \mathbb{R}^{|\mathcal{S}|\cdot|\mathcal{A}|} \ | \ \forall\rho^\prime\in H:\langle\zeta,\rho-\rho^\prime\rangle\geq0, \|\zeta\|=1\}$. Then, by \cite[Lemma 4.6]{dupuis1987large}, the system~\eqref{eq:rho:ode:pi} equals 
\begin{align}
    \label{eq:pi:decomp}
        \pi(\rho,v)=v+(\langle v,-\zeta(\rho,v)\rangle\lor0)\zeta(\rho,v),
\end{align} 
with $\zeta(\rho,v)\in \arg\max_{\zeta\in -N(\rho)}\langle v,-\zeta\rangle$.
    
The existence and uniqueness of a solution to the projected dynamical system are guaranteed by \citet[Theorem 2 (uniqueness), Theorem 3 (existence)]{dupuis1993dynamical}. Note that due to the discontinuity of the dynamics, a solution is defined as an absolutely continuous function that fulfills the ODE save on a set of Lebesgue measure zero. To prove the global asymptotic stability of the projected ODE~\eqref{eq: rhoode}, we next provide a Lyapunov argument based on~\eqref{eq:pi:decomp}.

Due to the first order optimality condition, $\rho^\star$ is an equilibrium of the well-defined solution to the projected dynamic system (\ref{eq: rhoode}). The function 
\begin{equation}\label{eq: rho-lyapunov}
    \phi(\rho(t)):=L(V^\star,\rho^\star)-L(\lambda(\rho(t)),\rho(t))
\end{equation}
serves as a strict Lyapunov function. The inequality $\phi(\rho(t))\geq 0$ holds since $L(\lambda(\rho),\rho)=\min_{V\in\mathbb R^{|\mathcal S|}}L(V,\rho)\leq\max_{\rho\in H}\min_{V\in\mathbb R^{|\mathcal S|}}L(V,\rho)=L(V^\star,\rho^\star)$.  
The equality only holds for $\rho=\rho^\star$ due to the unique saddle-point as established in Corollary \ref{cor: dual-reg}. Define $\Gamma(\rho(t)):=\nabla_\rho L(\lambda(\rho(t)),\rho(t))$ and with a slight abuse of notation $\zeta(\rho(t)):=\zeta(\rho(t),\Gamma(\rho(t)))$. Then, 
\begin{align*}
    \dot\phi(\rho(t))&=-\Gamma(\rho(t))^\top \dot\rho(t)\\
                &=-\Gamma(\rho(t))^\top (\Gamma(\rho(t))+(\langle \Gamma(\rho(t)),-\zeta(\rho(t))\rangle\lor0)\zeta(\rho(t)))\\
                &=-\|\Gamma(\rho(t))\|^2+(\langle \Gamma(\rho(t)),-\zeta(\rho(t))\rangle\lor0)\langle \Gamma(\rho(t)),\zeta(\rho(t))\rangle\\
                &=-\|\Gamma(\rho(t))\|^2-(\langle \Gamma(\rho(t)),-\zeta(\rho(t))\rangle\lor0)\langle \Gamma(\rho(t)),-\zeta(\rho(t))\rangle\leq0, 
\end{align*}
where the last inequality is strict except for the equilibrium point $\rho^\star$ that satisfies the first-order optimality condition. In the second equality, we apply the decomposition (\ref{eq:pi:decomp}). Note that $\rho^\star\in H$. Hence, the first-order optimality condition also holds for the projected system.

Lastly, we address the stability of the iterate sequence, property~\hyperref[prop-A5]{(A5)}.     
The dual variable is constrained to the compact set $H$. We, therefore, only have to ensure that $\sup_k\|V_k\|<\infty$ almost surely. 
We apply the scaling limit method outlined in \citet{borkar2000ode} and \citet[Theorem 4.1]{borkar2023stochastic}, which states that the global asymptotic stability of the origin for a certain scaling limit ODE is a sufficient condition for the almost sure boundedness of the iterates. 
The scaling limit for the value iterates for a fixed occupancy $\rho$ is defined as 
$$h_\infty(V(t)):=\lim_{c\rightarrow\infty}-\frac{1}{c}\nabla_VL(cV(t),\rho).$$ 
As can be seen from the gradient~\eqref{eq: Vgrad}, the scaling limit exists for all $V\in\mathbb{R}^{|\mathcal{S}|}$ and is given by 
$$h_\infty(V(t))=-\eta_V V(t),$$
with global asymptotic equilibrium at the origin for $\eta_V>0$. 
This shows the almost sure boundedness of the value iterates. 
Note that the stability of the value iterates is achieved through the added convex regularization and does not hold with $\eta_V=0$ due to the lack of curvature in the unregularized linear program and its Lagrangian formulation~\eqref{eq: clLP}. 
\hfill $\blacksquare$

\section{Asynchronous Almost Sure Convergence: Proof Details}
\label{sec:async:PropProof}
This section contains the proofs of the asynchronous algorithm's properties listed in Proposition~\ref{prop:async:SA}. The proofs are presented in a series of lemmas.
For~\hyperref[prop-B4]{(B4)} the proof of the synchronous convergence analysis still applies. 
The two options regarding the behavioral policy update only need to be addressed separately for property~\hyperref[prop-B5]{(B5)}.

\begin{lemma}\label{lemma:async:noisebias}
    The noise components and bias terms of the gradient estimators~\eqref{eq:async:stochgrad} satisfy property~\hyperref[prop-B1]{(B1)}. The conditional expectation of the buffer-based transition probability estimator is
    \begin{align}
        \label{eq: CondExp}
        \mathbb{E}[\widehat{\mathcal{P}}_k(s'|s, a)|\mathcal F_k^-]=\frac{\nu_{k-1}(s,a)}{\nu_{k}(s,a)} \mathcal{P}_{\mathcal{D}_{k-1}}(s' | s, a) + \frac{\nu_{k}(s,a)-\nu_{k-1}(s,a)}{\nu_{k}(s,a)} \mathcal{P}(s'|s,a).
    \end{align}
\end{lemma}   
\textbf{Proof}\hspace{2pt}
    We first derive the conditional expectation of the buffer-based transition estimator.
    Fix $(s,a)\in\mathcal D_{\mathrm{inc}}(s_k)$, such that at iteration $k$ there is a draw from the corresponding buffer (otherwise
    $\widehat{\mathcal P}_k(\cdot|s,a)\equiv 0$ by convention and the claims are trivial).
    Recall that $\xi_k(s,a)$ is drawn uniformly from the list $\mathcal D_k(s,a)$ and
    $\widehat{\mathcal P}_k(s'|s,a)=\mathbf 1\{\xi_k(s,a)=s'\}$.
    If $(s,a)=X_k=(s_{k-1},a_{k-1})$, then $\mathcal D_k(s,a)=\mathcal D_{k-1}(s,a)\cup\{s_k\}$ and
    $\nu_k(s,a)=\nu_{k-1}(s,a)+1$.
    Conditionally on $\mathcal F_k^-=\sigma(\mathcal F_{k-1},s_k,a_k)$, the multiset
    $\mathcal D_k(s,a)$ is fixed and the uniform draw satisfies
    \[
    \mathbb P\big(\xi_k(s,a)\in \mathcal D_{k-1}(s,a)| \mathcal F_k^-\big)=\frac{\nu_{k-1}(s,a)}{\nu_k(s,a)},\qquad
    \mathbb P\big(\xi_k(s,a)=s_k| \mathcal F_k^-\big)=\frac{1}{\nu_k(s,a)}.
    \]
    Moreover, conditional on $\mathcal F_k^-$ and on the event $\{\xi_k(s,a)\in\mathcal D_{k-1}(s,a)\}$,
    the law of $\xi_k(s,a)$ is the empirical distribution $\mathcal P_{\mathcal D_{k-1}}(\cdot|s,a)$.
    Therefore,
    \[
    \mathbb E[\widehat{\mathcal P}_k(s'|s,a)| \mathcal F_k^-]
    =\frac{\nu_{k-1}(s,a)}{\nu_k(s,a)}\,\mathcal P_{\mathcal D_{k-1}}(s'|s,a)
    +\frac{1}{\nu_k(s,a)}\,\mathcal P(s'|s,a).
    \]
    If $(s,a)\neq X_k$, then the buffer cell is not updated, that is, $\mathcal D_k(s,a)=\mathcal D_{k-1}(s,a)$ and
    $\nu_k(s,a)=\nu_{k-1}(s,a)$, hence
    $\mathbb E[\widehat{\mathcal P}_k(s'|s,a)| \mathcal F_k^-]=\mathcal P_{\mathcal D_{k-1}}(s'|s,a)$.
    Combining both cases yields~\eqref{eq: CondExp}.
    
    Based on this explicit conditional expectation, we turn to the noise and bias analysis. 
    The noise and bias sequences of the gradient estimators are defined in equations~\eqref{eq:SA:def:MG},~\eqref{eq:SA:def:bias}.    
    By construction $\mathbb E[M_k^{(1)}|\mathcal F_k^-]=0$. Since $\mathcal F_{k-1}\subset\mathcal F_k^-$, using the tower-property gives $\mathbb E[M_k^{(1)}|\mathcal F_{k-1}]=\mathbb E\left[\mathbb E[M_k^{(1)}|\mathcal F_k^-]|\mathcal F_{k-1}\right]=0$. For $M_k^{(2)}$, we use $s_k\sim\mathcal P(\cdot|X_k)$ and $X_k, V_{k-1}\in \mathcal F_{k-1}$ to get $\mathbb E[M_k^{(2)}|\mathcal F_{k-1}]=0$. This shows that the noise sequences $\{M_k^{(1)}\}, \{M_k^{(2)}\}$ are martingale difference sequences with respect to $\{\mathcal F_k\}$ and~\hyperref[prop-B1]{(B1)}-a holds.

    The bounds on the second moments of the martingale differences,~\hyperref[prop-B1]{(B1)}-b, can be derived similarly to the proof of Proposition~\ref{prop: sync:gradient:noise} upon noting that $\|\widehat{\mathcal P}_k(\cdot|s,a)\|_2\le 1$ and\\ $\|\mathbb E[\widehat{\mathcal P}_k(\cdot|s,a)|\mathcal F_k^-]\|_2\le 1$. 
    Recall
    \[
    \mathcal{E}_k
    =\gamma\sum_{(s,a)\in\mathcal{S}\times\mathcal{A}}\rho_{k-1}(s,a)
    \Big(\mathbb E\big[\widehat{\mathcal{P}}_k(\cdot| s,a)| \mathcal F_k^-\big]-\mathcal{P}(\cdot| s,a)\Big).
    \]
    Fix $(s,a)\in\mathcal X$ and $s'\in\mathcal S$. In case that $(s,a)\in\mathcal D_{\mathrm{inc}}(s_k)$, the estimator $\widehat{\mathcal P}_k(s'| s,a)=\mathbf 1\{\xi_k(s,a)=s'\}$ is a one-hot vector of a sample
    $\xi_k(s,a)$ drawn uniformly from the list $\mathcal D_k(s,a)$. Hence,
    \[
    \mathbb E\big[\widehat{\mathcal P}_k(s'| s,a)| \mathcal F_k^-\big]
    =\mathcal P_{\mathcal D_k}(s'| s,a).
    \]
    By convention, if $D_k(s,a)$ is not drawn from at time $k$, then $\widehat{\mathcal P}_k(\cdot| s,a)\equiv 0$ and the same
    identity holds with $\mathcal P_{\mathcal D_k}(\cdot| s,a)\equiv 0$.
    
    Now note that $\mathcal D_k(s,a)$ consists of the subsequence of next-states observed upon visits to $(s,a)$, that is,
    \[
    \mathcal D_k(s,a)=\{\, s_{\ell}\,:\ 1\le \ell\le k,\ X_\ell=(s,a)\,\},
    \qquad  |\mathcal D_k(s,a)|=\nu_k(s,a).
    \]
    Conditional on $\{X_\ell=(s,a)\}$, the next state satisfies $s_{\ell}\sim \mathcal P(\cdot| s,a)$.
    Under~\hyperref[prop-B5]{(B5)} we have $\nu_k(s,a)\to\infty$ a.s. for all $(s,a)\in\mathcal X$, and thus by the strong law of large numbers,
    \[
    \mathcal P_{\mathcal D_k}(\cdot| s,a)\ \xrightarrow[k\to\infty]{\text{a.s.}}\ \mathcal P(\cdot| s,a)
    \quad \text{for each fixed }(s,a)\in\mathcal X.
    \]
    Since $H$ is compact and $\rho_{k-1}\in H$ for all $k$, the weights $\rho_{k-1}(s,a)$ are uniformly bounded.
    Therefore, dominated convergence yields
    \begin{equation*}
    \big\|\mathcal E_k\big\|
    \le \gamma \sum_{(s,a)\in\mathcal X} \big|\rho_{k-1}(s,a)\big|\,
    \big\| \mathcal P_{\mathcal D_k}(\cdot| s,a)-\mathcal P(\cdot| s,a)\big\|
    \ \xrightarrow[k\to\infty]{\text{a.s.}}\ 0.
    \end{equation*}
    Moreover, boundedness of $\rho_{k-1}$ and $\|\mathcal P_{\mathcal D_k}(\cdot| s,a)\|_2\le 1$ imply that $\{\mathcal E_k\}$ is bounded uniformly. This proves~\hyperref[prop-B1]{(B1)}-c.

\hfill$\blacksquare$

\begin{lemma}\label{lemma:async:dynamics}
    The non-autonomous ODEs~\eqref{eq:async:ODE} satisfy properties~\hyperref[prop-B2]{(B2)} and~\hyperref[prop-B3]{(B3)}.    
\end{lemma}
\textbf{Proof}\hspace{2pt}
    The properties~\hyperref[prop-B2]{(B2)} and~\hyperref[prop-B3]{(B3)} regarding the equilibria of the dynamical systems, hold for the asynchronous projected stochastic ascent-descent scheme because we can relate the asymptotic behavior of any solution to the non-autonomous rescaled ODEs~\eqref{eq:async:ODE} to the limiting ODEs of the synchronous setting~\eqref{eq:sync:ODEs} following~\cite[Section 6.4]{borkar2023stochastic}. 
    
    To see that~\hyperref[prop-B2]{(B2)} and~\hyperref[prop-B3]{(B3)} hold, we specify strict Lyapunov functions for the ODEs~\eqref{eq:sync:ODEs} that keep their Lyapunov function properties for the non-autonomous rescaled ODEs. 
    We show the proof only for~\hyperref[prop-B2]{(B2)} as~\hyperref[prop-B3]{(B3)} follows with the same arguments.
    
    Let $\rho\in H$ be a fixed arbitrary dual variable. Let $\epsilon>0$ and $\Lambda^{(1)}(\cdot)$ be a measurable function with $\Lambda^{(1)}(t)\in\Omega^\epsilon_{|\mathcal{S}|}$ for all $t$. By definition, $\Lambda^{(1)}(t)$ is a diagonal matrix for each $t$ with  bounded diagonal entries $w_i(t)\in[\epsilon,1]$. The measurability and boundedness of $t\mapsto\Lambda^{(1)}(t)$ together with the Lipschitz continuity of $V\mapsto\nabla_VL(V,\rho)$ imply that the vector field $-\Lambda^{(1)}(t)\nabla_VL(V,\rho)$ is measurable in $t$ and locally Lipschitz in $V$. Hence, for every initial condition, the ODE~\eqref{eq:async:V:ODE} admits a (unique) Carath\'eodory solution; see, e.g., standard results in differential equations~\citep{filippov2013differential}.
    
    Define the Lyapunov function
    \[
    \Phi_\rho(V):= L(V,\rho)-L(\lambda(\rho),\rho),
    \]
    where $\lambda(\rho)$ is the unique minimizer of $V\mapsto L(V,\rho)$ (unique by strong convexity).
    Along any solution $V(\cdot)$ of~\eqref{eq:async:V:ODE}, using $\nabla_V\Phi_\rho(V)=\nabla_V L(V,\rho)$ we have
    \begin{align*}
    \frac{d}{dt}\Phi_\rho(V(t))
    &= \nabla_V L(V(t),\rho)^\top \dot V(t)
    = -\nabla_V L(V(t),\rho)^\top \Lambda^{(1)}(t)\nabla_V L(V(t),\rho) \\
    &= -\sum_{s\in\mathcal S} w_s(t)\big(\nabla_V L(V(t),\rho)(s)\big)^2
    \le -\epsilon \|\nabla_V L(V(t),\rho)\|_2^2 \;\le\; 0.
    \end{align*}
    Moreover, equality holds if and only if $\nabla_V L(V(t),\rho)=0$, that is, $V(t)=\lambda(\rho)$.
    By $\eta_V$-strong convexity of $V\mapsto L(V,\rho)$, the Polyak-\L{}ojasiewicz inequality holds:
    $\|\nabla_V L(V,\rho)\|_2^2 \ge 2\eta_V\,\Phi_\rho(V)$ for all $V$.
    Therefore,
    \[
    \frac{d}{dt}\Phi_\rho(V(t)) \;\le\; -2\epsilon\eta_V\,\Phi_\rho(V(t)),
    \]
    which implies $\Phi_\rho(V(t))\le e^{-2\epsilon\eta_V t}\Phi_\rho(V(0))$ and hence
    $V(t)\to \lambda(\rho)$ as $t\to\infty$. Thus $\lambda(\rho)$ is a globally asymptotically stable equilibrium
    for~\eqref{eq:async:V:ODE}, uniformly over measurable $\Lambda^{(1)}(\cdot)\in\Omega^\epsilon_{|\mathcal S|}$.
    
    For~\hyperref[prop-B3]{(B3)} the same arguments hold with the unchanged Lyapunov function used in the proof of condition~\hyperref[prop-A4]{(A4)}, see~\eqref{eq: rho-lyapunov}, since multiplying the vector field by a diagonal $\Lambda^{(2)}(t)\in\Omega^\epsilon_{|\mathcal{X}|}$ preserves negativity of the Lyapunov derivative.
\hfill$\blacksquare$
    
\begin{lemma}\label{lemma:B7}
    Under Assumption~\ref{ass:ergodicity}, the update component selections in Algorithm~\ref{alg2} fulfill property~\hyperref[prop-B5]{(B5)}.
\end{lemma}
\textbf{Proof}\hspace{2pt}
    Let Assumption~\ref{ass:ergodicity} hold.
    In Algorithm~\ref{alg2}, conditional on the history $\mathcal F_k$, the next state-action pair $(s_{k+1},a_{k+1})$ is generated by first sampling $s_{k+1}\sim \mathcal P(\cdot| s_k,a_k)$ and then sampling $a_{k+1}\sim \pi_k^b(\cdot| s_{k+1})$, where the behavioral policy $\pi_k^b$ is either fixed or (in the on-policy exploration variant) a deterministic function of $\rho_k$ via~\eqref{eq:dual:policy}, optionally $\epsilon_k$-mixed with the uniform policy.
    Therefore, for any $(s,a),(s',a')\in\mathcal X$,
    \[
    \mathbb P\left((s_{k+1},a_{k+1})=(s',a')| \mathcal F_k\right)
    =\mathcal P(s'| s_k,a_k)\,\pi_k^b(a'| s').
    \]
    In particular, in the on-policy exploration case where $\pi_k^b=\pi^{\rho_k}$ (or its $\epsilon_k$-mixture), this can be written as
    \[
    \mathbb P\left((s_{k+1},a_{k+1})=(s',a')| \mathcal F_k\right)
    =\mathcal Q_{(s_k,a_k),(s',a')}(\rho_k),
    \]
    which shows~\hyperref[prop-B5]{(B5)}-a.
    By Assumption~\ref{ass:ergodicity}, for any strictly exploratory policy, the induced Markov chain on $\mathcal X$ is irreducible and aperiodic. Since $\rho\in H$ implies a uniform lower bound on action probabilities (and $\epsilon_k$-mixing preserves strict exploration), the policies $\pi^b_\rho$ are strictly exploratory uniformly over $\rho\in H$. Hence, $\mathcal Q(\rho)$ defines an irreducible and aperiodic Markov chain on $\mathcal X$ for all $\rho\in H$ and~\hyperref[prop-B5]{(B5)}-b holds. 

    Recall the definition of $H$ with constants $0<C^L\le C^U$ such that $C^L< \rho(s,a)< C^U$ for all $(s,a)$ and all $\rho\in H$.
    Then $\sum_{b\in\mathcal A}\rho(s',b)\ge |\mathcal A|C^L$ for every $s'\in\mathcal S$, and for all
    $\rho_1,\rho_2\in H$ and $a\in\mathcal A$,
    \begin{align*}
        |\pi_{\rho_1}(a|s')-\pi_{\rho_2}&(a|s')|=\left|\frac{\rho_1(s',a)}{\sum_{a'\in\mathcal{A}}\rho_1(s',a')}-\frac{\rho_2(s',a)}{\sum_{a'\in\mathcal{A}}\rho_2(s',a')}\right|\\
        &\leq\left|\frac{\rho_1(s',a)-\rho_2(s',a)}{\sum_{a'\in\mathcal{A}}\rho_1(s',a')}\right|+\left|\rho_2(s',a)\frac{\sum_{a'\in\mathcal{A}}(\rho_2(s',a')-\rho_1(s',a'))}{\sum_{a'\in\mathcal{A}}\rho_1(s',a')\sum_{a'\in\mathcal{A}}\rho_2(s',a')}\right|\\
        &\leq\frac{1}{|\mathcal{A}|C^L}|\rho_1(s',a)-\rho_2(s',a)|+\frac{C^U}{|\mathcal{A}|^2(C^L)^2}\sum_{a'\in\mathcal{A}}|\rho_1(s',a')-\rho_2(s',a')|\\
        &\leq C_\mathcal{Q}\sum_{a'\in\mathcal{A}}|\rho_1(s',a')-\rho_2(s',a')|,
    \end{align*}
    with $C_\mathcal{Q}:=\tfrac{1}{|\mathcal{A}|C^L}+\tfrac{C^U}{|\mathcal{A}|^2(C^L)^2}$, where the last inequality uses $$|\rho_1(s',a)-\rho_2(s',a)|\le\sum_{a'}|\rho_1(s',a')-\rho_2(s',a')|.$$
    Therefore, for any $(s,a),(s',a')\in\mathcal X$,
    \begin{align*}
    \big|\mathcal Q_{(s,a),(s',a')}(\rho_1)-\mathcal Q_{(s,a),(s',a')}(\rho_2)\big|
    &=\mathcal P(s'| s,a)\,\big|\pi^{\rho_1}(a'| s')-\pi^{\rho_2}(a'| s')\big|\\
    &\le C_{\mathcal Q}\sum_{b\in\mathcal A}\big|\rho_1(s',b)-\rho_2(s',b)\big|.
    \end{align*}
    Since the bound holds uniformly for all $(s,a),(s',a')$, the kernel mapping $\rho\mapsto \mathcal Q(\rho)$ is Lipschitz (for example in the max-entry norm), and hence in any equivalent norm on the finite-dimensional space of kernels. This shows~\hyperref[prop-B5]{(B5)}-c.
\hfill $\blacksquare$

\section{Asynchronous Convergence Rate: Proof Details}

This section contains the proof details of the convergence rate theorem, Theorem~\ref{thm:B1}. The first section addresses Remark~\ref{rem:UGE:dual:policies} and Proposition~\ref{prop: collected-facts-A2}, followed by the proofs of Propositions~\ref{prop: primal-tracking} and~\ref{prop: dual-variable-error} and finally some additional technical lemmata collected in the third subsection.
\subsection{Proofs Related to Assumption~\ref{ass:UGE:dual:policies} and Proposition~\ref{prop: collected-facts-A2}} \label{ssec:Proof:convergence:rate}
We first address the sufficient condition for the uniform geometric ergodicity with common constants of the family of Markov chains induced by the dual iterate.

\begin{lemma}[Deterministic mixing policy implies Assumption~\ref{ass:UGE:dual:policies}]
\label{lemma:UGE:dual:policies}
Assume that\\ there exists a deterministic state-feedback policy
\(\pi_d:\mathcal S\to\mathcal A\) such that the resulting Markov chain on \(\mathcal S\) with
transition kernel
\[
\mathcal P_d(s'|s):=\mathcal P\bigl(s'|s,\pi_d(s)\bigr)
\]
is uniformly geometrically ergodic. Then Assumption~\ref{ass:UGE:dual:policies}(i)
and (ii) hold.
Specifically, let \(M\in\mathbb N\), \(\epsilon_0\in(0,1]\), and a probability measure
\(\upsilon\) on \(\mathcal S\) be such that
\[
\mathcal P_d^{\,M}(\cdot|s)\ge \epsilon_0\,\upsilon(\cdot),
\qquad \forall s\in\mathcal S.
\]
Further define
\[
\pi_{\min}:=\frac{C^L}{|\mathcal A|\,C^U},
\qquad
\epsilon:=\pi_{\min}^M\epsilon_0,
\qquad
\varrho:=(1-\epsilon)^{1/M}.
\]
Then the following hold.

\begin{enumerate}
    \item[(a)] For every \(\rho\in H\), the state-action kernel \(\mathcal P_\rho\) has a
    unique stationary distribution \(\mu_\rho\), and
    \[
    \sup_{x\in\mathcal X}
    \|\delta_x \mathcal P_\rho^k-\mu_\rho\|_{\mathrm{TV}}
    \le
    C_{\mathcal X}\varrho^k,
    \qquad k\ge 0,
    \]
    with the common constant \(C_{\mathcal X}:=\varrho^{-M}\).

    \item[(b)] For every sequence \(\{\rho_t\}_{t\in\mathbb N}\subset H\),
    \[
    \sup_{k\ge 0}
    \delta\!\left(
    \prod_{t=k}^{k+m_\star-1}\mathcal P_{\rho_t}
    \right)
    \le
    \kappa,
    \qquad
    m_\star:=M+1,
    \quad
    \kappa:=1-\epsilon.
    \]
\end{enumerate}
\end{lemma}
\textbf{Proof}\hspace{2pt}
By the equivalence of uniform geometric ergodicity and the multi-step Doeblin condition,
see \citet[Theorem 16.0.2]{ref:meyn2012:markov}, there exist
\(M\in\mathbb N\), \(\epsilon_0\in(0,1]\), and a probability measure \(\upsilon\) on
\(\mathcal S\) such that
\[
\mathcal P_d^{\,M}(\cdot|s)\ge \epsilon_0\,\upsilon(\cdot),
\qquad \forall s\in\mathcal S.
\]
For \(\rho\in H\), let \(\pi_\rho\) denote the dual-induced policy~\eqref{eq:dual:policy} and define the
corresponding state-transition kernel
\[
\tilde{\mathcal P}_\rho(s'|s)
:=
\sum_{a\in\mathcal A}\pi_\rho(a|s)\mathcal P(s'|s,a).
\]
Since \(C^L\le \rho(s,a)\le C^U\) for all \((s,a)\in\mathcal S\times\mathcal A\),
we have the uniform lower bound
\[
\pi_\rho(a|s)=\frac{\rho(s,a)}{\sum_{b\in\mathcal A}\rho(s,b)}
\ge
\frac{C^L}{|\mathcal A|\,C^U}
=
\pi_{\min},
\qquad \forall \rho\in H,\ s\in\mathcal S,\ a\in\mathcal A.
\]
Hence, for every \(\rho\in H\) and \(s,s'\in\mathcal S\),
\begin{align*}
\tilde{\mathcal P}_\rho(s'|s)
&=
\sum_{a\in\mathcal A}\pi_\rho(a|s)\mathcal P(s'|s,a)
\\
&\ge
\pi_\rho(\pi_d(s)|s)\,\mathcal P\bigl(s'|s,\pi_d(s)\bigr)
\\
&\ge
\pi_{\min}\,\mathcal P_d(s'|s).
\end{align*}
By monotonicity under composition, iterating the above inequality \(M\) times yields
\[
\tilde{\mathcal P}_\rho^{\,M}(\cdot|s)
\ge
\pi_{\min}^M\,\mathcal P_d^{\,M}(\cdot|s)
\ge
\pi_{\min}^M\epsilon_0\,\upsilon(\cdot)
=
\epsilon\,\upsilon(\cdot),
\qquad \forall \rho\in H,\ s\in\mathcal S.
\]
Therefore, for every \(\rho\in H\), the kernel \(\tilde{\mathcal P}_\rho^{\,M}\) satisfies
the same Doeblin minorization, and therefore
\[
\delta(\tilde{\mathcal P}_\rho^{\,M})\le 1-\epsilon,
\qquad \forall \rho\in H.
\]
We first use this to obtain uniform geometric ergodicity of the frozen state kernels.
Let \(k=mM+r\) with \(m\in\mathbb N_0\) and \(0\le r<M\). For any probability measures
\(\eta,\eta'\) on \(\mathcal S\), by submultiplicativity of Dobrushin's coefficient and
the bound \(\delta(Q)\le 1\) for any Markov kernel \(Q\),
\[
\|\eta \tilde{\mathcal P}_\rho^k-\eta' \tilde{\mathcal P}_\rho^k\|_{\mathrm{TV}}
\le
(1-\epsilon)^m\|\eta-\eta'\|_{\mathrm{TV}}
\le
\varrho^{-(M-1)}\varrho^k\|\eta-\eta'\|_{\mathrm{TV}}.
\]
In particular, each \(\tilde{\mathcal P}_\rho\) has a unique stationary law
\(\tilde\mu_\rho\) on \(\mathcal S\), and
\[
\sup_{s\in\mathcal S}
\|\delta_s\tilde{\mathcal P}_\rho^k-\tilde\mu_\rho\|_{\mathrm{TV}}
\le
\varrho^{-(M-1)}\varrho^k,
\qquad k\ge 0.
\]
We now lift this to the state-action chain on \(\mathcal X=\mathcal S\times\mathcal A\).
Define the kernels
\[
\mathcal T(s'|s,a):=\mathcal P(s'|s,a),
\qquad
\mathcal B_\rho\bigl(s,a\mid s'\bigr):=\mathbf 1_{\{s=s'\}}\pi_\rho(a|s).
\]
Then
\[
\mathcal P_\rho=\mathcal T\mathcal B_\rho,
\qquad
\tilde{\mathcal P}_\rho=\mathcal B_\rho\mathcal T.
\]
Therefore, if \(\tilde\mu_\rho\) is stationary for \(\tilde{\mathcal P}_\rho\), then
\[
\mu_\rho:=\tilde\mu_\rho\mathcal B_\rho
\]
is stationary for \(\mathcal P_\rho\), that is,
\[
\mu_\rho(s,a)=\tilde\mu_\rho(s)\pi_\rho(a|s).
\]
Moreover, for any \(x\in\mathcal X\) and any \(k\ge 1\),
\[
\delta_x\mathcal P_\rho^k
=
(\delta_x\mathcal T)\tilde{\mathcal P}_\rho^{\,k-1}\mathcal B_\rho.
\]
Using non-expansiveness of total variation under Markov kernels, we obtain
\begin{align*}
\|\delta_x\mathcal P_\rho^k-\mu_\rho\|_{\mathrm{TV}}
&=
\|(\delta_x\mathcal T)\tilde{\mathcal P}_\rho^{\,k-1}\mathcal B_\rho
-
\tilde\mu_\rho\mathcal B_\rho\|_{\mathrm{TV}}
\\
&\le
\|(\delta_x\mathcal T)\tilde{\mathcal P}_\rho^{\,k-1}-\tilde\mu_\rho\|_{\mathrm{TV}}
\\
&\le
\varrho^{-(M-1)}\varrho^{k-1}
=
\varrho^{-M}\varrho^k.
\end{align*}
Since \(\|\delta_x-\mu_\rho\|_{\mathrm{TV}}\le 1\le \varrho^{-M}\) for \(k=0\), this proves
assertion (a) with \(C_{\mathcal X}:=\varrho^{-M}\).

It remains to prove assetion (b). Fix an arbitrary sequence \(\{\rho_t\}_{t\in\mathbb N}\subset H\)
and a block starting at time \(i\). By the pointwise bound
\(\tilde{\mathcal P}_{\rho_t}\ge \pi_{\min}\mathcal P_d\) derived above, monotonicity under
composition gives
\[
\prod_{t=i}^{i+M-1}\tilde{\mathcal P}_{\rho_t}(\cdot|s)
\ge
\pi_{\min}^M\mathcal P_d^{\,M}(\cdot|s)
\ge
\epsilon\,\upsilon(\cdot),
\qquad \forall s\in\mathcal S.
\]
Hence,
\[
\delta\!\left(
\prod_{t=i}^{i+M-1}\tilde{\mathcal P}_{\rho_t}
\right)
\le
1-\epsilon.
\]
Now using again the factorization of the state-action kernels, we have
\[
\prod_{t=i}^{i+M}\mathcal P_{\rho_t}
=
\mathcal T
\left(
\prod_{t=i}^{i+M-1}\tilde{\mathcal P}_{\rho_t}
\right)
\mathcal B_{\rho_{i+M}}.
\]
Therefore, for any probability measures \(\mu,\mu'\) on \(\mathcal X\),
\begin{align*}
\left\|
\mu \prod_{t=i}^{i+M}\mathcal P_{\rho_t}
-
\mu' \prod_{t=i}^{i+M}\mathcal P_{\rho_t}
\right\|_{\mathrm{TV}}
&\le
\left\|
(\mu\mathcal T)
\prod_{t=i}^{i+M-1}\tilde{\mathcal P}_{\rho_t}
-
(\mu'\mathcal T)
\prod_{t=i}^{i+M-1}\tilde{\mathcal P}_{\rho_t}
\right\|_{\mathrm{TV}}
\\
&\le
(1-\epsilon)\|\mu\mathcal T-\mu'\mathcal T\|_{\mathrm{TV}}
\\
&\le
(1-\epsilon)\|\mu-\mu'\|_{\mathrm{TV}}.
\end{align*}
Taking the supremum over \(\mu,\mu'\) yields
\[
\delta\!\left(
\prod_{t=i}^{i+M}\mathcal P_{\rho_t}
\right)
\le
1-\epsilon.
\]
Since \(i\) was arbitrary, Assumption~\ref{ass:UGE:dual:policies}(ii) holds with
\[
m_\star=M+1,
\qquad
\kappa=1-\epsilon.
\]
This completes the proof.
\hfill $\blacksquare$

Next, we provide the proofs of the convergence rate-related properties of the iterates generated by Algorithm~\ref{alg2} with primal variable projection to $\mathcal V_r$, summarized in Proposition~\ref{prop: collected-facts-A2}. The visitation floor and replay-buffer bias decay results are of particular interest since they are the main differences of our asynchronous approach compared to \cite{ref:zeng2024:two}. 

Throughout, assume that the exploration of the state-action space in Algorithm~\ref{alg2} is performed based on the dual-induced strictly exploratory policies. The alternatives of the fixed exploratory behavioral policy and $\epsilon$-exploration follow along the same arguments but with simpler constants and are left out.

We start with the visitation floor property, stated in Proposition~\ref{prop: collected-facts-A2}(i).
\begin{lemma}
\label{lemma: visitation-floor-envelopes}
The following hold
\begin{enumerate}
\item[(i)]
There exists $p_\star\in(0,1)$ such that for every fixed $\rho\in H$,
$$
\mu_\rho(s,a)\ge p_\star,\qquad\forall (s,a)\in\mathcal X.
$$

\item[(ii)]
For every $\delta\in(0,1)$, there exists a deterministic burn-in time
$K(\delta)<\infty$ such that 
$$\mathbb P\left(\underline{\nu}_{k}\ge\frac{p_\star}{2} k\quad \text{for all }k\ge K(\delta)\right)\ge1-\delta.$$

\item[(iii)]
The envelopes $\alpha(\underline{\nu}_k),\ \beta(\underline{\nu}_k)$ are non-increasing in $k$ and for all $k\ge K(\delta)$,
$$\bar\alpha_k\le\alpha(\underline{\nu}_k)\le\alpha(\lfloor \tfrac{p_\star}{2}k\rfloor),
\qquad
\ \bar\beta_k\le \beta(\underline{\nu}_k)\le\beta(\lfloor\tfrac{p_\star}{2}k\rfloor)
$$ 
with probability of at least $1-\delta$.
\end{enumerate}
\end{lemma}
\textbf{Proof}\hspace{2pt}
(i) Since $\rho\in H$ implies $C^L\le \rho(s,a)\le C^U$, we have for each $s$,
$\pi_\rho(a| s)=\rho(s,a)/\tilde\rho(s)\ge C^L/(|\mathcal A|C^U)$. Under
Assumption~\ref{ass:UGE:dual:policies}, every state transition kernel corresponding to a policy $\pi_\rho$, induces a corresponding unique stationary law $\tilde \mu_{\rho}$ with all components strictly positive. 
The compactness of $H$ and continuity of $\rho\mapsto\pi_\rho$ and $\pi_\rho\mapsto \tilde \mu_{\rho}$ yield the existence of a uniform lower bound $c_S:=\inf_{\rho\in H}\min_{s\in\mathcal S} \tilde \mu_{\rho}(s)>0$.
Thus $\mu_\rho(s,a)=\tilde\mu_{\rho}(s)\pi_\rho(a| s)\ge c_S\cdot C^L/(|\mathcal A|C^U)=:p_\star$.

(ii)
To show the high-probability linear growth bound for all state-action pairs, we first establish lower bounds for the time-varying expected visits, then apply a concentration bound for inhomogeneous Markov chains by \cite{ref:paulin2015:concentration} to establish concentration for any time step $k$, and lastly apply a union bound to achieve the high-probability minimal growth bound for all steps $k\ge K(\delta)$.

We begin by deriving a block contraction of the time-inhomogeneous kernel sequence with Dobrushin's ergodic coefficient.

By \cite[Corollary 4.3.17]{ref:bremaud2020:markov} the coefficient satisfies the contraction inequality $\|\mu Q-\nu Q\|_{\mathrm{TV}}\le \delta(Q)\|\mu-\nu\|_{\mathrm{TV}}$ for any two probability measures $\mu,\ \nu$ on $\mathcal X$.

Let $x_0\in\mathcal X$ denote the initial state-action pair. Write $p_k:=\delta_{x_0}\prod_{t=0}^{k-1}\mathcal P_{\rho_t}$ with the time-inhomogeneous kernel sequence $\{\mathcal P_{\rho_k}\}$, and let
$$d_k:=\|p_k-\mu_{\rho_k}\|_{\mathrm{TV}}.$$ Under Assumption~\ref{ass:UGE:dual:policies}, we can choose constants $m_\star\in\mathbb N$ and $\kappa\in(0,1)$ such that
\[
\sup_{k\ge 0}\delta\left(\prod_{t=k}^{k+m_\star-1}\mathcal P_{\rho_t}\right)\le\kappa.
\]
Then, by the contraction inequality, for any probability measure $q$ on $\mathcal X$,
\[
\left\|(q-\mu_\rho)\prod_{t=k}^{k+m_\star-1}\mathcal P_{\rho_t}\right\|_{\mathrm{TV}}\ \le\  \delta\left(\prod_{t=k}^{k+m_\star-1}\mathcal P_{\rho_t}\right)\,\|q-\mu_\rho\|_{\mathrm{TV}} \le \kappa \|q-\mu_\rho\|_{\mathrm{TV}}.
\]
Using this inhomogeneous Markov chain block contraction and the triangle inequality, for all $k\ge0$,
\begin{equation}
\label{eq:block-recursion}
\begin{aligned}
d_{k+m_\star}&=\left\|p_k \prod_{t=k}^{k+m_\star-1}\mathcal P_{\rho_t}-\mu_{\rho_{k+m_\star}}\right\|_{\mathrm{TV}}\\
&\le \left\|(p_k -\mu_{\rho_k})\prod_{t=k}^{k+m_\star-1}\mathcal P_{\rho_t}\right\|_{\mathrm{TV}}
+\ \left\|\mu_{\rho_k}\prod_{t=k}^{k+m_\star-1}\mathcal P_{\rho_t}-\mu_{\rho_{k+m_\star}}\right\|_{\mathrm{TV}}\\
&\le \kappa d_k + \sum_{t=k}^{k+m_\star-1}\|\mu_{\rho_t}-\mu_{\rho_{t+1}}\|_{\mathrm{TV}}.
\end{aligned}
\end{equation}
The last inequality holds due to the following telescoping argument. Insert and subtract the intermediate measures and use the triangle inequality:
\begin{align*}
    \left\|\mu_{\rho_k}\prod_{t=k}^{k+m_\star-1}\mathcal P_{\rho_t}-\mu_{\rho_{k+m_\star}}\right\|_{\mathrm{TV}}
    &\le \sum_{t=k}^{k+m_\star-1}\left\|\mu_{\rho_t}\prod_{j=t}^{k+m_\star-1}\mathcal P_{\rho_j}-\mu_{\rho_{t+1}}\prod_{j=t+1}^{k+m_\star-1}\mathcal P_{\rho_j}\right\|_{\mathrm{TV}}\\
    &=\sum_{t=k}^{k+m_\star-1}\left\|(\mu_{\rho_t}\mathcal P_{\rho_t}-\mu_{\rho_{t+1}})\prod_{j=t+1}^{k+m_\star-1}\mathcal P_{\rho_j}\right\|_{\mathrm{TV}},   
\end{align*}
where the empty product is taken as the identity.
Using the stationary law property of $\mu_{\rho_t}\mathcal P_{\rho_t}=\mu_{\rho_t}$ and the non-expansivity of Markov kernels, each summand is bounded by the corresponding term $\|\mu_{\rho_t}-\mu_{\rho_{t+1}}\|$ and~\eqref{eq:block-recursion} follows.

The stationary distributions are Lipschitz continuous in the dual variable (see Lemma \ref{lemma:Lip-stationary}), so there exists $L_\mu>0$ such that
\(
\|\mu_{\rho_t}-\mu_{\rho_{t+1}}\|_{\mathrm{TV}}
\ \le\ L_\mu\|\rho_{t+1}-\rho_t\|_2.
\)
Therefore,
\[
\sum_{t=k}^{k+m_\star-1}\|\mu_{\rho_t}-\mu_{\rho_{t+1}}\|_{\mathrm{TV}}\le L_\mu\ \sum_{t=k}^{k+m_\star-1}\|\rho_{t+1}-\rho_t\|_2.
\]
Since 
\begin{equation}\label{eq: exp-nu}
    \mathbb E[\nu_k(s,a)]\ =\ \sum_{t=0}^{k-1} p_t(s,a),
\end{equation}
and by definition $p_k(s,a) \geq \mu_{\rho_k}(s,a) - d_k$, we next aim to bound $\sum_{t=0}^{k-1}d_t$ from above. For this, we partition the timeline into $m_\star$-residual classes.
For each $\ell\in\{0,\dots,m_\star-1\}$, define
\[
g^{(\ell)}_n:=d_{\ell+n m_\star},
\qquad
S^{(\ell)}_n:=\sum_{i=\ell+n m_\star}^{\ell+(n+1)m_\star-1}\|\rho_{i+1}-\rho_i\|_2,
\qquad
n=0,\dots,N_\ell,
\]
where $N_\ell:=\lfloor (k-1-\ell)/m_\star\rfloor$.  
Note that the blocks form a partition of $\{0,\dots,k-2\}$, hence
$\sum_{\ell=0}^{m^\star-1}\sum_{n=0}^{N_\ell} S_n^{(\ell)}=\sum_{t=0}^{k-2}\|\rho_{t+1}-\rho_t\|_2$ and 
$\sum_{\ell=0}^{m_\star-1}\sum_{n=0}^{N_\ell}g^{(\ell)}_n
   =\sum_{t=0}^{k-1}d_t$.

Then~\eqref{eq:block-recursion} implies
\[
g^{(\ell)}_{n+1}\ \le\ \kappa\,g^{(\ell)}_n\ +\ L_\mu\,S^{(\ell)}_n.
\]
Summing over $n=0,\dots,N_\ell-1$ gives
\[
\sum_{n=0}^{N_\ell-1} g^{(\ell)}_{n+1}
\ \le\
\kappa \sum_{n=0}^{N_\ell-1} g^{(\ell)}_{n}
\ +\
L_\mu \sum_{n=0}^{N_\ell-1} S^{(\ell)}_{n}.
\]
Since
\[
\sum_{n=0}^{N_\ell-1} g^{(\ell)}_{n+1} - \kappa \sum_{n=0}^{N_\ell-1} g^{(\ell)}_{n}
= (1-\kappa)\sum_{n=0}^{N_\ell-1} g^{(\ell)}_{n+1}\ +\ \kappa\big(g^{(\ell)}_{N_\ell}-g^{(\ell)}_{0}\big),
\]
dropping the nonnegative term $\kappa g^{(\ell)}_{N_\ell}$ yields
\[
(1-\kappa)\sum_{n=0}^{N_\ell-1} g^{(\ell)}_{n+1}\ \le\ \kappa\, g^{(\ell)}_{0} \ +\ L_\mu \sum_{n=0}^{N_\ell-1} S^{(\ell)}_{n}.
\]
Dividing by $(1-\kappa)$ and adding $g^{(\ell)}_{0}$ to both sides,
\[
\sum_{n=0}^{N_\ell} g^{(\ell)}_{n} \ \le\ \Bigl(1+\frac{\kappa}{1-\kappa}\Bigr) g^{(\ell)}_{0}\ +\ \frac{L_\mu}{1-\kappa}\sum_{n=0}^{N_\ell-1} S^{(\ell)}_{n}
\ =\ \frac{1}{1-\kappa}\, g^{(\ell)}_{0}\ +\ \frac{L_\mu}{1-\kappa}\sum_{n=0}^{N_\ell-1} S^{(\ell)}_{n}.
\]
Summing this bound over $\ell=0,\dots,m_\star-1$ yields
\begin{equation*}
\sum_{t=0}^{k-1} d_t
\ \le\
\frac{1}{1-\kappa}\sum_{\ell=0}^{m_\star-1} d_\ell
\ +\
\frac{L_\mu}{1-\kappa}\sum_{t=0}^{k-2}\|\rho_{t+1}-\rho_t\|_2.
\end{equation*}
Since $d_\ell\le1$ for all $\ell$, $\sum_{\ell=0}^{m_\star-1} d_\ell\le m_\star$, which gives the final bound
\[
\sum_{t=0}^{k-1} d_t
\ \le\ \frac{m_\star}{1-\kappa}
\ +\ \frac{L_\mu}{1-\kappa}\sum_{t=0}^{k-2}\|\rho_{t+1}-\rho_t\|_2.
\]
Inserting this bound into~\eqref{eq: exp-nu}, for any fixed $(s,a)\in\mathcal X$,
\[
\mathbb E[\nu_k(s,a)]\ =\ \sum_{t=0}^{k-1} p_t(s,a)
\ \ge\ \sum_{t=0}^{k-1}\big(\mu_{\rho_t}(s,a)-d_t\big)
\ \ge\ k\,p_\star\ -\ \frac{m_\star}{1-\kappa}
\ -\
\frac{L_\mu}{1-\kappa}\sum_{i=0}^{k-2}\|\rho_{i+1}-\rho_i\|_2.
\]
By the 1-Lipschitz property of projection in $\ell_2$,
$
\|\rho_{t+1}-\rho_t\|_2\ \le\ \beta(\nu_t(s_t,a_t))\|\hat h_t\|,
$
where $\hat h_t$ is the single-coordinate gradient from~\eqref{eq: rstochgrad_alg2}. Due to the projection step, the iterates $V_t$ are bounded, and there is a constant $B<\infty$ with $\|\hat h_t\|\le B$ for all $t\ge 0$; hence
$$
\|\rho_{t+1}-\rho_t\|_2\ \le\ B\,\beta(\nu_t(s_t,a_t)).
$$
Note that $\sum_{(s,a)\in\mathcal X} \nu_k(s,a)=k$. The specific choice of the base stepsize schedule $\beta_k=\beta_0/(k+1)$ gives the following bound on the total step length up to iteration $k$
$$
\sum_{t=0}^{k-1}\beta(\nu_t(s_t,a_t))
=\sum_{(s,a)\in\mathcal X}\sum_{t=0}^{\nu_k(s,a)-1}\frac{\beta_0}{t+1}
\ \le\ \beta_0\,|\mathcal X|\,\sum_{t=0}^{k-1}\frac{1}{t+1}.
$$
Using $\int_t^{t+1}\frac{dx}{x}\ge \frac{1}{t+1}$ for $t\ge 1$, we get
$$
\sum_{t=0}^{k-1}\frac{1}{t+1}
=1+\sum_{t=1}^{k-1}\frac{1}{t+1}
\le 1+\int_{1}^{k}\frac{dx}{x}
=1+\log(k),
$$
and therefore, 
$$\sum_{t=0}^{k-1}\beta(\nu_t(s_t,a_t))\le\beta_0\,|\mathcal X|(1+\log k).$$
Hence, for any $\epsilon\in(0,p_\star)$ there exists $K_0(\epsilon)$ such that for all $k\ge K_0(\epsilon)$,
\begin{equation}\label{eq: exp-nu-bound}
\mathbb E\!\left[\frac{\nu_k(s,a)}{k}\right]\ \ge\ p_\star-\frac{\epsilon}{2}.
\end{equation}
Fix $(s,a)\in\mathcal X$. We next establish the concentration lower bound for any $k\in\mathbb N$.
Denote the path of the process $\{X_k\}$ up to iteration $k$ by $x_{0:k-1}\in\mathcal X^k$. Define the normalized count function for the fixed $(s,a)$ pair by 
\[
\varphi:\mathcal X^k\rightarrow\mathbb R,\ x_{0:k-1}\mapsto\varphi(x_{0:k-1}):=\frac{1}{k}\sum_{t=0}^{k-1}\mathbf 1\{x_t=(s,a)\}.
\]
Then, $\varphi(x_{0:k-1})=\tfrac{\nu_k(s,a)}{k}$. 
Let $y_{0:k-1}\in\mathcal X^k$ denote a second path. It can be seen that 
$$\varphi(x_{0:k-1})-\varphi(y_{0:k-1})\le\sum_{t=0}^{k-1}c_t\mathbf{1}\{x_t\ne y_t\},\quad\text{with } c_t=\frac{1}{k},\ \sum_{t=0}^{k-1}c_t^2=\frac{1}{k}$$ 
and we can apply the one-sided McDiarmid inequality for time-inhomogeneous Markov chains \cite[Corollary 2.10]{ref:paulin2015:concentration} to $\varphi$ and obtain
\begin{subequations}
    \begin{equation}\label{eq:Paulin-tail}
    \mathbb P\Big(\frac{\nu_k(s,a)}{k}-\mathbb E\big[\frac{\nu_k(s,a)}{k}\big]\le -\eta\Big)
    \ \le\
    \exp\Big(-\,\frac{2k\,\eta^2}{\tau_{\min}}\Big),
    \qquad
    \tau_{\min}:=\inf_{0\le \zeta<1}\tau_\mathrm{inh}(\zeta)\Big(\tfrac{2-\zeta}{1-\zeta}\Big)^2,
    \end{equation}
where 
    \begin{equation}
    \label{eq:inh-mixing}
    \tau_\mathrm{inh}(\zeta):=\min\{\ell\in\mathbb N:\sup_{n\in\mathbb N}\sup_{x,y\in\mathcal X}\|\delta_x\prod_{t=n}^\ell\mathcal P_{\rho_t}-\delta_y\prod_{t=n}^\ell\mathcal P_{\rho_t}\|_{\mathrm{TV}}\le\zeta\},
    \end{equation}
\end{subequations}
see \citep[Definition~1.4]{ref:paulin2015:concentration}, is the inhomogeneous mixing time for $\{\mathcal P_{\rho_k}\}$. It holds that $\tau_{\min}<\infty$ with a simple bound in our setting deferred to Lemma~\ref{lemma: inhom-mixing-time-bound}.
From the mean lower bound in~\eqref{eq: exp-nu-bound}, there exists $K_0\big(\tfrac{p_\star}{2}\big)$ such that, for all $k\ge K_0\big(\tfrac{p_\star}{2}\big)$,
\[
\mathbb E\left[\tfrac{\nu_k(s,a)}{k}\right]\ \ge\ \tfrac{3}{4}p_\star.
\]
Hence, on the event $\{\tfrac{\nu_k(s,a)}{k}\le \tfrac{p_\star}{2}\}$ we have a downward deviation of at least
$\mathbb E[\tfrac{\nu_k(s,a)}{k}]-\tfrac{p_\star}{2}\ \ge\ \tfrac{p_\star}{4}$, and applying~\eqref{eq:Paulin-tail} yields
\begin{equation}\label{eq:Paulin-fixed-k}
\mathbb P\left(\tfrac{\nu_k(s,a)}{k}\le \tfrac{p_\star}{2}\right)
\ \le\
\exp\Big(-\frac{k\ p_\star^2}{8\,\tau_{\min}}\Big),
\qquad k\ \ge\ K_0\Big(\tfrac{p_\star}{2}\Big).
\end{equation}

Union bound over $(s,a)\in\mathcal X$ and sum the geometric tail over $k\ge K$:
$$
\mathbb P\Big(\exists k\ge K:\min_{(s,a)\in\mathcal{X}}\tfrac{\nu_k(s,a)}{k}<p_\star-\epsilon\Big)
\le|\mathcal X|\sum_{k=K}^\infty e^{-ck}
=\frac{|\mathcal X|}{1-e^{-c}}\,e^{-cK},\quad
c=\frac{p_\star^2}{8\,\tau_{\min}}.
$$
Thus, for any $\delta\in(0,1)$, taking
\begin{equation}\label{eq:K-delta-Paulin}
K(\delta)
\ :=\
\max\left\{
K_0\Big(\tfrac{p_\star}{2}\Big),\
\left\lceil \frac{1}{c}\,
\log\frac{|\mathcal X|}{(1-e^{-c})\,\delta}\right\rceil
\right\},
\end{equation}
ensures
\(
\mathbb P\big(\forall k\ge K(\delta):\ \underline{\nu}_k/k\ge p_\star/2\big)\ \ge\ 1-\delta.
\)
which completes the proof of (ii).

(iii)
On the minimum-visitation event in (ii), for all $k\ge K(\delta)$ and all $(s,a)\in\mathcal X$ we have $\nu_k(s,a)\ge (p_\star/2)\,k$ and
$\tilde\nu_k(s)\ge (p_\star/2)\,k$. Since $\alpha(\cdot)$ and $\beta(\cdot)$ are nonincreasing,

$$
\alpha\big(\tilde\nu_k(s)\big)\ \le\ \alpha\left(\Big\lfloor\tfrac{p_\star}{2}\,k\Big\rfloor\right),
\qquad
\beta\big(\nu_k(s,a)\big)\ \le\ \beta\left(\Big\lfloor\tfrac{p_\star}{2}\,k\Big\rfloor\right). 
$$
\hfill $\blacksquare$

Next, we prove the replay-buffer bias decay stated in Proposition~\ref{prop: collected-facts-A2}(ii). In the following, we use the shorthand notation $\bar{\mathcal E}_k(\rho)$ for the stationary law expectation of the bias vector given a fixed $\rho\in H$, that is,
\begin{equation}
\label{eq: stationary-law-bias}
    \bar{\mathcal E}_k(\rho):= \mathbb{E}_{s\sim\tilde{\mu}_\rho}[e_s\,\mathcal E_k(\rho)(s)].
\end{equation}

\begin{lemma}[Replay-buffer bias under high-probability visitation]\hspace{2pt}
\label{lemma:buffer-bias-HP}
Fix \(\delta \in (0,1)\), let \(\mathcal G_\delta\) denote the high-probability event of Lemma~\ref{lemma: visitation-floor-envelopes}, and let \(K(\delta)\) be the corresponding burn-in time. Define
\[
\Lambda(\delta)\ :=\ \sqrt{\log\!\Big(\frac{(2^{|\mathcal S|}-2)\,|\mathcal X|}{\delta}\Big)}.
\]
Then on \(\mathcal G_\delta\) and for all \(k\ge K(\delta)\),
\begin{equation}\label{eq:buffer-bias-bounds}
\big\|\bar{\mathcal E}_k(\rho)\big\|_\infty\ \le\ C_{\mathrm{buf}}(\delta)\frac{1}{\sqrt{k}},\quad C_{\mathrm{buf}}(\delta):=\frac{3\gamma |\mathcal X|C^U\Lambda(\delta)}{\sqrt{p_\star}}.
\end{equation}
In particular, \(\|\bar{\mathcal E}_k(\rho)\|_\infty=O\big(\sqrt{1/k}\big)\) on \(\mathcal G_\delta\), uniformly in $\rho,\,V$.
\end{lemma}
\textbf{Proof}\hspace{2pt}
For each fixed $(s,a)\in \mathcal X$ and corresponding buffer size \(|\mathcal D_k(s,a)|=\nu_k(s,a)\), the \(L^1\)-bound for a categorical law on \(|\mathcal S|\) outcomes \citep[Theorem 2.1]{ref:weissman2003:inequalities} gives for all $\epsilon>0$
\[
\mathbb P\!\left(\big\|\mathcal{P}_{\mathcal D_k}(\cdot| s,a)-\mathcal P(\cdot| s,a)\big\|_1\ge \epsilon\right)
\ \le\ (2^{|\mathcal S|}-2)\,\exp\!\Big(-\tfrac{\nu_k(s,a)\epsilon^2}{2}\Big).
\]

Condition on the \(\sigma\)-field generated by the realized counts \(\{\nu_k(s,a)\}_{(s,a)\in\mathcal X}\).
To union bound over the \(|\mathcal X|\) pairs, we choose
\[
\epsilon_{k}\ :=\ \sqrt{\frac{2}{\underline{\nu}_k}\ \log\!\Big((2^{|\mathcal S|}-2)\,|\mathcal X|\,/\delta\Big)},
\]
Then
\[
\mathbb P\Big(\bigcup_{(s,a)\in\mathcal X}
\{\|\mathcal P_{\mathcal D_k}(\cdot| s,a)-\mathcal P(\cdot| s,a)\|_1\ge \epsilon_{k}\}\ \Bigm|\ \sigma(\{\nu_k(s,a)\}_{(s,a)\in\mathcal X})\Big)
\ \le\ \delta.
\]

Therefore, with probability at least \(1-\delta\), simultaneously for all \((s,a)\in\mathcal X\),
\[
\big\|\mathcal P_{\mathcal D_k}(\cdot| s,a)-\mathcal P(\cdot| s,a)\big\|_1
\ \le\ \sqrt{\frac{2}{\underline{\nu}_k}\ \log\!\Big((2^{|\mathcal S|}-2)\,|\mathcal X|/\delta\Big)}
\ \le\ \sqrt{\frac{2}{\underline{\nu}_k}}\,\Lambda(\delta).
\]
On the high-probability visitation event \(\mathcal G_\delta\) of Lemma~\ref{lemma: visitation-floor-envelopes}, we have
\(\underline{\nu}_{k}\ge \tfrac{p_\star}{2}\,k\) for all \(k\ge K(\delta)\), hence
\begin{equation}
\label{eq:L1-kernel}
\max_{(s,a)\in\mathcal X}\ \big\|\mathcal P_{\mathcal D_k}(\cdot| s,a)-\mathcal P(\cdot| s,a)\big\|_1
\ \le\ \frac{2\Lambda(\delta)}{\sqrt{p_\star\,k}}.
\end{equation}
Now we pass from kernel error to the buffer-bias term. The stationary law expectations can be written as the Hadamard product between the buffer bias term defined in~\eqref{eq:SA:def:bias} and the strictly positive stationary state and state-action distributions induced by the dual policy, that is,

\begin{equation*}
    \bar{\mathcal E}_k(\rho)(s')=\tilde\mu_\rho(s')\, \gamma\sum_{(s,a)\in\mathcal S\times\mathcal A}\rho(s,a)\,\tfrac{\nu_{k-1}(s,a)}{\nu_k(s,a)}\left(\mathcal P_{\mathcal D_{k-1}}(s'| s,a)-\mathcal P(s'| s,a)\right).
 \end{equation*}
With the trivial bound $\|\tilde \mu_\rho\|_\infty\le1,$ and the increasing number of visits $\nu_k(s,a)\ge\nu_{k-1}(s,a)$, we can omit the preconditioning with the stationary distribution in the following norm bounds and find for \(\bar{\mathcal E}_k(\rho)\),
\begin{align*}
\big\|\bar{\mathcal E}_k(\rho)\big\|_\infty
\le& \gamma\ \Big\|\sum_{(s,a)\in\mathcal X}\rho(s,a)\big(\mathcal P_{\mathcal D_{k-1}}(\cdot|s,a)-\mathcal P(\cdot|s,a)\big)\Big\|_\infty\\
\le& \gamma\|\rho\|_1\ \max_{(s,a)\in\mathcal X}\big\|\mathcal P_{\mathcal D_{k-1}}(\cdot|s,a)-\mathcal P(\cdot|s,a)\big\|_1\le\frac{3\gamma |\mathcal X|C^U\Lambda(\delta)}{\sqrt{p_\star}}\frac{1}{\sqrt{k}},
\end{align*}
since \(\|\rho\|_1\le |\mathcal X|C^U\) by the restriction to \(\rho\in H\) and $\tfrac{2}{\sqrt{k-1}}\le\tfrac{2\sqrt{2}}{\sqrt{k}}\le\tfrac{3}{\sqrt{k}}$ for all $k\ge 2$. 
\hfill $\blacksquare$

We next address the regularity properties stated in Proposition~\ref{prop: collected-facts-A2}(iii). The $\eta_V$-strong convexity of the mapping $V\mapsto L(V,\rho)$ has been established already in equations~\eqref{eq: gradients} and the proof of Proposition~\ref{prop: ODEs}. The explicit Lipschitz constant of the best response mapping $\lambda$ is derived in the following.

\begin{lemma}
    \label{lemma: LipschitzLambda}
    The best response map $\lambda:H\rightarrow\mathbb{R}^{|\mathcal{S}|}$ is affine and hence globally Lipschitz on $H$. In particular, 
    \begin{equation}\label{eq:lambda:optimizer:form}
        \lambda(\rho)_{s'}=\frac{1}{\eta_V}\big(\tilde \rho(s')-\gamma\sum_{(s,a)\in\mathcal{S}\times\mathcal{A}}\rho(s,a)\mathcal{P}(s'|s,a)\big),
    \end{equation}
    and for all $\rho,\,\hat\rho\in H$, $$\|\lambda(\rho)-\lambda(\hat\rho)\|_2\le\frac{\sqrt{|\mathcal{S}||\mathcal{A}|(1+\gamma^2)}}{\eta_V}\|\rho-\hat\rho\|_2.$$
\end{lemma}
\textbf{Proof}\hspace{2pt}
    The displayed linear form~\eqref{eq:lambda:optimizer:form} of the optimizer follows from the first-order optimality condition derived in the proof of Proposition~\ref{prop: ODEs}. For the Lipschitz constant, we define the matrix $M\in\mathbb{R}^{|\mathcal S||\mathcal A|\times|\mathcal S|}$ with rows $$\text{row}_{sa}(M)=e^\top_s-\gamma\mathcal P(\cdot|s,a)^\top, \; \text{so that}\; \lambda(\rho)=\frac{1}{\eta_V}M^\top\rho.$$ Then $\|\lambda(\rho)-\lambda(\hat\rho)\|_2=\frac{1}{\eta_V}\|M^\top(\rho-\hat\rho)\|_2\leq\frac{\|M^\top\|_2}{\eta_V}\|\rho-\hat \rho\|_2$. Each row of $M$ satisfies $\|\text{row}_{sa}(M)\|_2^2\leq\|e^\top_s\|_2^2+\gamma^2\|\mathcal P(\cdot|s,a)\|_2^2\leq(1+\gamma^2)$. Therefore, the claim follows by noting $\|M^\top\|_2^2=\|M\|_2^2\leq\|M\|_F^2=\sum_{(s,a)\in\mathcal S\times\mathcal A}\|\text{row}_{sa}(M)\|_2^2\leq|\mathcal S||\mathcal A|(1+\gamma^2)$.
\hfill$\blacksquare$

We prove the L-smoothness of the reduced objective by showing the Lipschitz property for $\nabla f(\rho)$. 
Upon inserting the best response formula~\eqref{eq:lambda:optimizer:form} in the Lagrangian, we obtain 
$$
f(\rho)=L(\lambda(\rho),\rho)=-\frac{1}{2\eta_V}\|M^\top\rho\|^2_2+\langle\rho,r\rangle+\eta_\rho g(\rho),
$$ 
since $\Delta[V](s,a)=(-MV)(s,a)+r(s,a)$ and hence, $L(V,\rho)=\frac{\eta_V}{2}\|V\|^2_2-\langle\rho,MV\rangle+\langle\rho,r\rangle+\eta_\rho g(\rho)$. Then $\nabla f(\rho)=-\frac{MM^\top}{\eta_V}\rho+r+\eta_\rho\nabla g(\rho)$ and we obtain the following smoothness for the reduced objective:  For all$\,\rho,\hat\rho\in H$
\begin{equation*}
    \|\nabla f(\rho)-\nabla f(\hat\rho)\|_2\ \le\Big(\eta_\rho \frac{2}{C^L}+\frac{\|M\|^2_2}{\eta_V}\Big)\|\rho-\hat\rho\|_2\le\Big(\eta_\rho \frac{2}{C^L}+\frac{|\mathcal S||\mathcal A|(1+\gamma^2)}{\eta_V}\Big)\|\rho-\hat\rho\|_2,\ 
\end{equation*}
where $\tfrac{2}{C^L}$ is the Lipschitz constant of $\nabla g(\rho)$ on H as derived in the following Lemma, and where we have used the upper bound for $\|M\|_2^2$ from the proof of Lemma~\ref{lemma: LipschitzLambda}.

\begin{lemma}[Lipschitz gradient of $g$]
\label{lemma:lipschitz-grad-g}
For the entropy regularization term \(g(\rho)\), defined in~\eqref{eq:Entropy:regularizer}, it holds that, for all $\rho,\hat\rho\in H$,
$$
\|\nabla g(\rho)-\nabla g(\hat\rho)\|_2 \;\le\; \frac{2}{C^L}\,\|\rho-\hat\rho\|_2 .
$$
\end{lemma}
\textbf{Proof}\hspace{2pt}
The coordinate-wise gradient is
$$
\nabla g(\rho)(s,a)
= -\log\Big(\frac{\rho(s,a)}{\tilde\rho(s)}\Big)
= \log \tilde\rho(s) - \log\rho(s,a).
$$
The Hessian is block-diagonal across states. For a fixed $s$ and $a,a'\in \mathcal{A}$,
$$
\big[\nabla^2 g(\rho)\big]_{(s,a),(s,a')}
= \frac{1}{\tilde\rho(s)} - \frac{\mathbf{1}\{a=a'\}}{\rho(s,a)}.
$$
Hence the $|\mathcal A|\times |\mathcal A|$ block at state $s$ is
\[
H_s(\rho)\;=\; -\,\mathrm{diag}\,\Big(\tfrac{1}{\rho(s,a)}\Big)_{a\in \mathcal A} \;+\; \frac{1}{\tilde\rho(s)}\,\mathbf 1\mathbf 1^\top .
\]
Using the triangle inequality for operator norms and the facts
$\|\mathrm{diag}(u)\|_2=\|u\|_\infty$, $\|\mathbf 1\mathbf 1^\top\|_2=|\mathcal A|$, and $\tilde \rho(s)\ge{|\mathcal A|}C^L$ for $\rho\in H$, we get
$$
\|H_s(\rho)\|_2
\le\max_{a}\frac{1}{\rho(s,a)}+\frac{1}{\tilde\rho(s)}|\mathcal A|
\le\frac{1}{C^L} + \frac{|\mathcal{A}|}{|\mathcal{A}|C^L}
= \frac{2}{C^L}.
$$
Since $\nabla^2 g(\rho)$ is block-diagonal with blocks $H_s(\rho)$, it follows that
$$
\sup_{\rho}\,\|\nabla^2 g(\rho)\|_2 \;=\; \sup_{s,\rho}\,\|H_s(\rho)\|_2 \;\le\; \frac{2}{C^L}.
$$

Finally, 
$$
\|\nabla g(\rho)-\nabla g(\hat\rho)\|_2
\le \Big(\sup_{z}\|\nabla^2 g(z)\|_2\Big)\,\|\rho-\hat\rho\|_2
\le \frac{2}{C^L}\,\|\rho-\hat\rho\|_2.
$$
\hfill$\blacksquare$

The following Lemma states the strong concavity of the reduced objective on $H$.

\begin{lemma}[Strong concavity of the reduced objective]\label{lemma: strong-concavity-f} 
The reduced objective\\ $f(\rho) = \min_VL(V,\rho)$ is $\mu$-strongly concave on $H$ with any modulus $0<\mu\le\mu_\mathrm{opt}$, where 
\begin{equation}
\label{eq: strong-conc-f}
    \mu_{\mathrm{opt}}=\frac{1}{4}\Big[(A+B+C)-\sqrt{(A+B+C)^2-4AC}\Big] \;>\;0,
\end{equation}

with 
$$
A := \frac{\eta_\rho}{C^U}, \qquad 
B := \frac{|\mathcal{S}||\mathcal{A}|(1+\gamma^2)}{\eta_V}, \qquad
C := \frac{(1-\gamma)^2\,|\mathcal{A}|\,C^{L^{\,2}}}{\eta_V\,C^{U^{\,2}}}.
$$
\end{lemma}
\textbf{Proof}\hspace{2pt}
As shown previously, it holds that
$$
\nabla f(\rho) = -\tfrac{1}{\eta_V} M M^\top \rho+r+\eta_\rho \nabla g(\rho), 
\qquad 
\nabla^2 f(\rho) = -\tfrac{1}{\eta_V} M M^\top +\eta_\rho \nabla^2 g(\rho).
$$
Since \(g\) is concave, \(\nabla^2 g(\rho) \preceq 0\), and together with the negative semidefinite term \(- M M^\top\) this implies that \(f\) is concave.
For any $h \in\mathbb R^{|\mathcal S||\mathcal A|}$, 
\begin{equation}
\label{eq:hessian:f}
    -h^\top \nabla^2f(\rho)h=\frac{1}{\eta_V}\|M^\top h\|^2_2-\eta_\rho \ h^\top\nabla^2g(\rho) h. 
\end{equation}
For any \(h \in \mathbb R^{|\mathcal S||\mathcal A|}\), decompose \(h=u+v\) per state, where for each \(s\in\mathcal S\),
\(\sum_{a\in\mathcal A}u(s,a)=0\) and \(v(s,a)=q(s)\,\frac{\rho(s,a)}{\tilde\rho(s)}\) for some \(q(s)\in\mathbb R\).
With this decomposition, \(u\) redistributes mass across actions while keeping the state
marginal fixed, and \(v\) changes only the state marginals. In particular, $v^\top\nabla^2g(\rho)v=0$.
For any $u$ with $\sum_{a\in\mathcal A} u(s,a)=0$, it holds that for all $s\in\mathcal{S}$ and corresponding block matrices $H_s(\rho)$ in the Hessian of $g$, specified in the proof of Lemma~\ref{lemma:lipschitz-grad-g},
\[u(s,\cdot)^\top H_s(\rho)u(s,\cdot) = \sum_{a\in \mathcal A} \tfrac{u(s,a)^2}{\rho(s,a)} \ge\tfrac{1}{C^U}\,\sum_{a\in \mathcal A} u(s,a)^2.\]
Since the Hessian of g is block diagonal, it follows that uniformly on $H$,
\begin{equation}
    \label{eq:hessian:g:decomp}
    u^\top\big(-\eta_\rho\nabla^2 g(\rho)\big)u=\sum_{s\in\mathcal S}\big(u(s,\cdot)^\top H_s(\rho)u(s,\cdot)\big)\ge  \tfrac{\eta_\rho}{C^U} \|u\|_2^2=A\|u\|_2^2.
\end{equation}
It holds that $\|M^\top u\|_2\le \|M\|_2\|u\|_2$, and together with the operator bound $\|M\|^2_2\le|\mathcal S||\mathcal A|(1+\gamma^2)$ derived in the proof of Lemma~\ref{lemma: LipschitzLambda}, we have 
\begin{equation}
\label{eq:Mu-bound}
    \|M^\top u\|_2^2\le|\mathcal S||\mathcal A|(1+\gamma^2)\|u\|^2_2.
\end{equation}
Let again $\tilde{\mathcal P}_\rho(s'|s)=\sum_{a\in\mathcal A}\pi_\rho(a|s)\mathcal P(s'|s,a)$ denote the state transition kernel induced by the dual policy. And with some abuse of notation let $\bar P_\rho\in \mathbb R^{|\mathcal S|\times |\mathcal S|}$ denote the corresponding state transition matrix with $[\bar P_\rho]_{ss'}=\tilde{\mathcal P}_\rho(s'|s)$. It holds that
\[
M^\top v=\sum_{s\in\mathcal S, a\in\mathcal A}v(s,a)(e_s-\gamma \mathcal P(\cdot|s,a))=\sum_{s\in\mathcal S}\left(\sum_{a\in\mathcal A}v(s,a)\right)e_s-\gamma\sum_{s\in\mathcal S, a\in\mathcal A}v(s,a)\mathcal P(\cdot|s,a)
\] 
Since $\sum_{a\in\mathcal A}v(s,a)=q(s)$ and $$\sum_{a\in\mathcal A}v(s,a)\mathcal P(\cdot|s,a)=q(s)\sum_{a\in\mathcal A}\pi_\rho(a|s)\mathcal P(\cdot|s,a)=q(s)\bar{\mathcal{P}}_\rho(\cdot|s),$$ we have
\[M^\top v=(I-\gamma \bar P_{\rho}^\top)\ q.\]
Because $\bar P_{\rho}^\top$ is column-stochastic, $\|\bar P_{\rho}^\top q\|_1\le \|q\|_1$ and consequently,  $\|M^\top v\|_1\ge(1-\gamma)\|q\|_1\ge(1-\gamma)\|q\|_2$. Using $\|x\|_2\ge\|x\|_1/\sqrt{|\mathcal S|}$ for $x\in\mathbb{R}^{|\mathcal S|}$, 
$$
\|M^\top v\|_2^2 \ge \frac{1}{|\mathcal S|}\,\|M^\top v\|_1^2
\ge \frac{(1-\gamma)^2}{|\mathcal S|}\|q\|_2^2.
$$
Moreover,
$$
\|v\|_2^2=\sum_{s\in\mathcal S}q(s)^2\sum_{a\in\mathcal A}\Big(\frac{\rho(s,a)}{\tilde\rho(s)}\Big)^2
\le\frac{C^{U^2}}{|\mathcal A|C^{L^2}}{\|c\|_2}^2,
$$
because $\tilde\rho(s)\ge |\mathcal A|C^L$ and $\rho(s,a)\le C^U$. Together,

\begin{equation}
\label{eq:Mv-bound}
    {\|M^\top v\|_2}^2
\ge\frac{(1-\gamma)^2|\mathcal A|}{|\mathcal S|}\frac{C^{L^2}}{C^{U^2}}{\|v\|_2}^2.
\end{equation}
From Equation~\eqref{eq:hessian:f} we have
\[
-h^\top\nabla^2 f(\rho)h
= \frac{1}{\eta_V}\|M^\top(u+v)\|_2^2
    - \eta_\rho\,u^\top \nabla^2 g(\rho)\,u.
\]
For any \(\varepsilon\in(0,1)\), Young's inequality gives
\[
\|a+b\|_2^2 \ge (1-\varepsilon)\|a\|_2^2 - \Big(\frac{1}{\varepsilon}-1\Big)\|b\|_2^2,
\]
so with \(a=M^\top v\), \(b=M^\top u\),
\[
-h^\top\nabla^2 f(\rho)h
\;\ge\;
\frac{1}{\eta_V}\!\Big[
(1-\varepsilon)\|M^\top v\|_2^2
-\Big(\tfrac{1}{\varepsilon}-1\Big)\|M^\top u\|_2^2
\Big]
+\frac{\eta_\rho}{C^U}\|u\|_2^2.
\]
Substituting the operator bounds above yields
\[
-h^\top\nabla^2 f(\rho)h
\;\ge\;
\Big[A
      -B\Big(\tfrac{1}{\varepsilon}-1\Big)\Big]\!\|u\|_2^2
\;+\;
C(1-\varepsilon)\,\|v\|_2^2.
\]

Moreover, for all decompositions \(h=u+v\) we have
\(\|h\|_2^2 \le 2(\|u\|_2^2+\|v\|_2^2)\), hence
\[
-h^\top\nabla^2 f(\rho)h
\;\ge\;
2\,\mu(\varepsilon)\,\|h\|_2^2,
\]
where for every \(\varepsilon\in(0,1)\),
\[
\mu(\varepsilon)
=\frac{1}{2}\min\!\Bigg\{
A - B\Big(\tfrac{1}{\varepsilon}-1\Big),
\;
C(1-\varepsilon)
\Bigg\}.
\]

Maximizing over \(\varepsilon\in(0,1)\) amounts to maximizing \(\mu(\varepsilon)\),
which is achieved by equating the two arguments of the minimum. Solving
\[
A - B\Big(\tfrac{1}{\varepsilon}-1\Big)
\;=\;
C(1-\varepsilon)
\]
for \(\varepsilon\) and inserting back gives the optimal modulus
\[
\mu_{\mathrm{opt}}
=\max_{0<\varepsilon<1}\mu(\varepsilon)
=\frac{1}{4}\Big[(A+B+C)
  -\sqrt{(A+B+C)^2 - 4AC}\Big].
\]
Since \(A>0\) and \(C>0\) by construction, the discriminant
\((A+B+C)^2 - 4AC\) is strictly smaller than \((A+B+C)^2\), so the
bracket is positive and hence \(\mu_{\mathrm{opt}}>0\).
Thus \(-\nabla^2 f(\rho)\succeq \mu I\) on \(H\) for every \(0<\mu\le\mu_{\mathrm{opt}}\),
which proves strong concavity.
\hfill$\blacksquare$

Lastly, the explicit definitions of the deterministic gradient bounds are $$D:=\sqrt{|\mathcal S|}(\eta_V \tfrac{C_r+\eta_\rho U_G}{1-\gamma}+(1+\gamma )|\mathcal X|C^U)$$ and $$B:=\sqrt{|\mathcal S||\mathcal A|}\big(C_r+(1+\gamma)\tfrac{C_r+\eta_\rho U_G}{1-\gamma}+\eta_\rho\log(\tfrac{1}{\pi_{\min}})\big).$$ 
Then for all attainable $(\rho,V,X,\Xi)$ it holds that $\|\hat g(\rho,V,X,\Xi)\|_2\le D$, $\|\nabla_VL(V,\rho)\|_2\le D$, and $\|\nabla_\rho L(V,\rho)\|_2\le B$. This completes the proof of Proposition~\ref{prop: collected-facts-A2}.

\subsection{Proofs of the Propositions~\ref{prop: primal-tracking} and~\ref{prop: dual-variable-error}}
\label{ssec: contraction-proofs}

Given Assumption~\ref{ass:UGE:dual:policies}, we can formulate the limiting update rules by taking expectations under the stationary distribution of the dual-induced Markov chains.

Under the stationary distribution and using asynchronous indicator selection,
\begin{align*}
    \mathbb{E}_{X\sim \mu_\rho}[\widehat g_k(\rho,V;X,\Xi)]
    &= \operatorname{diag}(\tilde \mu_\rho)\,\nabla_V L(V,\rho)\ +\bar{\mathcal{E}}_k(\rho),\\
    \mathbb{E}_{X\sim \mu_\rho}[\widehat h_k(\rho,V;X)]
    &= \operatorname{diag}(\mu_\rho)\,\nabla_\rho L(V,\rho).
\end{align*}

\paragraph{Proof of Proposition~\ref{prop: primal-tracking}}
Fix some $\delta\in(0,1)$ and let $k\ge \max\{K(\delta),\mathcal K\}$. We argue on the high-probability event $\mathcal G_\delta$ in the following. Recall that on $\mathcal G_\delta$ the random stepsizes $\bar\alpha_k,\,\bar\beta_k$ have the following deterministic bounds $$\alpha_k\le\bar\alpha_k\le\alpha_{\mathrm{env},k},\quad\beta_k\le\bar\beta_k\le\beta_{\mathrm{env},k}.$$ 

\emph{(I.)}
To separate the two sources of randomness, namely the update index selection at iteration $k$, $X_k=(s_{k-1},a_{k-1})$ and $s_k$, and the buffer draw $\Xi_k$, we work with the pre-sampling filtration $\mathcal F_k^-$, defined as the $\sigma$-field containing the entire history up to and including the buffer update $\mathcal D_k$, but before sampling $\Xi_k$.
In particular, $\Xi_k$ is conditionally independent of the past given $\mathcal F_k^-$ and satisfies
$\Xi_k \sim \mathcal U(\mathcal D_k(\cdot))$.

Conditionally on $\mathcal F_k^-$, we draw a buffer sample $\Xi_k$ to form the estimator
\[
\widehat G_k:=\widehat g_k(\rho_{k-1},V_{k-1};X_k,\Xi_k).
\]
The shorthand notation of the update is
\[
V_{k+1}=\Pi_{\mathcal V_r}\big[V_{k}-\bar\alpha_{k+1}\,\widehat G_{k+1}\big], 
\]
Recall that $\lambda(\rho):=\arg\min_{V} L(V,\rho).$
By non-expansiveness of $\Pi_{\mathcal V_r}$ and expanding the square,
\begin{equation}
\label{eq: Primal-tracking-one-step-bound}
    \|V_{k+1}-\lambda(\rho_k)\|_2^2 
        \le \|V_k-\lambda(\rho_k)\|_2^2
        -2\bar\alpha_{k+1}\,\big\langle V_k-\lambda(\rho_k),\,\widehat G_{k+1}\big\rangle
        +\bar\alpha_{k+1}^2\,\|\widehat G_{k+1}\|_2^2.
\end{equation}
Taking the conditional expectation to remove the buffer noise, we set
\[
\widetilde G_k:=\mathbb E\big[\widehat G_k| \mathcal F_k^-\big].
\]
Correspondingly, in explicit notation
\[
\tilde g_k(\rho,V;X):=\mathbb E\big[\widehat g_k(\rho,V;X,\Xi)| \mathcal F_k^-\big].
\]
Further, we specify the associated martingale term
\[
\bar M_k:=\widehat G_k-\widetilde G_k,
\]
so that $\mathbb E[\bar M_k| \mathcal F_k^-]=0$ and 
$\mathbb E[\|\bar M_k\|_2^2| \mathcal F_k^-]\le \mathbb E[\|\widehat G_k\|_2^2|\mathcal F_k^-]\le D^2$.
Decompose $\widetilde G_k$ into a target term, a Markovian mismatch (depending only on the law of $X_k$), and the stationary replay bias:
\begin{equation}
\widetilde G_k
=\underbrace{\mathrm{diag}(\tilde\mu_{\rho_{k-1}})\,\nabla_V L(V_{k-1},\rho_{k-1})}_{\text{target at }\rho_{k-1}}
+\underbrace{\Delta^{X}_k}_{\text{mismatch of the index law}}
+\underbrace{\bar{\mathcal E}_k}_{\text{buffer bias at stationarity}},
\end{equation}
where
\begin{align}
\label{eq:index-law-mismatch}
    \Delta^{X}_k&:=\widetilde G_k-\mathbb E_{X\sim \mu_{\rho_{k-1}}}[\tilde g_{k}(\rho_{k-1},V_{k-1};X)],
\end{align}
 and we omit the dependence of the bias term under stationary-law (see~\ref{eq: stationary-law-bias}) on $\rho_{k-1}$ by setting
    $$\bar{\mathcal E}_k:=\bar{\mathcal E}_k(\rho_{k-1})=\mathbb E_{X\sim \mu_{\rho_{k-1}}}[\tilde g_k(\rho_{k-1},V_{k-1};X)]
-\mathrm{diag}(\tilde\mu_{\rho_{k-1}})\,\nabla_V L(V_{k-1},\rho_{k-1}).$$
Therefore,
\[
\widehat G_k
=\mathrm{diag}(\tilde\mu_{\rho_{k-1}})\,\nabla_V L(V_{k-1},\rho_{k-1})
+\Delta^{X}_k+\bar{\mathcal E}_k+\bar M_k.
\]
Taking the conditional expectation $\mathbb E[\cdot| \mathcal F_k^-]$ in the one-step bound~\eqref{eq: Primal-tracking-one-step-bound} removes the martingale term and yields
\begin{align*}
\mathbb E\big[\|V_{k+1}-\lambda(\rho_k)\|_2^2| \mathcal F_k^-\big]
&\le \|V_k-\lambda(\rho_k)\|_2^2
-2\bar\alpha_{k+1}\Big\langle V_k-\lambda(\rho_k),\,\mathrm{diag}(\tilde\mu_{\rho_k})\,\nabla_V L(V_k,\rho_k)\Big\rangle\\
&\quad -2\bar\alpha_{k+1}\Big\langle V_k-\lambda(\rho_k),\,\Delta^{X}_{k+1}\Big\rangle
-2\bar\alpha_{k+1}\Big\langle V_k-\lambda(\rho_k),\,\bar{\mathcal E}_{k+1}\Big\rangle\\
&~~~+\bar\alpha_{k+1}^2\,\mathbb E\big[\|\widehat G_{k+1}\|_2^2| \mathcal F_k^-\big].
\end{align*}
The remaining terms are controlled as follows.

\emph{(II.)}
Due to the $\eta_V$-strong convexity of the mapping $V\mapsto L(V,\rho_k)$, the gradient $\nabla_VL(\cdot,\rho_k)$ is $\eta_V$-strongly monotone and $\nabla_VL(\lambda(\rho_k),\rho_k)=0$. Hence with $\tilde\eta_V=\eta_V\,p_\star\,|\mathcal A|$, it holds that
\[
\langle V_k-\lambda(\rho_k),\mathrm{diag}(\tilde\mu_{\rho_k}) \nabla_VL(V_k,\rho_k)\rangle \;\ge\; \tilde\eta_V\,\|V_k-\lambda(\rho_k)\|_2^2.
\]
Expanding, using the above, and the boundedness of the gradient and its stochastic single entry estimate, $\|\widehat G_k\|_2^2\le D^2,\, \|\nabla_VL\|_2^2\le D^2$, gives
\begin{align*}
\mathbb E\big[\|V_{k+1}-\lambda(\rho_k)\|_2^2| \mathcal F_k^-&\big]
\le 
\|V_k-\lambda(\rho_k)\|_2^2
-2\bar\alpha_{k+1}\,\tilde\eta_V\,\|V_k-\lambda(\rho_k)\|_2^2 \\
& -2\bar\alpha_{k+1}\,\langle V_k-\lambda(\rho_k), \Delta^{X}_{k+1}\rangle-2\bar\alpha_{k+1}\,\langle V_k-\lambda(\rho_k), \bar{\mathcal{E}}_{k+1}\rangle
+\bar\alpha_{k+1}^2\,D^2.
\end{align*}

\emph{(III.)}\label{proof-Prop4_3-stepIII}
 The technical proof of the following bound on $-2\bar\alpha_{k+1}\langle V_k-\lambda(\rho_k),\Delta^{X}_{k+1}\rangle$ applies a comparison of a frozen-chain to the mixing of the inhomogeneous chain and is provided in Lemma~\ref{lemma: Markov-difference}. We obtain the following 
\begin{align*}
       -2\mathbb E\big[\bar{\alpha}_{k+1}\langle V_k-\lambda(\rho_k),\Delta^{X}_{k+1}\rangle | \mathcal G_\delta\big] &\le 2\mathbb E\big[\bar{\alpha}_{k+1}|\langle V_k-\lambda(\rho_k),\Delta^{X}_{k+1}\rangle ||\mathcal G_\delta\big]\\
       &\le C_1\tau_k^2\alpha_{env,k}\alpha_{env,k-\tau_k}\mathbb E\big[\|V_{k-\tau_k}-\lambda(\rho_{k-\tau_k})\|_2^2+1|\mathcal G_\delta\big],
\end{align*}
where $C_1=2(LDB+LB+LD+D^2+\tfrac{\gamma |\mathcal X|C^U}{\beta_0})$.

To address the replay-buffer bias, we apply the inequality $\langle 2\upsilon_1,\upsilon_2\rangle\le c\|\upsilon_1\|_2^2+\tfrac 1c\|\upsilon_2\|_2^2$ that holds for any vectors $\upsilon_1,\upsilon_2$ and $c>0$ with $\upsilon_1=-(V_k-\lambda(\rho_k))$, $\upsilon_2=\bar{\mathcal E}_{k+1}$ and $c=\tilde\eta_V$ to get
\begin{align*}
&-2\mathbb E\big[\bar\alpha_{k+1}\langle V_k-\lambda(\rho_k),\bar{\mathcal E}_{k+1}\rangle|\mathcal F_k^-,\mathcal G_\delta\big] \\
&\qquad \qquad \le
\tilde\eta_V\mathbb E\big[\bar\alpha_{k+1}\|V_k-\lambda(\rho_k)\|_2^2|\mathcal F_k^-,\mathcal G_\delta\big]
+\tfrac{1}{\tilde\eta_V}\,\mathbb E\big[\bar\alpha_{k+1}\|\bar{\mathcal E}_{k+1}\|_2^2|\mathcal F_k^-,\mathcal G_\delta\big]\\
& \qquad \qquad \le
\tilde\eta_V\mathbb E\big[\bar\alpha_{k+1}\|V_k-\lambda(\rho_k)\|_2^2|\mathcal F_k^-,\mathcal G_\delta\big]
+\tfrac{C_{\mathrm{buf}}(\delta)^2|\mathcal S|}{\tilde\eta_V}\tfrac{\alpha_{env,k}}{k},
\end{align*}
where the last inequality uses the bound on the bias under stationary law derived in Lemma~\ref{lemma:buffer-bias-HP}.

Recall the deterministic bounds on the random stepsize sequence $\bar\alpha_k$, $\alpha_k\le\bar\alpha_k\le\alpha_{\mathrm{env},k}$. Putting the bounds together and taking the conditional expectation gives us
\begin{align*}
    \mathbb E\big[\|V_{k+1}-\lambda(\rho_k)\|_2^2|\mathcal G_\delta\big]&\le
    (1-\tilde\eta_V\alpha_k)\mathbb E[\|V_k-\lambda(\rho_k)\|_2^2|\mathcal G_\delta]\\
    &~~~~~+C_1\tau_k^2\alpha_{\mathrm{env},k-\tau_k}\alpha_{\mathrm{env},k}(\mathbb E[\|V_{k-\tau_k}-\lambda(\rho_{k-\tau_k})\|_2^2| \mathcal G_\delta]+1)\\
    &~~~~~+\tfrac{C_{\mathrm{buf}}(\delta)^2|\mathcal S|}{\tilde\eta_V}\tfrac{\alpha_{\mathrm{env},k}}{k}+D^2\alpha_{\mathrm{env},k}^2
\end{align*}

\emph{(IV.)}
To obtain the desired contraction property of $\|V_{k+1}-\lambda(\rho_{k+1})\|_2^2$ we use that for any $c>0$
\[
\|V_{k+1}-\lambda(\rho_{k+1})\|_2^2
\le 
(1+c)\,\|V_{k+1}-\lambda(\rho_k)\|_2^2
+(1+\frac{1}{c})\,\|\lambda(\rho_{k+1})-\lambda(\rho_k)\|_2^2.
\]
With the Lipschitz continuity of the best response (see Lemma~\ref{lemma: LipschitzLambda}) and the bound on the dual update, we have 
\[
    \|\lambda(\rho_{k+1})-\lambda(\rho_k)\|_2^2\le L^2\|\rho_{k+1}-\rho_k\|_2^2\le L^2 B^2\beta_{\mathrm{env,k}}^2.
\]
With the adaptive choice of the inequality $c_k=\tfrac{\tilde\eta_V}{2}\alpha_{k}$ and the shorthand notation $y_k=\mathbb E\big[\|V_k-\lambda(\rho_k)\|_2^2|\mathcal G_\delta\big]$, we have the combined bound
\begin{align*}
    y_{k+1}\;\le&\;(1+\tfrac{\tilde\eta_V}{2}\alpha_{k})(1-\tilde\eta_V\alpha_{k})y_k
    +(1+\tfrac{\tilde\eta_V}{2}\alpha_{k}) C_1\tau_k^2\alpha_{\mathrm{env},k-\tau_k}\alpha_{\mathrm{env},k}(y_{k-\tau_k}+1)\\
    &~+(1+\tfrac{\tilde\eta_V}{2}\alpha_{k})\tfrac{C_{\mathrm{buf}}(\delta)^2|\mathcal S|}{\tilde\eta_V}\tfrac{\alpha_{\mathrm{env},k}}{k}+
    (1+\tfrac{\tilde\eta_V}{2}\alpha_{k})D^2\alpha_{\mathrm{env},k}^2\\
    &~+L^2B^2\beta_{\mathrm{env},k}^2+\tfrac{2L^2B^2\beta_{\mathrm{env},k}^2}{\tilde\eta_V\,\alpha_{k}}
\end{align*}
To simplify the above bound, we use $(1+\tfrac{\tilde\eta_V}{2}\alpha_{k})(1-\tilde\eta_V\alpha_{k})\;\le\;1-\tfrac{\tilde\eta_V}{2}\alpha_{k}$ and the following bound on $y_{k-\tau_k}$ in the second summand
\begin{align*}
    y_{k-\tau_k}&=\mathbb E\big[\|V_k-\lambda(\rho_k)-(V_k-V_{k-\tau_k})+\big(\lambda(\rho_k)-\lambda(\rho_{k-\tau_k})\big)\|_2^2|\mathcal G_\delta\big]\\
        &\le 3(y_k+\mathbb E[\|V_k-V_{k-\tau_k}\|_2^2|\mathcal G_\delta]+\mathbb E[\|\lambda(\rho_k)-\lambda(\rho_{k-\tau_k})\|_2^2|\mathcal G_\delta])\\
        &\le 3y_k+3D^2\tau_k^2\alpha_{\mathrm{env},k-\tau_k}^2+3L^2B^2\tau_k^2\beta_{\mathrm{env},k-\tau_k}^2,
\end{align*}
where the second inequality is a consequence of the Cauchy-Schwarz inequality, and the last inequality follows from Lipschitz continuity of the best response and the bounded changes to the variables over the mixing interval. Using the condition $\tfrac{\tilde\eta_V}{2}\alpha_k\le\tfrac{1}{2}$ for $k\ge\max\{K(\delta),\mathcal K\}$ the inequality simplifies to
\begin{align*}
    y_{k+1}\;\le&\;(1-\tfrac{\tilde\eta_V}{2}\alpha_{k})y_k
    + \tfrac{9}{2}C_1\tau_k^2\alpha_{\mathrm{env},k-\tau_k}\alpha_{\mathrm{env},k}y_k\\
    &~+\tfrac{9}{2}C_1\tau_k^2\alpha_{\mathrm{env},k-\tau_k}\alpha_{\mathrm{env},k}(D^2\tau_k^2\alpha_{\mathrm{env},k-\tau_k}^2+L^2B^2\tau_k^2\beta_{\mathrm{env},k-\tau_k}^2+\tfrac13)\\
    &~+\tfrac{3}{2}\tfrac{C_{\mathrm{buf}}(\delta)^2|\mathcal S|}{\tilde\eta_V}\tfrac{\alpha_{\mathrm{env},k}}{k}+
    \tfrac{3}{2}D^2\alpha_{\mathrm{env},k}^2+L^2B^2\beta_{\mathrm{env},k}^2+\tfrac{2L^2B^2\beta_{\mathrm{env},k}^2}{\tilde\eta_V\,\alpha_{k}}
\end{align*}
By the definition of $\mathcal K$ in~\eqref{eq: K_ZDR}, for all $k\ge\max\{K(\delta),\mathcal K\}$, $C_1C_\mathrm{env}\tau_k^2\alpha_{\mathrm{env},k-\tau_k}\le\tfrac{\tilde\eta_V}{18}$ and $\tau_k\alpha_{\mathrm{env},k-\tau_k}\le\min\{\tfrac{1}{\sqrt{3}D},\,\tfrac{1}{\sqrt{3}LB}\}$. Inserting the bounds simplifies the inequality to 
\begin{equation}
    \label{eq: aux-var-contr}
    \begin{aligned}
    y_{k+1}\;\le&\;(1-\tfrac{\tilde\eta_V}{4}\alpha_{k})y_k
    +\tfrac{9}{2}C_1\tau_k^2\alpha_{\mathrm{env},k-\tau_k}\alpha_{\mathrm{env},k}+\tfrac{3}{2}\tfrac{C_{\mathrm{buf}}(\delta)^2|\mathcal S|}{\tilde\eta_V}\tfrac{\alpha_{\mathrm{env},k}}{k}\\
    &~+
    \tfrac{3}{2}D^2\alpha_{\mathrm{env},k}^2+L^2B^2\beta_{\mathrm{env},k}^2+\tfrac{2L^2B^2}{\tilde\eta_V}\tfrac{\beta_{\mathrm{env},k}^2}{\alpha_k},
    \end{aligned}
\end{equation}
which finishes the primal variable contraction proof.
\hfill $\blacksquare$

\paragraph{Proof of Proposition~\ref{prop: dual-variable-error}}
    Fix \(\delta\in(0,1)\) and let \(k\ge \max\{K(\delta),\mathcal K\}\). We work on the
    high-probability event \(\mathcal G_\delta\). Since \(\tau_k\ge 1\), the burn-in condition
    \(
    \tau_k\alpha_{\mathrm{env},k-\tau_k}\le \tfrac{1}{\sqrt{3}LB}
    \)
    implies
    \(
    \alpha_{\mathrm{env},k-\tau_k}\le \tfrac{1}{\sqrt{3}LB}.
    \)
    For the chosen two-timescale stepsizes, \(\beta_t\le \alpha_t\) for all sufficiently large \(t\), hence
    \[
    \bar\beta_k\le \beta_{\mathrm{env},k}\le \beta_{\mathrm{env},k-\tau_k}\le \alpha_{\mathrm{env},k-\tau_k}\le \frac{1}{\sqrt{3}LB}\le\frac{1}{L},
    \]
    where the last inequality follows from \(B\ge 1\).
    Similarly to the proof of Proposition~\ref{prop: primal-tracking}, we set
    \[
    \widehat H_{k}:=\widehat h_{k}(\rho_{k-1},V_{k-1};X_k),
    \qquad
    \widetilde H_{k}:=\mathbb E[\widehat H_{k}\mid \mathcal F_{k-1}].
    \]
    Then, by taking the expectation under the stationary distribution in~\eqref{eq:async:rho:recursion}, we have
    \[
    \mathbb E_{X\sim\mu_{\rho_{k-1}}}\Big[\widehat h_{k}(\rho_{k-1},V_{k-1};X)\Big]=\operatorname{diag}(\mu_{\rho_{k-1}})\nabla_\rho L(V_{k-1},\rho_{k-1}).
    \]
    Using this stationary mean-field identity and \(\nabla f(\rho)=\nabla_\rho L(\lambda(\rho),\rho)\), we decompose
    \[
    \widetilde H_{k+1}
    =
    D(\rho_k)\nabla f(\rho_k)
    +\Delta^V_{k+1}
    +\bar\Delta^X_{k+1},
    \]
    where
    \[
    D(\rho_k):=\operatorname{diag}(\mu_{\rho_k}),\quad
    \Delta^V_{k+1}
    :=
    \mathbb E_{X\sim \mu_{\rho_k}}
    \Big[
    \widehat h_{k+1}(\rho_k,V_k;X)-\widehat h_{k+1}(\rho_k,\lambda(\rho_k);X)
    \Big],
    \]
    and
    \[
    \bar\Delta^X_{k+1}
    :=
    \widetilde H_{k+1}
    -
    \mathbb E_{X\sim \mu_{\rho_k}}
    \big[
    \widehat h_{k+1}(\rho_k,\lambda(\rho_k);X)
    \big].
    \]
    To isolate the core finite-time contraction mechanism, we carry out the following one-step estimate in the regime where the clipping is inactive after burn-in. 
    By \(L\)-smoothness of \(f\), 
    \begin{align*}
    f(\rho_{k+1})
    &\ge
    f(\rho_k)
    +
    \bar\beta_{k+1}
    \langle \nabla f(\rho_k),\widehat H_{k+1}\rangle
    -
    \frac{L}{2}\bar\beta_{k+1}^2\|\widehat H_{k+1}\|_2^2.
    \end{align*}
    By Proposition~\ref{prop: collected-facts-A2}(iv), \(\|\widehat H_{k+1}\|^2_2\le B^2\). 
    Substituting the decomposition of \(\widetilde H_{k+1}\), we obtain
    \begin{align*}
    f(\rho_{k+1})
    &\ge
    f(\rho_k)
    +
    \bar\beta_{k+1}
    \langle \nabla f(\rho_k),D(\rho_k)\nabla f(\rho_k)\rangle
    \\
    &\qquad
    +
    \bar\beta_{k+1}
    \langle \nabla f(\rho_k),\Delta^V_{k+1}+\bar\Delta^X_{k+1}+(\widehat H_{k+1}-\widetilde H_{k+1})\rangle
    +
    \frac{LB^2}{2}\bar\beta_{k+1}^2 .
    \end{align*}
    Rearranging gives
    \begin{align}
    f(\rho^\star)-f(\rho_{k+1})
    &\le
    f(\rho^\star)-f(\rho_k)-\bar\beta_{k+1}
    \langle \nabla f(\rho_k),D(\rho_k)\nabla f(\rho_k)\rangle\nonumber
    \\
    &\qquad
    -
    \bar\beta_{k+1}
    \langle \nabla f(\rho_k),\Delta^V_{k+1}+\bar\Delta^X_{k+1}+(\widehat H_{k+1}-\widetilde H_{k+1})\rangle
    +
    \frac{L B^2}{2}\bar\beta_{k+1}^2.
    \label{eq:prop44-intermediate}
    \end{align}
    Using \(\mu_{\rho_k}(s,a)\ge p_\star\) for all \((s,a)\in\mathcal X\) 
    \begin{align}
    -\langle \nabla f(\rho_k),D(\rho_k)\nabla f(\rho_k)\rangle 
    &\le -p_\star\|\nabla f(\rho_k) \|_2^2\le-2p_\star\mu_{\mathrm{opt}}\bigl(f(\rho^\star)-f(\rho_k)\bigr).
    \label{eq:prop44-intermediate-boundA}
    \end{align}
    By Lipschitz continuity of
    \(\widehat h_k(\rho,\cdot;X)\) in \(V\),
    \(
    \|\Delta^V_{k+1}\|_2
    \le
    L\|V_k-\lambda(\rho_k)\|_2.
    \)
    Using Young's inequality with parameter \(p_\star\),
    \[
    -\langle \nabla f(\rho_k),\Delta^V_{k+1}\rangle
    \le
    \frac{p_\star}{2}\|\nabla f(\rho_k)\|_2^2
    +
    \frac{1}{2p_\star}\|\Delta^V_{k+1}\|_2^2,
    \]
    and therefore
    \begin{align}
    -\bar\beta_{k+1}\langle \nabla f(\rho_k),\Delta^V_{k+1}\rangle
    \le
    \frac{p_\star}{2}\bar\beta_{k+1}\|\nabla f(\rho_k)\|_2^2
    +
    \frac{L^2}{2p_\star}\bar\beta_{k+1}\|V_k-\lambda(\rho_k)\|_2^2.
    \label{eq:prop44-intermediate-boundB}
    \end{align}
    By the Polyak--\L ojasiewicz inequality for \(\mu_{\mathrm{opt}}\)-strongly concave \(f\),
    \[
    \tfrac{1}{2}\|\nabla f(\rho)\|_2^2
    \ge
    \mu_{\mathrm{opt}}\bigl(f(\rho^\star)-f(\rho)\bigr).
    \] 
    Combining~\eqref{eq:prop44-intermediate-boundB} with~\eqref{eq:prop44-intermediate-boundA} and the Polyak--\L ojasiewicz inequality and inserting the bounds into Equation~\eqref{eq:prop44-intermediate} gives
    \begin{align}
    f(\rho^\star)-f(\rho_{k+1})
    &\le
    (1-p_\star \mu_{\mathrm{opt}}\bar\beta_{k+1})\bigl(f(\rho^\star)-f(\rho_k)\bigr)
    +\frac{L^2}{2p_\star}
    \bar\beta_{k+1}\|V_k-\lambda(\rho_k)\|_2^2\nonumber\\
    &\qquad
    -
    \bar\beta_{k+1}
    \langle \nabla f(\rho_k),\bar\Delta^X_{k+1}+(\widehat H_{k+1}-\widetilde H_{k+1})\rangle
    +
    \frac{L B^2}{2}\bar\beta_{k+1}^2.
    \label{eq:prop44-main-step-unprojected}
    \end{align}
    To address the Martingale difference term \(\widehat H_{k+1}-\widetilde H_{k+1}\), note that 
    \[
    \mathbb E[\langle\nabla f(\rho_k),\widehat H_{k+1}-\widetilde H_{k+1}\rangle|\mathcal F_{k}]=0
    \]
    Taking expectation conditional on \(\mathcal G_\delta\), using
    \(\beta_k\le\bar\beta_{k+1}\le \beta_{\mathrm{env},k}\), and invoking
    the Markovian mismatch bound from Lemma~\ref{lemma:dual-markov-mismatch},
    \[
    -\mathbb E\!\left[
    \langle \nabla f(\rho_k),\bar\Delta^X_{k+1}\rangle
    |
    \mathcal G_\delta
    \right]
    \le
    12L^2B^3\tau_k^2\beta_{\mathrm{env},k-\tau_k},
    \]
    we arrive at
    \begin{align*}
    \mathbb E[f(\rho^\star)-f(\rho_{k+1})| \mathcal G_\delta]
    &\le
    \bigl(1-p_\star\mu_{\mathrm{opt}}\beta_k\bigr)
    \mathbb E[f(\rho^\star)-f(\rho_k)| \mathcal G_\delta]
    +
    \frac{L^2\beta_{\mathrm{env},k}}{2p_\star}
    \mathbb E[\|V_k-\lambda(\rho_k)\|_2^2| \mathcal G_\delta]
    \\
    &\qquad
    +
    12L^2B^3\tau_k^2\beta_{\mathrm{env},k-\tau_k}\beta_{\mathrm{env},k}
    +
    \frac{L}{2}B^2\beta_{\mathrm{env},k}^2.
    \end{align*}
    The tower property with \(\mathcal G_\delta\) eliminates the Martingale difference term.
    Finally, since \(\tau_k\ge 1\) and \(\beta_{\mathrm{env},k-\tau_k}\ge \beta_{\mathrm{env},k}\),
    \(
    \beta_{\mathrm{env},k}^2
    \le
    \tau_k^2\beta_{\mathrm{env},k-\tau_k}\beta_{\mathrm{env},k},
    \)
    so the last term can be absorbed into the same remainder term. Thus, 
    \begin{align*}
    \mathbb E[f(\rho^\star)-f(\rho_{k+1})| \mathcal G_\delta]
    &\le
    \bigl(1-p_\star\mu_{\mathrm{opt}}\beta_k\bigr)
    \mathbb E[f(\rho^\star)-f(\rho_k)| \mathcal G_\delta]
    +
    \frac{L^2}{2p_\star}\beta_{\mathrm{env},k}
    \mathbb E[\|V_k-\lambda(\rho_k)\|_2^2| \mathcal G_\delta]
    \\
    &\qquad
    +
    \frac{25L^2B^3}{2}
    \tau_k^2\beta_{\mathrm{env},k-\tau_k}\beta_{\mathrm{env},k},
    \end{align*}
    which is the claimed recursion.
\hfill $\blacksquare$

\subsection{Auxiliary Results and Proofs} \label{ssec:Aux_results:convergence:rate}
This section contains additional auxiliary results required for the convergence rate proof as well as the proof of Corollary~\ref{cor:value-suboptimality}.

We first show the Lipschitz-continuity of the mapping $\pi\mapsto V^{\pi}_{ur}$. 
\begin{lemma}[Lipschitz continuity of the unregularized value in the policy]
\label{lemma:value-lipschitz-policy}
For \\a stationary policy \(\pi\in\Pi\), let \(V^\pi_{ur}\) denote the unregularized discounted
value function. Then, for any two stationary policies \(\pi,\hat\pi\in\Pi\),
\[
\|V^\pi_{ur}-V^{\hat\pi}_{ur}\|_\infty
\le
\frac{C_r}{(1-\gamma)^2}
\max_{s\in\mathcal S}\|\pi(\cdot|s)-\hat\pi(\cdot|s)\|_1.
\]
\end{lemma}
\noindent
\textbf{Proof}\hspace{2pt}
For a stationary policy \(\pi\), define
\[
r_\pi(s):=\sum_{a\in\mathcal A}\pi(a|s)r(s,a),
\qquad
[\mathcal P_\pi]_{ss'}
:=
\sum_{a\in\mathcal A}\pi(a|s)\mathcal P(s'|s,a).
\]
Then, \(V^\pi_{ur}\) satisfies the Bellman equation
\[
V^\pi_{ur}=r_\pi+\gamma \mathcal P_\pi V^\pi_{ur}.
\]
Subtracting the Bellman equations for \(\pi\) and \(\hat\pi\) yields
\[
V^\pi_{ur}-V^{\hat\pi}_{ur}
=
\gamma \mathcal P_\pi\bigl(V^\pi_{ur}-V^{\hat\pi}_{ur}\bigr)
+
(r_\pi-r_{\hat\pi})
+
\gamma(\mathcal P_\pi-\mathcal P_{\hat\pi})V^{\hat\pi}_{ur}.
\]
Hence,
\[
V^\pi_{ur}-V^{\hat\pi}_{ur}
=
(I-\gamma \mathcal P_\pi)^{-1}
\Bigl[
(r_\pi-r_{\hat\pi})
+
\gamma(\mathcal P_\pi-\mathcal P_{\hat\pi})V^{\hat\pi}_{ur}
\Bigr].
\]
Since \(\mathcal P_\pi\) is a stochastic matrix, we have
\(\|\mathcal P_\pi^t v\|_\infty\le \|v\|_\infty\) for all \(t\ge 0\) and
\(v\in\mathbb R^{|\mathcal S|}\). Hence, using the Neumann-series expansion,
\[
(I-\gamma \mathcal P_\pi)^{-1}v
=
\sum_{t=0}^\infty \gamma^t \mathcal P_\pi^t v,
\]
it follows that
\[
\|(I-\gamma \mathcal P_\pi)^{-1}v\|_\infty
\le
\sum_{t=0}^\infty \gamma^t \|\mathcal P_\pi^t v\|_\infty
\le
\sum_{t=0}^\infty \gamma^t \|v\|_\infty
=
\frac{1}{1-\gamma}\|v\|_\infty,
\qquad \forall\, v\in\mathbb R^{|\mathcal S|}.
\]
See also~\citet[Theorem 6.1.1]{puterman1994markov} for the corresponding discounted-value representation.
Moreover, for every \(s\in\mathcal S\),
\[
|r_\pi(s)-r_{\hat\pi}(s)|
=
\left|
\sum_{a\in\mathcal A}(\pi(a|s)-\hat\pi(a|s))r(s,a)
\right|
\le
C_r \|\pi(\cdot|s)-\hat\pi(\cdot|s)\|_1,
\]
and therefore
\[
\|r_\pi-r_{\hat\pi}\|_\infty
\le
C_r \max_{s\in\mathcal S}\|\pi(\cdot|s)-\hat\pi(\cdot|s)\|_1.
\]
Likewise, for any \(v\in\mathbb R^{|\mathcal S|}\) and \(s\in\mathcal S\),
\begin{align*}
|((\mathcal P_\pi-\mathcal P_{\hat\pi})v)(s)|
&=
\left|
\sum_{a\in\mathcal A}(\pi(a|s)-\hat\pi(a|s))
\sum_{s'\in\mathcal S}\mathcal P(s'|s,a)v(s')
\right| \\
&\le
\|v\|_\infty \|\pi(\cdot|s)-\hat\pi(\cdot|s)\|_1.
\end{align*}
Hence,
\[
\|(\mathcal P_\pi-\mathcal P_{\hat\pi})v\|_\infty
\le
\|v\|_\infty
\max_{s\in\mathcal S}\|\pi(\cdot|s)-\hat\pi(\cdot|s)\|_1.
\]
Since \(0\le r(s,a)\le C_r\), we also have
\[
\|V^{\hat\pi}_{ur}\|_\infty \le \frac{C_r}{1-\gamma}.
\]
Combining the above bounds gives
\begin{align*}
\|V^\pi_{ur}-V^{\hat\pi}_{ur}\|_\infty
&\le
\frac{1}{1-\gamma}
\left(
\|r_\pi-r_{\hat\pi}\|_\infty
+
\gamma \|(\mathcal P_\pi-\mathcal P_{\hat\pi})V^{\hat\pi}_{ur}\|_\infty
\right) \\
&\le
\frac{1}{1-\gamma}
\left(
C_r + \gamma \frac{C_r}{1-\gamma}
\right)
\max_{s\in\mathcal S}\|\pi(\cdot|s)-\hat\pi(\cdot|s)\|_1 \\
&=
\frac{C_r}{(1-\gamma)^2}
\max_{s\in\mathcal S}\|\pi(\cdot|s)-\hat\pi(\cdot|s)\|_1,
\end{align*}
which proves the claim.
\hfill$\blacksquare$

\paragraph{Proof of Corollary~\ref{cor:value-suboptimality}.}
By Proposition~\ref{thm:Li-et-al}, the optimal dual solution \(\rho^\star\) induces the optimal
regularized policy, that is,
\[
\pi_{\rho^\star}=\pi_r^\star.
\]
By the definition of \(L_\pi\) in Corollary~\ref{cor:value-suboptimality}, whose existence
follows directly from the Lipschitz estimate derived in the proof of
Lemma~\ref{lemma:B7}, we have
\[
\max_{s\in\mathcal S}
\|\pi_{\rho_k}(\cdot|s)-\pi_{\rho^\star}(\cdot|s)\|_1
\le
L_\pi \|\rho_k-\rho^\star\|_2.
\]
Applying Lemma~\ref{lemma:value-lipschitz-policy} with
\(\pi=\pi_{\rho_k}\) and \(\hat\pi=\pi_{\rho^\star}=\pi_r^\star\), we obtain
\[
\|V^{\pi_{\rho_k}}_{ur}-V^{\pi_r^\star}_{ur}\|_\infty
\le
\frac{C_r}{(1-\gamma)^2}
\max_{s\in\mathcal S}
\|\pi_{\rho_k}(\cdot|s)-\pi_{\rho^\star}(\cdot|s)\|_1
\le
L_V\|\rho_k-\rho^\star\|_2.
\]
Next, by Proposition~\ref{prop: optdual-reg} and the specialization
\(0\le G(\pi(\cdot|s))\le \log|\mathcal A|\) for the Shannon entropy regularizer, it holds for every
\(s\in\mathcal S\) that
\[
0
\le
V^\star_{ur}(s)-V^{\pi_r^\star}_{ur}(s)
\le
\frac{\eta_\rho\log|\mathcal A|}{1-\gamma}
=
B_{\mathrm{reg}}.
\]
Therefore,
\begin{align*}
V^\star_{ur}(s)-V^{\pi_{\rho_k}}_{ur}(s)
&=
\bigl(V^\star_{ur}(s)-V^{\pi_r^\star}_{ur}(s)\bigr)
+
\bigl(V^{\pi_r^\star}_{ur}(s)-V^{\pi_{\rho_k}}_{ur}(s)\bigr) \\
&\le
B_{\mathrm{reg}}
+
\|V^{\pi_r^\star}_{ur}-V^{\pi_{\rho_k}}_{ur}\|_\infty.
\end{align*}
Since \(V^\star_{ur}(s)-V^{\pi_{\rho_k}}_{ur}(s)\ge 0\), this implies
\[
0\le\left(
V^\star_{ur}(s)-V^{\pi_{\rho_k}}_{ur}(s)-B_{\mathrm{reg}}
\right)_+
\le
\|V^{\pi_r^\star}_{ur}-V^{\pi_{\rho_k}}_{ur}\|_\infty
\le
L_V\|\rho_k-\rho^\star\|_2.
\]
Squaring both sides and taking conditional expectations on \(\mathcal G_\delta\) yields
\[
\mathbb E\!\left[
\left(
V^\star_{ur}(s)-V^{\pi_{\rho_k}}_{ur}(s)-B_{\mathrm{reg}}
\right)_+^2
\;\middle|\;
\mathcal G_\delta
\right]
\le
L_V^2
\mathbb E\!\left[
\|\rho_k-\rho^\star\|_2^2
\;\middle|\;
\mathcal G_\delta
\right].
\]
The final rate statement follows immediately from Theorem~\ref{thm:B1}.
\hfill$\blacksquare$

\begin{lemma}[Direct bound on $\tau_{\mathrm{inh}}(\zeta)$ and $\tau_{\min}$]
\label{lemma: inhom-mixing-time-bound}
Let Assumption~\ref{ass:UGE:dual:policies}(ii) hold. Let $\tau_{\min}$ be
as defined in~\eqref{eq:Paulin-tail} and let $\tau_{\mathrm{inh}}(\zeta)$ be the
inhomogeneous mixing time of $\{P_{\rho_k}\}$ as defined in~\eqref{eq:inh-mixing}.
Then, for any $\zeta\in(0,1)$,
$$\tau_{\mathrm{inh}}(\zeta)
\le
(m_\star-1)
+
m_\star\left\lceil\frac{\log(1/\zeta)}{\log(1/\kappa)}
\right\rceil,$$
and consequently
$$
\tau_{\min}=\inf_{0\le \zeta<1}\tau_{\mathrm{inh}}(\zeta)\Bigl(\frac{2-\zeta}{1-\zeta}\Bigr)^2
\le
9\Biggl[(m_\star-1)+m_\star\left\lceil
\frac{\log 2}{\log(1/\kappa)}\right\rceil\Biggr].
$$
\end{lemma}
\textbf{Proof}\hspace{2pt}
By Assumption~\ref{ass:UGE:dual:policies}(ii), there exist constants
$m_\star\in\mathbb N$ and $\kappa\in(0,1)$ such that for every sequence
$\{\rho_t\}_{t\in\mathbb N}\subset H$,
\[
\sup_{k\ge 0}
\delta\!\left(
\prod_{t=k}^{k+m_\star-1} P_{\rho_t}
\right)
\le
\kappa.
\]
Fix an initial time index $i\ge 0$, and let $k=q\,m_\star+r$ with
$q\in\mathbb N_0$ and $0\le r<m_\star$. Write the product of $k$ consecutive
kernels as
\[
\prod_{t=i}^{i+k-1} P_{\rho_t}
=
\left(
\prod_{b=0}^{q-1}
\prod_{t=i+b m_\star}^{i+(b+1)m_\star-1} P_{\rho_t}
\right)
\left(
\prod_{t=i+q m_\star}^{i+q m_\star+r-1} P_{\rho_t}
\right).
\]
Using the submultiplicativity of Dobrushin's ergodic coefficient,
\[
\delta(Q_1Q_2)\le \delta(Q_1)\delta(Q_2),
\]
together with the bound $\delta(Q)\le 1$ for every Markov kernel $Q$, we obtain
\begin{align*}
\delta\!\left(
\prod_{t=i}^{i+k-1} P_{\rho_t}
\right)
&\le
\prod_{b=0}^{q-1}
\delta\!\left(
\prod_{t=i+b m_\star}^{i+(b+1)m_\star-1} P_{\rho_t}
\right)
\cdot
\delta\!\left(
\prod_{t=i+q m_\star}^{i+q m_\star+r-1} P_{\rho_t}
\right)
\\
&\le
\kappa^q.
\end{align*}
Hence, for all $x,y\in\mathcal X$,
\[
\left\|
\delta_x \prod_{t=i}^{i+k-1} P_{\rho_t}
-
\delta_y \prod_{t=i}^{i+k-1} P_{\rho_t}
\right\|_{\mathrm{TV}}
\le
\delta\!\left(
\prod_{t=i}^{i+k-1} P_{\rho_t}
\right)
\le
\kappa^q.
\]
Therefore, if $\kappa^q\le \zeta$, then any product of length
$k=q\,m_\star+r$ contracts total variation distance by at most $\zeta$.
Choosing
\[
q
=
\left\lceil
\frac{\log(1/\zeta)}{\log(1/\kappa)}
\right\rceil
\]
yields
\[
\tau_{\mathrm{inh}}(\zeta)
\le
(m_\star-1)
+
m_\star
\left\lceil
\frac{\log(1/\zeta)}{\log(1/\kappa)}
\right\rceil.
\]

For $\tau_{\min}$, we proceed exactly as in the current argument and evaluate at
$\zeta=\tfrac12$:
\[
\tau_{\min}
=
\inf_{0\le \zeta<1}
\tau_{\mathrm{inh}}(\zeta)
\Bigl(\frac{2-\zeta}{1-\zeta}\Bigr)^2
\le
\tau_{\mathrm{inh}}(1/2)\Bigl(\frac{2-1/2}{1-1/2}\Bigr)^2
=
9\,\tau_{\mathrm{inh}}(1/2).
\]
Combining this with the previous bound gives
\[
\tau_{\min}
\le
9\Biggl[
(m_\star-1)
+
m_\star
\left\lceil
\frac{\log 2}{\log(1/\kappa)}
\right\rceil
\Biggr],
\]
which proves the claim.
\hfill $\blacksquare$

The Lipschitz property of the dual-induced transition kernel has been established in the proof of Lemma~\ref{lemma:B7}. We next provide a proof of the Lipschitz continuity of the dual-induced stationary distribution which is applied in the proof of the visitation floor lemma. 

\begin{lemma}[Lipschitz stationary law under Assumption~\ref{ass:UGE:dual:policies}]\hspace{2pt}\label{lemma:Lip-stationary}
Let $d,\hat d$ be probability distributions on $\mathcal X$. There exists a constant $L_\mu<\infty$ (depending only on $C_{\mathcal X}$ and $\varrho$) such that for all $\rho,\hat\rho\in H$,
\begin{equation}
\label{eq:stationary-perturbation}
\|\mu_\rho-\mu_{\hat\rho}\|_{\mathrm{TV}}
\;\le\; L_\mu\,\|\mathcal P_\rho-\mathcal P_{\hat\rho}\|_{\infty}
\;\le\; L_\mu\,C_{\mathcal Q}\,
\max_{s'\in\mathcal S}\sum_{a'\in\mathcal A}\big|\rho(s',a')-\hat\rho(s',a')\big|.
\end{equation}
Moreover, for all $\rho,\hat\rho\in H$,
\begin{equation}
\label{eq:kernel-perturbation}
\| d\mathcal P_\rho - \hat d\mathcal P_{\hat\rho}\|_{\mathrm{TV}}
\le \|d-\hat d\|_{\mathrm{TV}}
+ \|\mathcal P_\rho-\mathcal P_{\hat\rho}\|_{\infty}.
\end{equation}
\end{lemma}
\textbf{Proof}\hspace{2pt}
Under Assumption~\ref{ass:UGE:dual:policies}, for each $\rho\in H$ there is a unique stationary distribution $\mu_\rho$ and applying \cite[Corollary 3.1]{ref:mitrophanov2005:sensitivity} we have

$$
\|\mu_\rho-\mu_{\hat\rho}\|_{\mathrm{TV}}
\le L_\mu \|\mathcal P_\rho-\mathcal P_{\hat\rho}\|_{\infty},
$$
where $L_\mu:=\log_\varrho(C_\mathcal{X}^{-1})+\frac{1}{1-\varrho}$. It remains to control $\|\mathcal P_\rho-\mathcal P_{\hat\rho}\|_{\infty}$ in terms of $\rho-\hat\rho$, which has already been shown in the proof of Lemma~\ref{lemma:B7}, where the Lipschitz property of the $\mathcal P_\rho$ has been established. This shows~\eqref{eq:stationary-perturbation}. 
The perturbation bound~\eqref{eq:kernel-perturbation} follows from the triangle inequality applied as 
\[
\| d\mathcal P_\rho - \hat d\mathcal P_{\hat\rho} \|_{\mathrm{TV}}
\le \| d\mathcal P_\rho - \hat d\mathcal P_\rho \|_{\mathrm{TV}}
+ \| \hat d\mathcal P_\rho - \hat d\mathcal P_{\hat\rho} \|_{\mathrm{TV}}.
\]
For the first term, total variation is non-expansive under any fixed Markov kernel $K$, that is, 
$\| dK - \hat dK \|_{\mathrm{TV}} \le \| d - \hat d \|_{\mathrm{TV}}$,
giving the first summand in~\eqref{eq:kernel-perturbation}.
For the second term,
$\| \hat d(\mathcal P_\rho-\mathcal P_{\hat\rho}) \|_{\mathrm{TV}}
\le \sup_{x\in\mathcal X}\|\mathcal P_\rho(\cdot|x)-\mathcal P_{\hat\rho}(\cdot|x)\|_{\mathrm{TV}}
=\|\mathcal P_\rho-\mathcal P_{\hat\rho}\|_{\infty}$,
which is~\eqref{eq:kernel-perturbation}.
Note that
\begin{align*}
\|\mathcal P_\rho(\cdot|x)-\mathcal P_{\hat\rho}(\cdot|x)\|_{\mathrm{TV}}
&=\tfrac{1}{2}\sum_{s'\in \mathcal S,a'\in \mathcal A}\mathcal P(s'|x)\big|\pi_\rho(a'|s')-\pi_{\hat\rho}(a'|s')\big|\\
&\le \tfrac{1}{2}\sup_{s'\in \mathcal S}\sum_{a'\in\mathcal A}|\pi_\rho(a'|s')-\pi_{\hat\rho}(a'|s')|.
\end{align*}
Applying Lemma~\ref{lemma:B7} gives the bound on the policy difference by $C_{\mathcal Q}$ times the $\rho$-difference per state.
\hfill$\blacksquare$

The following lemma provides the bound on the term $\mathbb E[\langle V_k-\lambda(\rho_k),\Delta^{X}_{k+1}\rangle|\mathcal G_\delta]$ that appears in step~\hyperref[proof-Prop4_3-stepIII]{(III.)} of the proof of Proposition~\ref{prop: primal-tracking}. Recall that $\Delta^{X}_{k+1}$, defined in~\eqref{eq:index-law-mismatch}, denotes the mismatch between the law driving the update index selection and the ideal, mixed law.

\begin{lemma}
\label{lemma: Markov-difference}
    Let $K(\delta),\,\mathcal G_\delta$ be as in Prop.~\ref{prop: collected-facts-A2} and $\mathcal{K}$ as in~\eqref{eq: K_ZDR}. Moreover, let the Lipschitz envelope $L$ and the gradient bounds $B,\,D$ be as in Prop.~\ref{prop: primal-tracking}.
    For any $k\ge\max\{K(\delta),\mathcal K\}$, it holds that 
    \[
    \mathbb E[\langle V_k-\lambda(\rho_k),\Delta^{X}_{k+1}\rangle|\mathcal G_\delta]\le C_1\tau_k^2\alpha_{\mathrm{env,k-\tau_k}}\mathbb E[\|V_{k-\tau_k}-\lambda(\rho_{k-\tau_k})\|_2^2+1|\mathcal G_\delta],
    \]
    where $C_1=2(LDB+LB+LD+D^2+\tfrac{\gamma |\mathcal X|C^U}{\beta_0})$.
\end{lemma}
\textbf{Proof}\hspace{2pt}
The lemma is an adjusted version of~\cite[Lemma 2]{ref:zeng2024:two}, with the difference lying in the boundedness of the gradient in our projected setting. The proof is also closely related to theirs, with changes occurring at the points where a uniform bound on the gradient can be applied. 

The idea for handling the Markovian mismatch is to compare the inhomogeneous index process $\{X_t\}$ to an auxiliary chain frozen at parameter $\rho_{k-\tau_k}$; see also the proof of~\cite[Lemma 1]{ref:zeng2024:two}. Specifically, define
\[
\{\widetilde X_{k-\tau_k+t}\}_{t=0,\ldots,\tau_k}
\quad\text{by}\quad
\widetilde X_{k-\tau_k}=X_{k-\tau_k},\qquad
\widetilde X_{k-\tau_k+t+1}\sim \mathcal P_{\rho_{k-\tau_k}}(\,\cdot\,|\,\widetilde X_{k-\tau_k+t}).
\]
Using $z_k:=V_{k}-\lambda(\rho_{k})$, we define the following shorthand notation 
\begin{align*}
    T_1&=\mathbb E[\langle V_k-V_{k-\tau_k},\mathbb E\big[\widetilde G_{k+1}| \mathcal F_k^-\big]-\mathbb E_{X\sim \mu_{\rho_k}}[\tilde g_{k+1}(\rho_k,V_k;X)]\rangle| \mathcal G_\delta]\\
    T_2&=\mathbb E[\langle \lambda(\rho_{k-\tau_k})-\lambda(\rho_k) ,\mathbb E\big[\widetilde G_{k+1}| \mathcal F_k^-\big]-\mathbb E_{X\sim \mu_{\rho_k}}[\tilde g_{k+1}(\rho_k,V_k;X)]\rangle| \mathcal G_\delta]\\
    T_3&=\mathbb E[\langle z_{k-\tau_k}, \mathbb E[\widetilde G_{k+1}|\mathcal F_k^-]-\mathbb E[\tilde g_{k+1}(\rho_{k-\tau_k},V_{k-\tau_k};X_k)|\mathcal F_k^-]\rangle| \mathcal G_\delta]\\ 
    T_4&=\mathbb E[\langle z_{k-\tau_k}, \mathbb E[\tilde g_{k+1}(\rho_{k-\tau_k},V_{k-\tau_k};X_k)|\mathcal F_{k}^-]-\mathbb E[\tilde g_{k+1-\tau_k}(\rho_{k-\tau_k},V_{k-\tau_k};X_k)|\mathcal F_{k-\tau_k}^-]\rangle| \mathcal G_\delta]\\ 
    T_5&=\\
    &\mathbb E[\langle z_{k-\tau_k}, \mathbb E[\tilde g_{k+1-\tau_k}(\rho_{k-\tau_k},V_{k-\tau_k};X_k)|\mathcal F_{k-\tau_k}^-]-\mathbb E[\tilde g_{k+1-\tau_k}(\rho_{k-\tau_k},V_{k-\tau_k};\widetilde X_k)|\mathcal F_{k-\tau_k}^-]\rangle| \mathcal G_\delta]\\
    T_6&=\\
    &\mathbb E[\langle z_{k-\tau_k}, \mathbb E[\tilde g_{k+1-\tau_k}(\rho_{k-\tau_k},V_{k-\tau_k};\widetilde X_k)|\mathcal F_{k-\tau_k}^-]-\mathbb E_{\bar X\sim\mu_{\rho_{k-\tau_k}}}[\tilde g_{k+1-\tau_k}(\rho_{k-\tau_k},V_{k-\tau_k};\bar X)]\rangle| \mathcal G_\delta]\\
    T_7&=\\
    &\mathbb E[\langle z_{k-\tau_k}, \mathbb E_{\bar X\sim\mu_{\rho_{k-\tau_k}}}[\tilde g_{k+1-\tau_k}(\rho_{k-\tau_k},V_{k-\tau_k};\bar X)]-\mathbb E_{X\sim\mu_{\rho_{k}}}[\tilde g_{k+1-\tau_k}(\rho_{k-\tau_k},V_{k-\tau_k};X)]\rangle| \mathcal G_\delta]\\
    T_8&=\mathbb E[\langle z_{k-\tau_k}, \mathbb E_{X\sim\mu_{\rho_{k}}}[\tilde g_{k+1-\tau_k}(\rho_{k-\tau_k},V_{k-\tau_k};X)]-\mathbb E_{X\sim\mu_{\rho_{k}}}[\tilde g_{k+1}(\rho_{k},V_{k};X)]\rangle| \mathcal G_\delta]
\end{align*}

With the above and the application of the tower property, we get the following decomposition
\begin{align*}
    \mathbb E[\langle z_k,\Delta^{X}_{k+1}\rangle|\mathcal G_\delta]&=\mathbb E[\langle z_k,\mathbb E\big[\widetilde G_{k+1}| \mathcal F_k^-\big]
-\mathbb E_{X\sim \mu_{\rho_k}}[\tilde g_{k+1}(\rho_k,V_k;X)]\rangle|\mathcal G_\delta]\\
        &=\mathbb E[\langle V_k-V_{k-\tau_k},\mathbb E\big[\widetilde G_{k+1}| \mathcal F_k^-\big]
-\mathbb E_{X\sim \mu_{\rho_k}}[\tilde g_{k+1}(\rho_k,V_k;X)]\rangle|\mathcal G_\delta]\\
        &~~~~~+\mathbb E[\langle \lambda(\rho_{k-\tau_k})-\lambda(\rho_k),\mathbb E\big[\widetilde G_{k+1}| \mathcal F_k^-\big]
-\mathbb E_{X\sim \mu_{\rho_k}}[\tilde g_{k+1}(\rho_k,V_k;X)]\rangle|\mathcal G_\delta]\\
        &~~~~~+\mathbb E[\langle z_{k-\tau_k},\mathbb E\big[\widetilde G_{k+1}| \mathcal F_k^-\big]
-\mathbb E_{X\sim \mu_{\rho_k}}[\tilde g_{k+1}(\rho_k,V_k;X)]\rangle|\mathcal G_\delta]\\
&=T_1+T_2+T_3+T_4+T_5+T_6+T_7+T_8
\end{align*}

We proceed by bounding the above terms $T_1,\dots,T_8$. First, using the boundedness of the stochastic gradients $\|\widehat G_k\|_2\le D$ and $\bar \alpha_k\le\alpha_{\mathrm{env},k}$ on $\mathcal G_\delta$, we have 
\begin{align*}
    T_1&\le \mathbb E[\| V_k-V_{k-\tau_k}\|_2\|\mathbb E\big[\widetilde G_{k+1}| \mathcal F_k^-\big]
-\mathbb E_{X\sim \mu_{\rho_k}}[\tilde g_{k+1}(\rho_k,V_k;X)]\|_2| \mathcal G_\delta]\\
        &\le \mathbb E[\underbrace{\| V_k-V_{k-\tau_k}\|_2}_{\le\sum_{t=k+1-\tau_k}^{k}\bar\alpha_t\|\widehat G_t\|_2}\,\big(\underbrace{\|\mathbb E\big[\widetilde G_{k+1}| \mathcal F_k^-\big]\|_2+\|\mathbb E_{X\sim \mu_{\rho_k}}[\tilde g_{k+1}(\rho_k,V_k;X)]\|_2}_{\le 2D}\big)| \mathcal G_\delta]\\
        &\le 2D^2\tau_k\alpha_{\mathrm{env},k-\tau_k}.
\end{align*}
Similarly, using the boundedness of the dual update $\|\rho_{k+1}-\rho_k\|_2\le\bar\beta_kB$ and the Lipschitz continuity of the best response,
\begin{align*}
    T_2&\le\mathbb E[\| \underbrace{\lambda(\rho_{k-\tau_k})-\lambda(\rho_k)\|_2}_{\le L\|\rho_{k-\tau_k}-\rho_k\|_2}\,\underbrace{\|\mathbb E\big[\widetilde G_{k+1}| \mathcal F_k^-\big]-\mathbb E_{X\sim \mu_{\rho_k}}[\tilde g_{k+1}(\rho_k,V_k;X)]\|_2}_{\le 2 D}| \mathcal G_\delta]\\
    &\le 2DL\mathbb E[\sum_{t=k+1-\tau_k}^{k}\bar\beta_t\sup_{V,\rho}\|\nabla_\rho L(V,\rho)\|_2|\mathcal G_\delta]
    \le 2LDB\tau_k\beta_{\mathrm{env},k-\tau_k}.
\end{align*}
Next, use the Lipschitz continuity of $\nabla_VL$ with joint Lipschitz constant $L$ to get
\begin{align*}
    T_3&\le\mathbb E[\| z_{k-\tau_k}\|_2\,\| \mathbb E[\widetilde G_{k+1}|\mathcal F_k^-]-\mathbb E[\tilde g_{k+1}(\rho_{k-\tau_k},V_{k-\tau_k};X_k)|\mathcal F_k^-]\|_2| \mathcal G_\delta]\\ 
    &\le L\mathbb E[\| z_{k-\tau_k}\|_2\,\big(\|V_k-V_{k-\tau_k}\|_2+\|\rho_k-\rho_{k-\tau_k}\|_2\big)| \mathcal G_\delta]\\ 
    &\le L\tau_k(D\alpha_{\mathrm{env},k-\tau_k}+B\beta_{\mathrm{env},k-\tau_k})\mathbb E[\|z_{k-\tau_k}\|_2|\mathcal G_\delta]\\
    &\le L(D+B)\tau_k\alpha_{\mathrm{env},k}\,\mathbb E[\|z_{k-\tau_k}\|_2^2+1|\mathcal G_\delta],
\end{align*}
where the last inequality uses $\beta_k\le\alpha_k$. 

The term $T_4$ addresses the drift in the buffer over the mixing interval. Using the linearity of the gradient $\nabla_VL$ in the transition probabilities (see~\ref{eq: Vgrad}), we have 
 \begin{align*}
    T_4=\mathbb E[\langle z_{k-\tau_k}&, \mathbb E[\tilde g_{k+1}(\rho_{k-\tau_k},V_{k-\tau_k};X_k)|\mathcal F_{k}^-]-\mathbb E[\tilde g_{k+1-\tau_k}(\rho_{k-\tau_k},V_{k-\tau_k};X_k)|\mathcal F_{k-\tau_k}^-]\rangle| \mathcal G_\delta]\\ 
    &\le\gamma\mathbb E[\|z_{k-\tau_k}\|_2\|\sum_{x\in\mathcal X}\rho_{k-\tau_k}(x)\big(\mathcal P_{\mathcal D_k}(\cdot|x)-\mathcal P_{\mathcal D_{k-\tau_k}}(\cdot|x)\big)\|_2 |\mathcal G_\delta]\\
    &\le2\gamma\mathbb E[\|z_{k-\tau_k}\|_2\|\rho_{k-\tau_k}\|_1\|\mathcal P_{\mathcal D_k}-\mathcal P_{\mathcal D_{k-\tau_k}}\|_\infty |\mathcal G_\delta]\\
    &\le2\gamma |\mathcal X|C^U\mathbb E[\|z_{k-\tau_k}\|_2\|\mathcal P_{\mathcal D_k}-\mathcal P_{\mathcal D_{k-\tau_k}}\|_\infty |\mathcal G_\delta],
\end{align*}
where for the last inequality we recall the bound $\|\rho\|_1\le |\mathcal X|C^U$ for all $\rho\in H$ from Lemma~\ref{lemma:buffer-bias-HP}.
Fix a state-action pair $x=(s,a)\in\mathcal X$, and set $m:=\nu_{k-\tau_k}(x)$. Let $v_x\le\tau_k$ be the number of new samples for $x$ in the window, and let $u_s$ be the number of transitions from $x$ to $s$ in the window. Let $c_s\in\mathbb N$ be such that $\mathcal P_{\text{old}}(s|x):=\frac{c_s}{m}$ and
$\mathcal P_{\text{new}}(s|x):=\frac{c_s+u_s}{m+v_x}$ are the empirical next-state distributions, then
\[
\|\mathcal P_{\text{new}}(\cdot|x)-\mathcal P_{\text{old}}(\cdot|x)\|_{\mathrm{TV}}
=\tfrac12\sum_{j\in\mathcal S}\Big|\tfrac{c_j+u_j}{m+v_x}-\tfrac{c_j}{m}\Big|
\ \le\ \frac{v_x}{m}.
\]
Therefore,
\[
\|\mathcal P_{\mathcal D_k}-\mathcal P_{\mathcal D_{k-\tau_k}}\|_{\infty}
\ \le\ \max_{x\in\mathcal X}\frac{v_x}{\nu_{k-\tau_k}(x)}
\ \le\ \frac{\tau_k}{\min_{x\in\mathcal X}\nu_{k-\tau_k}(x)}.
\]
On $\mathcal G_\delta$, it holds that $\nu_t(x)\ge\tfrac{p_\star}{2}t$ for $t\ge K(\delta)$. 
By the definition of $\beta_{\mathrm{env},k-\tau_k}=\tfrac{\beta_0}{\lfloor p_\star/2\,(k-\tau_k)\rfloor+1}$, we have $$\tfrac{1}{p_\star/2(k-\tau_k)}\le(1+\tfrac{2}{p_\star \mathcal K})\tfrac{1}{p_\star/2(k-\tau_k)+1}\le\tfrac{2}{\beta_0}\beta_{\mathrm{env},k-\tau_k}.$$ 
Consequently, 
\begin{align*}
        T_4&\le \tfrac{2\gamma |\mathcal X|C^U}{\beta_0}\tau_k\beta_{\mathrm{env},k-\tau_k}\mathbb E[\|z_{k-\tau_k}\|_2^2+1|\mathcal G_\delta].
\end{align*}
The next term, $T_5$, compares the inhomogeneous chain $\{X_k\}$ over the mixing window of length $\tau_k$ with the frozen chain $\{\widetilde X_k\}$. 
Let us introduce the following notation for the conditional probability laws 
$$p_t(\cdot)=\mathbb P(X_t=\cdot|\mathcal F_{t-1}^-)\qquad \text{and} \qquad \tilde p_t(\cdot)=\mathbb P(\widetilde X_t=\cdot|\mathcal F_{t-1}^-).$$ 
By construction of the frozen chain, $\widetilde X_{k-\tau_k}=X_{k-\tau_k}$, hence $\|p_{k-\tau_k}-\tilde p_{k-\tau_k}\|_{\mathrm{TV}}=0$.
Using the tower property for $\mathcal F_{k-\tau_k}^-\subset\mathcal F_{k-1}^-$ and the definition of the total variation distance, $\|p_t-\tilde p_t\|_\mathrm{TV}=\tfrac{1}{2}\sup_{h:\mathcal X\rightarrow [-1,1]}\Bigl|\int_\mathcal X h(x)(p_t-\tilde p_t)(dx)\Bigl|$, 
\begin{align*}
        \|\mathbb E[&\tilde g_{k+1-\tau_k}(\rho_{k-\tau_k},V_{k-\tau_k};X_k)|\mathcal F_{k-\tau_k}^-]-\mathbb E[\tilde g_{k+1-\tau_k}(\rho_{k-\tau_k},V_{k-\tau_k};\widetilde X_k)|\mathcal F_{k-\tau_k}^-]\|_2\\
        &= \|\mathbb E\big[\mathbb E[\tilde g_{k+1-\tau_k}(\rho_{k-\tau_k},V_{k-\tau_k};X_k)-\tilde g_{k+1-\tau_k}(\rho_{k-\tau_k},V_{k-\tau_k};\widetilde X_k)|\mathcal F_{k-1}^-]|\mathcal F_{k-\tau_k}^-\big]\|_2\\
        &=\Bigl\|\mathbb E\left[\int_{\mathcal X}\tilde g_{k+1-\tau_k}(\rho_{k-\tau_k},V_{k-\tau_k};x)(p_k-\tilde p_k(dx)|\mathcal F_{k-\tau_k}^-\right]\Bigl\|_2\\
        &\le 2D\mathbb E[\|p_k-\tilde p_k\|_\mathrm{TV}| \mathcal F_{k-\tau_k}^-],
\end{align*}
where we used Jensen's inequality and the bound $\|\int g\,d(p-\tilde p)\|_2 \le 2\sup_x\|g(x)\|_2\,\|p-\tilde p\|_{\mathrm{TV}}$ together with $\sup_x\|\tilde g_{k+1-\tau_k}(\cdot;\,x)\|_2\le D$. 
The one-step laws satisfy $p_t=p_{t-1}\mathcal P_{\rho_{t-1}},\qquad \tilde p_t=\tilde p_{t-1}\mathcal P_{\rho_{k-\tau_k}}$, therefore, with~\eqref{eq:kernel-perturbation} and the Lipschitz bound on the kernels it holds that $$\|p_t-\tilde p_t\|_\mathrm{TV}\le \|p_{t-1}-\tilde p_{t-1}\|_\mathrm{TV}+\|\mathcal P_{\rho_{t-1}}-\mathcal P_{\rho_{k-\tau_k}}\|_\infty\le\|p_{t-1}-\tilde p_{t-1}\|_\mathrm{TV}+L\|\rho_{t-1}-\rho_{k-\tau_k}\|_2.$$
Iterating from $t=k$ to $k-\tau_k+1$ and using $\|p_{k-\tau_k}-\tilde p_{k-\tau_k}\|_{\mathrm{TV}}=0$ gives 
$$\|p_k-\tilde p_k\|_\mathrm{TV}\le L\sum_{t=k-\tau_k}^{k-1}\|\rho_t-\rho_{k-\tau_k}\|_2.$$

Moreover, on $\mathcal G_\delta$ the dual increments for $t\ge k-\tau_k$ satisfy
$\|\rho_{t+1}-\rho_t\|_2\le B\bar\beta_{t+1}\le B\beta_{\mathrm{env},k-\tau_k}$, hence
$\|\rho_t-\rho_{k-\tau_k}\|_2 \le B(t-(k-\tau_k))\beta_{\mathrm{env},k-\tau_k}$ and thus
$\sum_{t=k-\tau_k}^{k-1}\|\rho_t-\rho_{k-\tau_k}\|_2 \le B\tau_k^2\beta_{\mathrm{env},k-\tau_k}$.
Combining the bounds yields
$$2D\,\mathbb E[\|p_k-\tilde p_k\|_{\mathrm{TV}}\mid \mathcal F_{k-\tau_k}^-]
\le 2LDB\,\tau_k^2\beta_{\mathrm{env},k-\tau_k}.$$
Consequently,
$$
T_5\le 2LDB\,\tau_k^2\beta_{\mathrm{env},k-\tau_k}\,\mathbb E[\|z_{k-\tau_k}\|_2\mid\mathcal G_\delta]
\le LDB\,\tau_k^2\beta_{\mathrm{env},k-\tau_k}\,\mathbb E[\|z_{k-\tau_k}\|_2^2+1\mid\mathcal G_\delta].
$$

\noindent
The term $T_6$ compares expectation under the fixed law chain to its stationary law expectation, and by definition of $\tau_k=\tau(\beta_{\mathrm{env},k})$ we have
\begin{align*}
    T_6&\le 2D\mathbb E[\|z_{k-\tau_k}\|_2\mathbb E[\underbrace{\,\sup_{x\in \mathcal X}\|\delta_x(\mathcal P_{\rho_{k-\tau_k}})^{\tau_k}-\mu_{\rho_{k-\tau_k}}\|_{\mathrm{TV}}}_{\le \beta_{\mathrm{env},k}}|\mathcal F_{k-\tau_k}^-]|\mathcal G_\delta]\\
    &\le D\beta_{\mathrm{env},k-\tau_k}\mathbb E[\|z_{k-\tau_k}\|_2^2+1|\mathcal G_\delta].
\end{align*}
To bound the stationary law drift comparison term $T_7$, we use the Lipschitz continuity of the stationary law, derived in Lemma~\ref{lemma:Lip-stationary}, and obtain
\begin{align*}
    T_7&\le\mathbb E[\|z_{k-\tau_k}\|_2\,\| \mathbb E_{\bar X\sim\mu_{\rho_{k-\tau_k}}}[\tilde g_{k+1-\tau_k}(\rho_{k-\tau_k},V_{k-\tau_k};\bar X)]\\
    &~~~~~~~~~~~~~~~~~~~~~~~~~~~~~~~~~-\mathbb E_{X\sim\mu_{\rho_{k}}}[\tilde g_{k+1-\tau_k}(\rho_{k-\tau_k},V_{k-\tau_k};X)]\|_2| \mathcal G_\delta]\\
    &\le 2D\mathbb E[\|z_{k-\tau_k}\|_2\,\|\mu_{\rho_{k-\tau_k}}-\mu_{\rho_{k}}\|_{\mathrm{TV}}|\mathcal G_\delta]\\
    &\le LDB\tau_k\beta_{\mathrm{env},k-\tau_k}\mathbb E[\|z_{k-\tau_k}\|_2^2+1|\mathcal G_\delta],
\end{align*}
where the last inequality combines the Lipschitz bound with the bounded dual update over the mixing window, similar to the derivation in the bound to $T_2$.

Finally, to bound $T_8$, we have by the Lipschitz continuity of the primal gradient
\begin{align*}
    T_8&\le \mathbb E[\| z_{k-\tau_k}\|_2\, \|\mathbb E_{X\sim\mu_{\rho_{k}}}[\tilde g_{k+1-\tau_k}(\rho_{k-\tau_k},V_{k-\tau_k};X)-\tilde g_{k+1}(\rho_{k},V_{k};X)]\|_2| \mathcal G_\delta]\\
    &\le L\mathbb E[\| z_{k-\tau_k}\|_2\, \big(\|V_k-V_{k-\tau_k}\|_2+\|\rho_k-\rho_{k-\tau_k}\|_2\big)| \mathcal G_\delta]\\
    &\le L(D+B)\tau_k\alpha_{\mathrm{env},k-\tau_k}\mathbb E[\| z_{k-\tau_k}\|_2^2+1|\mathcal G_\delta],
\end{align*}
where the final inequality follows by the same arguments as used for the bound in $T_3$. Putting it all together, we find with $\beta_k\le\alpha_k$ that 
\begin{align*}
        \mathbb E[\langle z_k,\Delta^{X}_{k+1}\rangle|\mathcal G_\delta]&\le\\
        \big(2L(D+B)&+LDB\tau_k+D+LDB+\tfrac{2\gamma |\mathcal X|C^U}{\beta_0}\big)\tau_k\alpha_{\mathrm{env},k}\,\mathbb E[\|z_{k-\tau_k}\|_2^2+1|\mathcal G_\delta]\\
        &~~~~~~~~~+2D(LB+D)\tau_k\alpha_{\mathrm{env},k-\tau_k}\\
        &\le \underbrace{2(LDB+LB+LD+D^2+\tfrac{\gamma |\mathcal X|C^U}{\beta_0})}_{C_1}\tau_k^2\alpha_{\mathrm{env},k-\tau_k}\mathbb E[\|z_{k-\tau_k}\|_2^2+1|\mathcal G_\delta],
\end{align*}
which finishes the proof.
\hfill $\blacksquare$

The bound on the Markov mismatch term in Proposition~\ref{prop: dual-variable-error} is given by~\cite[Lemma~1]{ref:zeng2024:two} which is restated here for completeness, without proof. We adapted the statement to our notation using \(\bar\beta_k\le\beta_{\mathrm{env},k}\) on \(\mathcal G_\delta\).
\begin{lemma}
\label{lemma:dual-markov-mismatch}
Fix \(\delta\in(0,1)\) and let \(k\ge \max\{K(\delta),\mathcal K\}\). On the event
\(\mathcal G_\delta\)
\[
-\mathbb E\!\left[
\left\langle \nabla f(\rho_k),\,\Delta^X_{k+1}\right\rangle
\;\middle|\;
\mathcal G_\delta
\right]
\le
12L^2B^3\,\tau_k^2\,\beta_{\mathrm{env},k-\tau_k}.
\]
\end{lemma}

\section{Additional Proofs}
\label{app:additional:proofs}
\paragraph{Proof of Corollary~\ref{cor: dual-reg}} 
The proof for the entropy-regularized MDP is provided in \cite[Section 3]{ying2020note}. We adapt the proof to apply to a general strongly convex regularizer. Although the proof closely follows \cite{ying2020note}, we include it to ensure the paper remains self-contained.
Throughout this proof we adopt the shorthand notation $v_s$ or $v_{sa}$ for $v(s)$ and $v(s,a)$, respectively for vectors in $\mathbb{R}^{|\mathcal S|}$ and $\mathbb{R}^{|\mathcal S||\mathcal A|}$ or corresponding mappings from $\mathcal S$ or $\mathcal S\times\mathcal A$ to $\mathbb R$.

We show the correspondence of the solution of the regularized Bellman optimality operator and the optimal primal variable in (\ref{eq: dual-reg LP}). We begin by formulating the corresponding primal linear program and then verify that its solution is the unique fixed point of the regularized Bellman optimality operator (\ref{eq:regBellman:opt}). Decompose $\rho_{sa}=\tilde{\rho}_s\pi_{sa}$ with $\tilde{\rho}_s=\sum_{a\in\mathcal{A}}\rho_{sa}$ and the policy $\pi=\pi_\rho\in\Pi$. For a homogeneous regularizer $G$ we have $G(\rho_s)=G(\tilde{\rho}_s\pi_s)=\tilde{\rho}_sG(\pi_s)$ for all $s\in\mathcal{S}$. Due to the convex regularizing term, we can replace the condition $\rho\geq0$ of the regularized primal-dual problem with the strict inequality $\rho>0$.

When considering the problem~\eqref{eq: dual-reg LP}, we recall the fact that the primal variable $V$ can without loss of generality be restricted to a compact set $\mathcal{V}_r\subset \mathbb{R}_+^{|\mathcal{S}|}$, where 
\begin{align}
\label{eq: V-reg-Set}
    \mathcal{V}_r = \{v\in\mathbb{R}_+^{|\mathcal{S}|} : v_i \leq \frac{C_r+\eta_\rho U_G}{1-\gamma}, \ \forall i=1,\dots |\mathcal{S}|\}.
\end{align}
Inserting the decomposition above into (\ref{eq: dual-reg LP}) and using the convex-concave structure of the regularized Lagrangian, which allows us to use Sion's minimax theorem~\citep{Sion:1958}, we obtain

\begin{align*}
    &\min_{V\in \mathcal{V}_r}\max_{\pi\in\Pi,\tilde{\rho}_s>0}\mu^\top V+\sum_{s\in\mathcal{S}}\tilde{\rho}_s\left(\sum_{a\in\mathcal{A}}\pi_{sa}\Delta[V]_{sa}+\eta_\rho G(\pi_s)\right)\\
    &\Leftrightarrow \min_{V\in \mathcal{V}_r}\left(\mu^\top V+\max_{\pi\in\Pi,\tilde{\rho}_s>0}\sum_{s\in\mathcal{S}}\tilde{\rho}_s\left(\sum_{a\in\mathcal{A}}\pi_{sa}(r_{sa}+\gamma\sum_{s'\in\mathcal{S}} V_{s^{\prime}} \mathcal{P}(s'|s,a))+\eta_\rho G(\pi_s)-V_s\right)\right)\\
    &\Leftrightarrow \min_{V\in \mathcal{V}_r}\left(\mu^\top V+\sup_{\tilde{\rho}_s>0}\sum_{s\in\mathcal{S}}\tilde{\rho}_s\max_{\pi_s\in\Delta_\mathcal{A}}\left(\sum_{a\in\mathcal{A}}\pi_{sa}(r_{sa}+\gamma\sum_{s'\in\mathcal{S}} V_{s^{\prime}} \mathcal{P}(s'|s,a))+\eta_\rho G(\pi_s)-V_s\right)\right).
\end{align*}

Then the primal linear programming formulation of (\ref{eq: dual-reg LP}) is 
\begin{align*}
    P_r: \begin{cases}
        \min\limits_{V\in\mathbb{R}^{|\mathcal{S}|}} \mu^\top V\\
       \ \ \text{s.t.} \ \max\limits_{\pi_s\in\Delta_\mathcal{A}}\sum_{a\in\mathcal{A}}\pi_{sa}(r_{sa}+\gamma\sum_{s'\in\mathcal{S}} V_{s^{\prime}} \mathcal{P}(s'|s,a))+\eta_\rho G(\pi_s)\leq V_s \quad \forall s\in\mathcal{S}.
    \end{cases}
\end{align*}
We note that the constraint $V\in \mathcal{V}_r$ can be relaxed to $V\in\mathbb{R}^{|\mathcal{S}|}$ in $P_r$ without loss of generality. Indeed, if we suppose that there is a feasible candidate $V$ of $P_r$ such that $V\notin \mathcal{V}_r$, we can construct a feasible candidate with a smaller objective value by projecting it onto $\mathcal{V}_r$.

The constraint of the primal LP is equivalently expressed as $T_{\star,G}V\leq V$ with the regularized Bellman optimality operator defined in (\ref{eq:regBellman:opt}). Let $\widehat V$ be a solution to $P_r$. Then by primal feasibility, the slack variable $S:=\widehat V-T_{\star,G}\widehat V\geq0$. Assume $S\neq0$. Note that $T_{\star,G}$ is a monotone operator, \cite[Prop. 2]{geist2019theory}, that is, for $V_1\geq V_2$, $T_{\star,G}V_1\geq T_{\star,G}V_2$. Therefore $\widetilde V:=T_{\star,G}\widehat V$ fulfills $T_{\star,G}\widetilde V\leq T_{\star,G}\widehat V=\widetilde V$, that is, $\widetilde V$ is primal feasible. For $S\neq0$, there exists a state $s$ with $\widetilde V_s<\widehat V_s$ and hence $\mu^\top \widetilde V<\mu^\top \widehat V$. This contradicts the optimality of $\widehat V$, therefore $S=0$ holds. The solution to the primal problem is the unique fixed point of the regularized Bellman optimality operator.  

The dual problem for the regularized MDP is 
\begin{align*}
    D_r: \begin{cases}
        \sup\limits_{\rho>0} \ \ \sum_{(s,a)\in\mathcal{S}\times\mathcal{A}}r_{sa}\rho_{sa}+\eta_\rho\sum_{s\in\mathcal{S}}\tilde{\rho}_sG(\pi_s)\\
        \text{s.t.} \ \tilde\rho_{s'}-\gamma\sum_{(s,a)\in\in\mathcal{S}\times\mathcal{A}}\rho_{sa}\mathcal P(s'|s,a)=\mu_{s'} \quad \forall s'\in\mathcal{S}.
    \end{cases}
\end{align*}
Inserting the decomposition $\rho_{sa}=\tilde{\rho}_s\pi_{s a}$ the dual constraint is $$\tilde{\rho}_{s'}-\gamma \sum_{(s,a)\in\mathcal{S}\times\mathcal{A}} \mathcal P(s'|s,a) \pi_{sa} \tilde{\rho}_{s}=\mu_{s'}.$$ In vector notation with $\sum_{a\in\mathcal{A}} \mathcal P(s'|s,a) \pi_{sa}=\mathcal P_{s s'}^{\pi}=\left(\mathcal P^\pi\right)_{s' s}^\top $ the constraint can be rearranged to 
$$\tilde{\rho}=\left(I-\gamma\left(\mathcal P^\pi\right)^{\top}\right)^{-1} \mu.$$ 
Denote this as $\tilde{\rho}^\pi$, highlighting the dependence on $\pi$. By inserting the optimal $\tilde{\rho}^\pi$ into the objective, we can reformulate the regularized dual problem as
\begin{align*}
&\max_{\pi\in\Pi} \sum_{s\in\mathcal{S}} \tilde{\rho}_s^\pi \underbrace{\sum_{a\in\mathcal{A}} r_{sa} \pi_{sa}}_{=r_s^\pi}+\eta_\rho \sum_{s\in\mathcal{S}} \tilde{\rho}^\pi_s G\left(\pi_s\right)\\
&=\max_{\pi\in\Pi}  \sum_{s\in\mathcal{S}} \tilde{\rho}_s^\pi\left(r_s^\pi+\eta_\rho G\left(\pi_s\right)\right)\\
&=\max_{\pi\in\Pi}  (\tilde{\rho}^\pi)^\top \left(r^\pi+\eta_\rho G\left(\pi\right)\right)\\
&=\max_{\pi\in\Pi} \mu^\top \left(I-\gamma \mathcal P^\pi\right)^{-1}\left(r^\pi+\eta_\rho G\left(\pi\right)\right)
\end{align*}
The regularized Bellman operator states $V^\pi=r^\pi+\eta_\rho G(\pi) +\gamma \mathcal P^\pi V^\pi$ or equivalently 
$V^\pi=\left(I-\gamma \mathcal P^\pi\right)^{-1}\left(r^\pi+\eta_\rho G(\pi)\right)$. Hence, the regularized dual problem $D_r$ is equivalent to
$$\max_{\pi\in\Pi}  \mu^\top  V\quad \text{s.t.}\quad V=T_{\pi, G} V=r^\pi+\eta_\rho G(\pi)+\gamma \mathcal P^\pi V.$$
Since $\mu$ is an arbitrary positive vector, this is the policy gradient problem with the unique solution $\pi^\star_r=\arg\max_{\pi\in\Pi}  T_{\pi,G}V^\pi_r$. By the decomposition of $\rho$, we have $\pi^\star_{sa}=\frac{\rho^\star_{sa}}{\tilde{\rho}^{\pi^\star}_s}$, that is, $\pi_{\rho^\star_r}=\pi^\star_r$. 
\hfill $\blacksquare$

\paragraph{Proof of Lemma~\ref{lemma: dualBound}} We begin by establishing an upper bound for $\rho^\star$. The optimal primal and dual solutions are linked by the first-order optimality condition of the primal variable $V$. We then use the bound on $V^\star$ to obtain a bound on $\rho^\star$. 
Moreover, thanks to Theorem~\ref{thm:Li-et-al}, we know that the tuple $(\pi^\star$, $V^\star)$ is a fixed point of the regularized Bellman operator~\eqref{eq:regBellman:opt} and therefore fulfills the regularized Bellman equation $$V^\star=r_{\pi^\star}+\gamma \mathcal P^{\pi^\star} V^\star+\eta_\rho G_{\pi^\star},$$ with $(r_{\pi})_s:=\sum_{a\in\mathcal{A}}r_{sa}\pi_{sa}$, $(\mathcal P^\pi)_{ss'}=\sum_{a\in\mathcal{A}}\pi_{sa}\mathcal{P}(s'| s,a)$ and $(G_\pi)_s=G(\pi(\cdot| s))$. Therefore, 
\begin{align*}
    \|V^\star\|_\infty&=\|r_{\pi^\star}+\gamma \mathcal P^{\pi^\star} V^\star+\eta_\rho G_{\pi^\star}\|_\infty\\
    &\leq \|r_{\pi^\star}\|_\infty+\gamma \underbrace{\|\mathcal P^{\pi^\star}\|_\infty}_{=1}\| V^\star\|_\infty+\eta_\rho \|G_{\pi^\star}\|_\infty\\
    &\leq C_r+\gamma\| V^\star\|_\infty+\eta_\rho U_G.\\
\end{align*}    
In the last inequality, we applied the upper bounds on the reward and regularizer. Rearranging terms we get
\begin{align}
\label{eq: primal-var-bound}
    \|V^\star\|_\infty&\leq\frac{C_r+\eta_\rho U_G}{1-\gamma}.
\end{align}
Based on the gradient (\ref{eq: Vgrad}), the first-order optimality conditions for~\eqref{eq: regLP} state 
\begin{equation}
\label{eq:CU-derivation}
\begin{aligned}
& \tilde\rho_{s^{\prime}}^\star-\gamma \sum_{(s,a)\in\mathcal{S}\times\mathcal{A}} \rho_{sa}^\star \mathcal{P}\left(s^{\prime} | s,a\right)=\eta_V V_{s^{\prime}}^\star \leqslant \frac{\eta_V (C_r+\eta_\rho U_G)}{1-\gamma} \\
& \Rightarrow \quad \sum_{s^{\prime}\in\mathcal{S}} \tilde\rho_{s^{\prime}}^\star \leqslant \frac{|\mathcal{S}| \eta_V (C_r+\eta_\rho U_G)}{1-\gamma}+\gamma \sum_{(s,a)\in\mathcal{S}\times\mathcal{A}} \rho_{sa}^\star \underbrace{\sum_{s^{\prime}\in\mathcal{S}} \mathcal{P}\left(s^{\prime} | s, a\right)}_{=1} \\
& \Leftrightarrow \quad \sum_{(s,a)\in\mathcal{S}\times\mathcal{A}} \rho_{s a}^\star \leqslant \frac{|\mathcal{S}| \eta_V (C_r+\eta_\rho U_G)}{(1-\gamma)^2}.
\end{aligned}
\end{equation}
Combined with non-negativity, we find that for $C^{\mathrm{max}}:=\frac{|\mathcal{S}| \eta_V (C_r+\eta_\rho U_G)}{(1-\gamma)^2}$ and $C^U:=2C^{\mathrm{max}}$ it holds that $\left\|\rho^\star\right\|_{\infty}\leq\|\rho^\star\|_1\leq C^{\mathrm{max}}< C^U<\infty$.
The lower bound for $\rho^\star$ follows from the normalized first-order optimality conditions for $\rho$ based on the gradient (\ref{eq: rgrad}), stating that the optimal policy takes the form of a Boltzmann policy
\begin{align*}
    \pi_{sa}^\star&=\frac{\rho^\star_{sa}}{\tilde{\rho}^\star_s}=\frac{\exp(\frac{1}{\eta_\rho}\Delta[V^\star]_{sa})}{\sum_{a\in\mathcal{A}}\exp(\frac{1}{\eta_\rho}\Delta[V^\star]_{sa})}\\
    &=\frac{\exp(\frac{1}{\eta_\rho}(-V_s+r_{sa}+\gamma\sum_{s'\in\mathcal{S}}\mathcal{P}(s'| s,a)V^\star_{s'})}{\sum_{a\in\mathcal{A}}\exp(\frac{1}{\eta_\rho}(-V_s+r_{sa}+\gamma\sum_{s'\in\mathcal{S}}\mathcal{P}(s'| s,a)V^\star_{s'})}.
\end{align*}
As we derived in~\eqref{eq: primal-var-bound}, the optimal regularized value function is bounded, and together with $r_{sa}\in[0,C_r]$ we get $$\Delta[V^\star]_{sa}=-V^\star_s+r_{sa}+\gamma\sum_{s'\in\mathcal{S}}\mathcal{P}(s'| s,a)V^\star_{s'}\geq-\|V^\star\|_\infty$$ and $$\Delta[V^\star]_{sa}\leq C_r+\gamma\|V^\star\|_\infty.$$ Hence, $$\min_{(s,a)\in\mathcal{S}\times\mathcal{A}}\pi^\star_{sa}\geq\frac{\exp(-\frac{1}{\eta_\rho}\|V^\star\|_\infty)}{|\mathcal{A}|\exp(\frac{1}{\eta_\rho}(C_r+\gamma\|V^\star\|_\infty))}=\frac{1}{|\mathcal{A}|}\exp(-\frac{1}{\eta_\rho}[C_r+(1+\gamma)\|V^\star\|_\infty]).$$ With $\|V^\star\|_\infty\leq\frac{C_r+\eta_\rho U_G}{1-\gamma}$ we obtain the lower bound for the optimal policy 
\begin{equation*}
\min_{(s,a)\in\mathcal{S}\times\mathcal{A}}\pi^\star_{sa}\geq\frac{1}{|\mathcal{A}|}\exp(-\frac{2\eta_\rho^{-1}C_r+(1+\gamma)U_G}{1-\gamma})=:C^L_1. 
\end{equation*}
Next, we establish a lower bound for $\tilde{\rho}$. Note that as a consequence of (\ref{eq: Vgrad}), $$\tilde{\rho}^\star=\eta_V(I-\gamma (\mathcal P^{\pi^\star})^\top )^{-1}V^\star$$ with $(P^\pi)_{ss'}:=\sum_{a\in\mathcal{A}}\pi_{sa}\mathcal{P}(s'|s,a)$. With the Neumann series, we can expand the inverse and obtain 
\begin{align*}
    \tilde{\rho}^\star&=\eta_V\left(\sum_{k=0}^\infty\gamma^k((\mathcal P^{\pi^\star})^\top )^k\right)V^\star=\eta_VV^\star+\eta_V\left(\sum_{k=1}^\infty\gamma^k((\mathcal P^{\pi^\star})^\top )^k\right)V^\star\geq\eta_VV^\star.
\end{align*}
Therefore, for all $s'\in\mathcal{S}$ it holds that
\begin{align*}
    \tilde{\rho}_{s'}^\star&\geq\eta_VV^\star_{s'}\\
    &=\eta_V\underbrace{(r_{\pi^\star}+\gamma P^{\pi^\star}V^\star)_{s'}}_{\geq 0}+\eta_V\eta_\rho \underbrace{(G_{\pi^\star})_{s'}}_{=-\sum_{a\in\mathcal{A}}\pi^\star_{s' a}\log\pi^\star_{s' a}}\\
    &\geq-\eta_V\eta_\rho\sum_{a\in\mathcal{A}}\pi^\star_{s' a}\log\pi^\star_{s' a}\\
    &\geq \min_{s\in\mathcal{S}}-\eta_V\eta_\rho\sum_{a\in\mathcal{A}}\pi^\star_{sa}\log\pi^\star_{sa}.
\end{align*}
The equality in the second line holds since $V^\star$ is the fixed point of the regularized Bellman operator (\ref{eq:regBellman}). The above derivation states that $\tilde{\rho}^\star_{s'}$ is bounded from below by the minimal entropy of the optimal policy, where the minimum is over the states $s\in\mathcal{S}$. Using the Boltzmann form of the optimal policy, we now turn to bound the minimal entropy.

Define $\Delta_{\max}(s):=\max_{a\in\mathcal{A}}\Delta[V^\star]_{sa}$ and $\Delta_{\min}(s):=\min_{a\in\mathcal{A}}\Delta[V^\star]_{sa}$. We then denote the maximum difference as $\Delta_s:=\Delta_{\max}(s)-\Delta_{\min}(s)$. Set $K_1:=|\mathcal{A}|-1$, $K_2(\Delta_s):=\exp(-\frac{1}{\eta_\rho}\Delta_s)$. The minimum entropy of the softmax policy is obtained when in state $s$ one action gets the ``high" weight, $\Delta_{\max}(s)$, and the other $K_1$ states get ``low" weight $\Delta_{\min}$.
In that case, the policy assigns the probabilities 
\begin{align*}
    p_{\max}&=\frac{\exp(\frac{1}{\eta_\rho}\Delta_{\max}(s))}{\exp(\frac{1}{\eta_\rho}\Delta_{\max}(s))+K_1\exp(\frac{1}{\eta_\rho}\Delta_{\min}(s))}=\frac{1}{1+K_1\exp(-\frac{1}{\eta_\rho}\Delta_s)}=\frac{1}{1+K_1K_2(\Delta_s)}\\
    p_{\min}&=\frac{K_2(\Delta_s)}{1+K_1K_2(\Delta_s)}.
\end{align*}
The corresponding entropy is
$$(G_\pi)_s=-(p_{\max} \log p_{\max}+K_1 p_{\min} \log p_{\min})=\log(1+K_1K_2(\Delta_s))+\frac{K_1}{\eta_\rho}\frac{\Delta_sK_2(\Delta_s)}{1+K_1K_2(\Delta_s)}.$$ This term is monotonically decreasing in $\Delta_s$. We know that $\Delta_s\leq C_r+(1+\gamma)\|V^\star\|_\infty\leq \frac{2C_r+(1+\gamma)\eta_\rho U_G}{1-\gamma}=:\bar\Delta$. Therefore, we can bound the minimum entropy of the optimal policy from below by inserting the upper bound for $\Delta_s$ and find $$\tilde{\rho}_{s'}^\star\geq \min_{s\in\mathcal{S}}-\eta_V\eta_\rho\sum_{a\in\mathcal{A}}\pi^\star_{sa}\log\pi^\star_{sa}\geq \eta_V\eta_\rho\log(1+K_1K_2(\bar\Delta))+\eta_V\frac{K_1\bar\Delta K_2(\bar\Delta)}{1+K_1K_2(\bar\Delta)}=:C^L_2.$$ Applying the decomposition $\rho^\star_{sa}=\pi^\star_{sa}\tilde{\rho}^\star_s$, we find that $\rho^\star_{sa}\geq C^L_1C^L_2=: C^{\mathrm{min}}$ and for $C^L:=\tfrac{1}{2}C^{\mathrm{min}}$ the Lemma statement now follows.

\hfill $\blacksquare$

\bibliography{refs}

@article{chen2016,
  title={Stochastic primal-dual methods and sample complexity of reinforcement learning},
  author={Chen, Yichen and Wang, Mengdi},
  journal={arXiv preprint arXiv:1612.02516},
  year={2016}
}

@InProceedings{ref:ozdaglar-23a,
  title = 	 {Revisiting the Linear-Programming Framework for Offline {RL} with General Function Approximation},
  author =       {Ozdaglar, Asuman E. and Pattathil, Sarath and Zhang, Jiawei and Zhang, Kaiqing},
  booktitle = 	 {International Conference on Machine Learning},
  pages = 	 {26769--26791},
  year = 	 {2023},
}

@article{Sion:1958,
  author  = {Sion, Maurice},
  title   = {On general minimax theorems},
  journal = {Pacific Journal of Mathematics},
  volume  = {8},
  number  = {1},
  pages   = {171--176},
  year    = {1958},
}

@article{kushner2003,
  title={Convergence with Probability One: Correlated Noise},
  author={Kushner, Harold J. and Yin, George},
  journal={Stochastic Approximation and Recursive Algorithms and Applications},
  pages={161--212},
  year={2003},
  publisher={Springer}
}

@article{dupuis1993dynamical,
  title={Dynamical systems and variational inequalities},
  author={Dupuis, Paul and Nagurney, Anna},
  journal={Annals of Operations Research},
  volume={44},
  pages={7--42},
  year={1993},
  publisher={Springer}
}

@article{dupuis1987large,
  title={Large deviations analysis of reflected diffusions and constrained stochastic approximation algorithms in convex sets},
  author={Dupuis, Paul},
  journal={Stochastics: An International Journal of Probability and Stochastic Processes},
  volume={21},
  number={1},
  pages={63--96},
  year={1987},
  publisher={Taylor \& Francis}
}

@book{borkar2023stochastic,
  title={Stochastic Approximation: A Dynamical Systems Viewpoint},
  author={Borkar, Vivek S.},
  edition={Second},
  year={2023},
  publisher={Springer Nature Singapore}
}

@book{puterman1994markov,
author = {Puterman, Martin L.},
title = {{M}arkov Decision Processes: Discrete Stochastic Dynamic Programming},
year = {1994},
isbn = {0471619779},
publisher = {John Wiley \& Sons, Inc.},
edition = {1st},
}

@article{borkar1997stochastic,
  title={Stochastic approximation with two time scales},
  author={Borkar, Vivek S.},
  journal={Systems \& Control Letters},
  volume={29},
  number={5},
  pages={291--294},
  year={1997},
  publisher={Elsevier}
}

@article{borkar2000ode,
  title={The ODE method for convergence of stochastic approximation and reinforcement learning},
  author={Borkar, Vivek S. and Meyn, Sean P.},
  journal={SIAM Journal on Control and Optimization},
  volume={38},
  number={2},
  pages={447--469},
  year={2000},
  publisher={SIAM}
}

@book{meyn2022control,
  title={Control Systems and Reinforcement Learning},
  author={Meyn, Sean P.},
  year={2022},
  publisher={Cambridge University Press}
}

@article{ref:Prashant-20,
 author = {Mehta, Prashant G. and Meyn, Sean P.},
 journal = {arXiv preprint},
 title = {Convex {Q}-Learning, Part 1: Deterministic Optimal Control},
 volume = {2008.03559},
 year = {2020}
}

@article{ref:lu2023convex,
  title={Convex {Q}-learning in a stochastic environment: Extended version},
  author={Lu, Fan and Meyn, Sean P.},
  journal={arXiv preprint arXiv:2309.05105},
  year={2023}
}

@article{li2024accelerating,
  title={Accelerating Primal-Dual Methods for Regularized {M}arkov Decision Processes},
  author={Li, Haoya and Yu, Hsiang-Fu and Ying, Lexing and Dhillon, Inderjit S.},
  journal={SIAM Journal on Optimization},
  volume={34},
  number={1},
  pages={764--789},
  year={2024},
  publisher={SIAM}
}

@article{ying2020note,
  title={A note on optimization formulations of {M}arkov decision processes},
  author={Ying, Lexing and Zhu, Yuhua},
  journal={arXiv preprint arXiv:2012.09417},
  year={2020}
}

@inproceedings{geist2019theory,
  title={A theory of regularized {M}arkov decision processes},
  author={Geist, Matthieu and Scherrer, Bruno and Pietquin, Olivier},
  booktitle={International Conference on Machine Learning},
  pages={2160--2169},
  year={2019}
}

@article{ref:neu-17,
  title={A unified view of entropy-regularized {M}arkov decision processes},
  author={Neu, Gergely and Jonsson, Anders and G{\'o}mez, Vicen{\c{c}}},
  journal={arXiv preprint arXiv:1705.07798},
  year={2017}
}

@inproceedings{ref:Dai-18:boosting,
title={Boosting the Actor with Dual Critic},
author={Bo Dai and Albert Shaw and Niao He and Lihong Li and Le Song},
booktitle={International Conference on Learning Representations},
year={2018},
}

@inproceedings{ref:chen2018:scalable,
  title={Scalable bilinear pi-learning using state and action features},
  author={Chen, Yichen and Li, Lihong and Wang, Mengdi},
  booktitle={International Conference on Machine Learning},
  pages={834--843},
  year={2018}
}

@book{ref:kakade2003:sample,
  title={On the Sample Complexity of Reinforcement Learning},
  author={Kakade, Sham M.},
  year={2003},
  publisher={University of London, University College London}
}

@article{ref:konda1999:actor,
  title={Actor-critic--type learning algorithms for {M}arkov decision processes},
  author={Konda, Vijaymohan R. and Borkar, Vivek S.},
  journal={SIAM Journal on Control and Optimization},
  volume={38},
  number={1},
  pages={94--123},
  year={1999},
  publisher={SIAM}
}

@article{ref:perkins2013:asynchronous,
  title={Asynchronous stochastic approximation with differential inclusions},
  author={Perkins, Steven and Leslie, David S.},
  journal={Stochastic Systems},
  volume={2},
  number={2},
  pages={409--446},
  year={2013},
  publisher={INFORMS}
}

@inproceedings{ref:lee2019:stochastic,
  title={Stochastic primal-dual {Q}-learning algorithm for discounted {M}{D}{P}s},
  author={Lee, Donghwan and He, Niao},
  booktitle={2019 {A}merican {C}ontrol {C}onference},
  pages={4897--4902},
  year={2019}
}

@article{ref:borkar1998:asynchronous,
  title={Asynchronous stochastic approximations},
  author={Borkar, Vivek S.},
  journal={SIAM Journal on Control and Optimization},
  volume={36},
  number={3},
  pages={840--851},
  year={1998},
  publisher={SIAM}
}

@book{filippov2013differential,
  title={Differential equations with discontinuous righthand sides: control systems},
  author={Filippov, Aleksei F.},
  volume={18},
  year={2013},
  publisher={Springer Science \& Business Media}
}

@article{ref:watkins1992:q,
  title={Q-learning},
  author={Watkins, Christopher and Dayan, Peter},
  journal={Machine learning},
  volume={8},
  pages={279--292},
  year={1992},
  publisher={Springer}
}

@article{ref:levine2020:offline,
  title={Offline reinforcement learning: Tutorial, review, and perspectives on open problems},
  author={Levine, Sergey and Kumar, Aviral and Tucker, George and Fu, Justin},
  journal={arXiv preprint arXiv:2005.01643},
  year={2020}
}

@book{ref:arrow1958:studies,
  title={Studies in Linear and Non-Linear Programming},
  author={Arrow, Kenneth J. and Hurwicz, Leonid and Uzawa, Hirofumi and Chenery, Hollis B. and Johnson, Selmer and Karlin, Samuel},
  volume={2},
  year={1958},
  publisher={Stanford University Press Stanford}
}

@inproceedings{ref:dai2018:sbeed,
  title={{SBEED}: Convergent reinforcement learning with nonlinear function approximation},
  author={Dai, Bo and Shaw, Albert and Li, Lihong and Xiao, Lin and He, Niao and Liu, Zhen and Chen, Jianshu and Song, Le},
  booktitle={International Conference on Machine Learning},
  pages={1125--1134},
  year={2018}
}

@article{ref:Manne-60,
 ISSN = {00251909, 15265501},
 author = {Alan S. Manne},
 journal = {Management Science},
 number = {3},
 pages = {259--267},
 publisher = {INFORMS},
 title = {Linear Programming and Sequential Decisions},
 volume = {6},
 year = {1960}
}

@Inbook{ref:Borkar-2002,
author="Borkar, Vivek S.",
title="Convex Analytic Methods in {M}arkov Decision Processes",
bookTitle="Handbook of {M}arkov Decision Processes: Methods and Applications",
year="2002",
publisher="Springer US",
address="Boston, MA",
pages="347--375",
isbn="978-1-4615-0805-2",
doi="10.1007/978-1-4615-0805-2_11",
}

@book{ref:Hernandez-96,
  title={Discrete-Time {M}arkov Control Processes: Basic Optimality Criteria},
  author={Hern{\'a}ndez-Lerma, Onesimo and Lasserre, Jean B.},
  isbn={9780387945798},
  lccn={95037683},
  series={Applications of Mathematics Series},
  year={1996},
  publisher={Springer}
}

@inproceedings{ref:gabbianelli2024:offline,
  title={Offline primal-dual reinforcement learning for linear {M}{D}{P}s},
  author={Gabbianelli, Germano and Neu, Gergely and Papini, Matteo and Okolo, Nneka M.},
  booktitle={International Conference on Artificial Intelligence and Statistics},
  pages={3169--3177},
  year={2024}
}

@article{ref:jaakkola1993:convergence,
  title={Convergence of stochastic iterative dynamic programming algorithms},
  author={Jaakkola, Tommi and Jordan, Michael and Singh, Satinder},
  journal={Advances in Neural Information Processing Systems},
  volume={6},
  year={1993}
}

@inproceedings{ref:haarnoja2018:soft,
  title={Soft actor-critic: Off-policy maximum entropy deep reinforcement learning with a stochastic actor},
  author={Haarnoja, Tuomas and Zhou, Aurick and Abbeel, Pieter and Levine, Sergey},
  booktitle={International Conference on Machine Learning},
  pages={1861--1870},
  year={2018}
}

@inproceedings{ref:schulman2015:trust,
  title={Trust region policy optimization},
  author={Schulman, John and Levine, Sergey and Abbeel, Pieter and Jordan, Michael and Moritz, Philipp},
  booktitle={International Conference on Machine Learning},
  pages={1889--1897},
  year={2015}
}

@article{ref:schulman2017:equivalence,
  title={Equivalence between policy gradients and soft {Q}-learning},
  author={Schulman, John and Chen, Xi and Abbeel, Pieter},
  journal={arXiv preprint arXiv:1704.06440},
  year={2017}
}

@article{ref:tsitsiklis1994:asynchronous,
  title={Asynchronous stochastic approximation and {Q}-learning},
  author={Tsitsiklis, John N.},
  journal={Machine learning},
  volume={16},
  pages={185--202},
  year={1994},
  publisher={Springer}
}

@inproceedings{ref:hong2024:primal,
  title={A Primal-Dual Algorithm for Offline Constrained Reinforcement Learning with Linear {M}{D}{P}s},
  author={Hong, Kihyuk and Tewari, Ambuj},
  booktitle={International Conference on Machine Learning},
  pages={18711--18737},
  year={2024}
}

@inproceedings{ref:zhan2022:offline,
  title={Offline reinforcement learning with realizability and single-policy concentrability},
  author={Zhan, Wenhao and Huang, Baihe and Huang, Audrey and Jiang, Nan and Lee, Jason},
  booktitle={Conference on Learning Theory},
  pages={2730--2775},
  year={2022}
}

@article{ref:chen2024:primal-dual,
    author = {Chen, Yi and Dong, Jing and Wang, Zhaoran},
    title = {A primal-dual approach to constrained {M}arkov Decision Processes},
    journal ={Management Science},
    year = {2024}
}

@article{ref:borkar2005:actor,
  title={An actor-critic algorithm for constrained {M}arkov decision processes},
  author={Borkar, Vivek S.},
  journal={Systems \& control letters},
  volume={54},
  number={3},
  pages={207--213},
  year={2005},
  publisher={Elsevier}
}

@article{ref:li2023:double,
  title={Double duality: Variational primal-dual policy optimization for constrained reinforcement learning},
  author={Li, Zihao and Liu, Boyi and Yang, Zhuoran and Wang, Zhaoran and Wang, Mengdi},
  journal={Journal of Machine Learning Research},
  volume={24},
  number={385},
  pages={1--43},
  year={2023}
}

@article{ref:towers2024gymnasium,
  title={Gymnasium: A standard interface for reinforcement learning environments},
  author={Towers, Mark and Kwiatkowski, Ariel and Terry, Jordan and Balis, John U. and De Cola, Gianluca and Deleu, Tristan and Goulao, Manuel and Kallinteris, Andreas and Krimmel, Markus and KG, Arjun and others},
  journal={arXiv preprint arXiv:2407.17032},
  year={2024}
}

@inproceedings{ref:bas-serrano2021:logistic,
  title={Logistic {Q}-learning},
  author={Bas-Serrano, Joan and Curi, Sebastian and Krause, Andreas and Neu, Gergely},
  booktitle={International conference on artificial intelligence and statistics},
  pages={3610--3618},
  year={2021},
  organization={PMLR}
}

@article{ref:zeng2024:two,
  title={A two-time-scale stochastic optimization framework with applications in control and reinforcement learning},
  author={Zeng, Sihan and Doan, Thinh T. and Romberg, Justin},
  journal={SIAM Journal on Optimization},
  volume={34},
  number={1},
  pages={946--976},
  year={2024},
  publisher={SIAM}
}

@article{ref:weissman2003:inequalities,
  title={Inequalities for the l1 deviation of the empirical distribution},
  author={Weissman, Tsachy and Ordentlich, Erik and Seroussi, Gadiel and Verdu, Sergio and Weinberger, Marcelo J.},
  journal={Hewlett-Packard Labs, Tech. Rep},
  pages={125},
  year={2003}
}

@article{ref:mitrophanov2005:sensitivity,
  title={Sensitivity and convergence of uniformly ergodic {M}arkov chains},
  author={Mitrophanov, Alexander Y.},
  journal={Journal of Applied Probability},
  volume={42},
  number={4},
  pages={1003--1014},
  year={2005},
  publisher={Cambridge University Press}
}

@article{ref:paulin2015:concentration,
  title={Concentration inequalities for {M}arkov chains by {M}arton couplings and spectral methods},
  author={Paulin, Daniel},
  year={2015},
  journal={Electronic Journal of Probability}
}

@book{ref:meyn2012:markov,
  title={Markov Chains and Stochastic Stability},
  author={Meyn, Sean P. and Tweedie, Richard L.},
  year={2012},
  publisher={Springer Science \& Business Media}
}

@article{ref:dobrushin1956:central,
  title={Central limit theorem for nonstationary {M}arkov chains},
  author={Dobrushin, Roland L.},
  journal={Theory of Probability \& Its Applications},
  volume={1},
  number={1},
  pages={65--80},
  year={1956},
  publisher={SIAM}
}

@book{ref:bremaud2020:markov,
  title={Markov Chains: Gibbs fields, Monte Carlo Simulation and Queues},
  author={Br{\'e}maud, Pierre},
  edition={2},
  year={2020},
  publisher={Springer},
  address={Cham},
  isbn={978-3030459819}
}

@phdthesis{ref:BasSerrano2022:LagragianDuality,
  author    = {Bas-Serrano, Joan},
  title     = {Lagragian Duality for Efficient Large-Scale Reinforcement Learning},
  school    = {Universitat Pompeu Fabra},
  type      = {PhD thesis},
  address   = {Barcelona, Spain},
  year      = {2022},
  note      = {Date of defense: 2022-06-28},
  pages     = {126},
}

@article{ref:wang2017:primal,
  title={Primal-Dual $\pi $ Learning: Sample Complexity and Sublinear Run Time for Ergodic Markov Decision Problems},
  author={Wang, Mengdi},
  journal={arXiv preprint arXiv:1710.06100},
  year={2017}
}

\end{document}